\documentclass[11pt,a4paper]{amsart}

\usepackage{amscd,amsthm,amsfonts,amssymb,amsmath,comment,supertabular,mathtools}
\usepackage{mathrsfs}
\usepackage[colorlinks=true]{hyperref}
\usepackage{fullpage}
\usepackage[all]{xy}
\usepackage{xcolor}

    \newcommand{\BA}{{\mathbb {A}}} 
    \newcommand{\BC}{{\mathbb {C}}} 
    \newcommand{\BE}{{\mathbb {E}}} \newcommand{\BF}{{\mathbb {F}}}
    \newcommand{\BG}{{\mathbb {G}}} \newcommand{\BH}{{\mathbb {H}}}

     \newcommand{\BP}{{\mathbb {P}}}
    \newcommand{\BQ}{{\mathbb {Q}}} \newcommand{\BR}{{\mathbb {R}}}

     \newcommand{\BZ}{{\mathbb {Z}}}

    \newcommand{\CA}{{\mathcal {A}}} \newcommand{\CB}{{\mathcal {B}}}
    \newcommand{\CC}{{\mathcal {C}}} 
     
    \newcommand{\CG}{{\mathcal {G}}} 
    \newcommand{\CI}{{\mathcal {I}}} 
     \newcommand{\CL}{{\mathcal {L}}}
    \newcommand{\CO}{{\mathcal {O}}} 
     \newcommand{\CR}{{\mathcal {R}}}
    \newcommand{\CS}{{\mathcal {S}}} 
     
     \newcommand{\CX}{{\mathcal {X}}}

     \newcommand{\RH}{{\mathrm {H}}}

    \newcommand{\fa}{{\mathfrak{a}}} \newcommand{\fb}{{\mathfrak{b}}}

     \newcommand{\fp}{{\mathfrak{p}}}
    \newcommand{\fq}{{\mathfrak{q}}}

     \newcommand{\fP}{{\mathfrak{P}}}

     \newcommand{\fX}{{\mathfrak{X}}}

    \newcommand{\ab}{{\mathrm{ab}}}

    \newcommand{\Aut}{{\mathrm{Aut}}}

    \newcommand{\corank}{{\mathrm{cork}}}
    
    \newcommand{\can}{{\mathrm{can}}}

     \renewcommand{\div}{{\mathrm{div}}}
    \newcommand{\End}{{\mathrm{End}}}

    \newcommand{\Gal}{{\mathrm{Gal}}} \newcommand{\GL}{{\mathrm{GL}}}
    
    \newcommand{\Hom}{{\mathrm{Hom}}}
    
    \renewcommand{\Im}{{\mathrm{Im}}}

    \newcommand{\Ker}{{\mathrm{Ker}}}

    \newcommand{\ord}{{\mathrm{ord}}} 
     \newcommand{\Pic}{\mathrm{Pic}}

    \renewcommand{\mod}{\ \mathrm{mod}\ }

    \newcommand{\Sel}{{\mathrm{Sel}}}

    \newcommand{\Spec}{{\mathrm{Spec\,}}} 
    
    \newcommand{\SU}{{\mathrm{SU}}}

    \newcommand{\tr}{{\mathrm{tr}}}\newcommand{\tor}{{\mathrm{tor}}}

        \newcommand{\Tr}{\mathrm{Tr}}

    \newcommand{\sD}{\mathscr{D}}    
            
    \newcommand{\sU}{\mathscr{U}}    
    \newcommand{\sV}{\mathscr{V}}

        \newcommand{\sP}{\mathscr{P}}

\newcommand{\matrixx}[4]{\begin{pmatrix}
#1 & #2 \\ #3 & #4
\end{pmatrix} }        

\DeclareFontFamily{U}{wncy}{}
\DeclareFontShape{U}{wncy}{m}{n}{<->wncyr10}{}
\DeclareSymbolFont{mcy}{U}{wncy}{m}{n}
\DeclareMathSymbol{\Sha}{\mathord}{mcy}{"58}

    \newcommand{\wt}{\widetilde}
    \newcommand{\wh}{\widehat}

    \newcommand{\ov}{\overline}

    \newcommand{\lra}{\longrightarrow}\newcommand{\lla}{\longleftarrow}
    \newcommand{\ra}{\rightarrow}

\newcommand{\Cor}[1]{}

    \theoremstyle{plain}
    \newtheorem{thm}{Theorem}[section] \newtheorem{coro}[thm]{Corollary}
    \newtheorem{lem}[thm]{Lemma}  \newtheorem{prop}[thm]{Proposition}
    \newtheorem {conj}[thm]{Conjecture} \newtheorem{defn}[thm]{Definition}
    
     \theoremstyle{definition}

\theoremstyle{remark} \newtheorem{remark}{Remark}[section]
\theoremstyle{remark} 
\theoremstyle{remark} 

    \numberwithin{equation}{section}

\newcommand{\RM}{\mathrm{RM}}

\newcommand{\Prob}{\mathrm{Prob}}
\newcommand{\Ort}{\mathrm{Ort}}
\newcommand{\Sym}{\mathrm{Sym}}
\newcommand{\Uni}{\mathrm{Uni}}
\newcommand{\LC}{\mathtt{LC}}
\newcommand{\Sesq}{\mathrm{Sesq}}

\newcommand{\sym}{\mathrm{sym}}

\newcommand{\shom}[2]{\underline{\Hom}(#1,#2)}
\newcommand{\bhom}{{\bf \large Hom}}
 
 \newcommand{\Inv}{\mathrm{Inv}}
 \newcommand{\Ab}{\mathfrak{Ab}}
 \newcommand{\qbinomial}[2]{\left[\begin{matrix}#1\\#2\end{matrix}\right]}
 \newcommand{\type}{\mathrm{Type}}
 \newcommand{\CTP}[2]{\langle #1,#2\rangle_{\mathrm{CT}}}
 \newcommand{\LP}[2]{\langle #1,#2\rangle_{\lambda}}
  \newcommand{\poLP}[2]{\langle #1,#2\rangle_{\lambda,\fp_0}}

\begin{document}
\title[Selmer groups]{Quadratic spaces and Selmer groups of abelian varieties with multiplication}
\author[J. Shu]{Jie Shu}
\address{School of Mathematical Sciences, Tongji University, Shanghai 200092,  P. R. China}
\email{shujie@tongji.edu.cn}

\begin{abstract}
For certain symmetric isogeny $\lambda: A\ra A^\vee$ of abelian varieties over a global field $F$, B. Poonen and E. Rains put an orthogonal quadratic structure on $\RH^1(\BA_F,A[\lambda])$ and realize the Selmer group $\Sel_\lambda(A)$ as an intersection of two maximal isotropic subspaces of $\RH^1(\BA_F,A[\lambda])$. With this understanding of Selmer groups, they expect to model the Selmer groups of elliptic curves and Jacobian varieties of hyperelliptic curves as the intersections of random maximal isotropic subspaces of  orthogonal spaces.  We extend this phenomenon to abelian varieties with multiplication and discuss the Shafarevich-Tate groups. 
\end{abstract}

\subjclass[2020]{Primiary 11G10; Secondary 11E04, 11G05, 14K15, 14K22}

\thanks{}

\maketitle


\section{Introduction}
D. R. Heath-Brown  \cite{Heath-Brown93,HB94}, P. Swinnerton-Dyer \cite{SD08} and D. Kane \cite{Kane13} obtained the distribution for the variable $s(E)=\dim_{\BF_2}\Sel_2(E)/\Im(E(\BQ)_\tor)$ as $E$ varies over quadratic twists of certain elliptic curves over $\BQ$ so that
\[\Prob(s(E)=d)=\left(\prod_{j\geq 0}\frac{1}{1+2^{-j}}\right)\left(\prod_{j=1}^d\frac{2}{2^j-1}\right).\]
B. Poonen and E. Rains \cite{PR12} observed that this distribution coincides with the limit probability distribution of the dimension of the intersection of a fixed maximal isotropic subspace with a random one in an orthogonal quadratic $\BF_p$-space of dimension $2n$ as $n\ra \infty$. The precise distribution they get is, for non-negative integers $d$,
\[\sD_p^\Ort(d)=\left(\prod_{j\geq 0}\frac{1}{1+p^{-j}}\right)\left(\prod_{j=1}^d\frac{p}{p^j-1}\right).\]
Meanwhile they also realized the Selmer group $\Sel_p(E)$ as the intersection of two maximal isotropic subspaces in certain infinite-dimensional orthogonal space. Precisely, using Mumford's theta group and ideas of Zarhin \cite[\S 2]{Zar74}, they equipped the ad\'elic space $\RH^1(\BA_\BQ,E[p])$ as a quadratic $\BF_p$-space. Then arithmetic duality theorems are applied to show that the Kummer images $\prod_{p\leq \infty}E(\BQ_p)/pE(\BQ_p)$ and the global cohomology group $\RH^1(\BQ,E[p])$ are maximal isotropic subspaces of $\RH^1(\BA_\BQ,E[p])$. The Selmer group $\Sel_p(E)$ is the intersection of these two subspaces.  Indeed, they worked this out for symmetric isogenies $\lambda:A\ra A^\vee$ of abelian varieties over global fields. With this understanding of Selmer groups they made the following conjecture.
\begin{conj}[Poonen-Rains]\label{PRC}
Fix a global field $F$ and a prime $p$. As $E$ varies over all elliptic curves over $F$ ordered by height, 
\[\Prob(\dim_{\BF_p}\Sel_pE=d)=\sD_p^\Ort(d).\]
Consequently the average size of $\Sel_p(E)$ for all elliptic curves $E/F$ is $p+1$.
\end{conj}
Analogous conjectures for $n$-Selmer groups of elliptic curves and conjectures for Jacobian varieties of hyperelliptic curves are also made (cf. \cite[Conjecture 1.4, 1.7 and 1.8]{PR12}). Many results along the direction of Conjecture \ref{PRC} are known. Besides the works of Heath-Brown \cite{Heath-Brown93,HB94}, Swinnerton-Dyer \cite{SD08}, D. Kane \cite{Kane13}, G. Yu \cite{Yu05,Yu06}, A. J. de Jong \cite{dJ02}, M. Bhargava and A. Shankar \cite{BS13a,BS13b,BS15a,BS15b} and B. Mazur and K. Rubin \cite{MR10} those are elaborated in \cite{PR12}, we point out some other related works since then:
\begin{itemize}
\item Klagsbrun, Mazur and Rubin \cite{KMR14} set up a Markov process to study the distribution of $2$-Selmer ranks in the family of  quadratic twists of an elliptic curve over an arbitrary number field $K$. Under the hypothesis that $\Gal(K(E[2])/K)\cong S_3$, they described the distribution of the $2$-Selmer ranks in terms of the distribution $\sD_2^\Ort$, adjusted by the disparity factor.  Their method also applies to $p$-Selmer ranks of  a $2$-dimensional self-dual Galois $\BF_p$-representation twisted by characters of order $p$.

\item Generalizing the work \cite{PR12}, Bhargava, Kane, Lenstra, Poonen and Rains modeled the $p^\infty$-Selmer group $\Sel_{p^\infty}(E)$ of elliptic curves over a global field $k$ as intersections of two submodules in $(\BQ_p/\BZ_p)^{2n}$ induced by maximal isotropic subspaces in $\BZ_p^{2n}$ as $n\ra \infty$. This gives a conjectural probability distribution of  the sequence 
\[0\ra E(k)\otimes \BQ_p/\BZ_p\ra \Sel_{p^\infty}(E)\ra \Sha(E)[p^\infty]\ra 0\]
(as a random variable) with values in the set of isomorphism classes of short exact sequences of co-finite $\BZ_p$-modules (cf. \cite[Conjecture 1.3]{BKLPR15}).  This single conjecture would explain many of the known theorems and conjectures on ranks, Selmer groups and Shafarevich-Tate groups of elliptic curves.

\item T. Feng, A. Landesman and Rains \cite{FLR23} proved a version of the BKLPR heuristics for Selmer groups of elliptic curves over $\BF_q(t)$. For the family of all elliptic curves over a global function field, N. Achenjang \cite{Ach23} proved the average size of  $2$-Selmer groups is at most $3+O(1/q)$. J. Ellenberg and A. Landesman \cite{EL25} also proved a version of the BKLPR heuristics for Selmer groups of quadratic twist families of abelian varieties over global function fields.

\item Vastly generalizing the works of Heath-Brown \cite{Heath-Brown93,HB94}, Swinnerton-Dyer \cite{SD08}, D. Kane \cite{Kane13},  A. Smith described the distribution of $2^\infty$-Selmer groups in quadratic twist families of elliptic curves over $\BQ$, and use this to deduce the minimalist conjecture for quadratic twist families \cite{Smith22,Smith25}. J. Pan and Y. Tian also obtained the distribution of $2$-Selmer groups of certain quadratic twist families of elliptic curves over $\BQ$ with an emphasis on their Markov nature \cite{PT25}.  Also Smith and K. Poymans described the distribution of $3$-Selmer groups in cubic twist families of the elliptic curves $y^2=x^3+d$ over $\BQ$ in terms of $\sD_3^\Ort$ \cite{KS24}.
\end{itemize}

For what kind of family of abelian varieties over a global field $F$, the general $n$-Selmer groups would distribute as predicted by the BKLPR heuristics? Besides the Kummer images and global cohomologies should sit properly in the quadratic space, a reasonable feature is that $\End_F(A)=\BZ$ for general abelian varieties in this family.  To explain, let's assume $F=\BQ$ and suppose a commutative ring $\CO$ embeds into $\End_\BQ(A)$ for all abelian varieties in the family. The Kummer images $\prod_{\ell\leq \infty}A(\BQ_\ell)/pA(\BQ_\ell)$ and the global cohomology group $\RH^1(\BQ,A[p])$ are simultaneously maximal isotropic $\BF_p$-subspaces and $\CO$-submodules of $\RH^1(\BA_\BQ,A[p])$. If $\CO$ is strictly larger than $\BZ$, for general $p$, maximal isotropic subspaces attained from Kummer images and global cohomologies would be substantially fewer than the totality of all maximal isotropic subspaces. This may violate the BKLPR heuristics. This reasonable feature is implicit in all the known cases consistent with  the BKLPR heuristics. First a general elliptic curve over global fields has endomorphism ring isomorphic to $\BZ$. As for the conjectures \cite[Conjecture 1.7 and 1.8]{PR12} on hyperelliptic Jacobians, by results of Zarhin \cite{Zar00, Zar01}, if the characteristic of the base field is not $2$, a general hyperelliptic Jacobian is simple with $\BZ$ as its endomorphism ring. In \cite{EL25}, Ellenberg and Landesman requires the abelian variety $A$ has multiplicative reduction with toric part of dimension 1 over some point of the base curve $C$ which also implies $\End_C(A)=\BZ$. 

The present work, generalizing the work of Poonen-Rains \cite{PR12}, addresses this phenomenon where abelian varieties are endowed with multiplication by an order $\CO$ of a number field. Let $\fp$ be an invertible prime ideal of $\CO$ with residue field $k$. Generally, the ad\'elic spaces $\RH^1(\BA_F, A[\fp])$ acquire quadratic $k$-structures of arbitrary types, i.e.  orthogonal, symplectic, unitary and split unitary types.  Moreover, the Selmer group $\Sel_\fp(A)$ is the intersection of two maximal isotropic subspaces, i.e.   the Kummer image and the image of the global cohomology.

\subsection{Quadratic spaces and quadratic distributions}
 Let $k$ be a finite field of cardinality $q=p^e$ with an involution $\sigma$. According to Bak \cite[\S 1]{Bak81}, using the idea of form parameters, we classify all nondegenerate quadratic spaces in the hermitian category of locally compact $k$-spaces. According to different symmetric parameters $\lambda=\pm 1$ and form parameters $\Lambda$, all these nondegenerate quadratic spaces are classified into orthogonal, symplectic and unitary spaces. We refer to Section \ref{quadratic-spaces} for precise definitions and details.

Let $V$ be a nondegenerate quadratic $k$-space of dimension $2n$ which is metabolic, i.e.  contains a maximal isotropic subspace.  The form parameter $\Lambda$ of $V$ is an $\BF_p$-subspace of $k$ and denote $|\Lambda|=q^\delta$. The probability distribution of the rank of the intersection of a fixed maximal isotropic subspace with a random one depends only on the form parameter $\Lambda$ of the quadratic space $V$ and will converges to the limit distribution $\sD_q^\delta$ defined as follows, as $n\ra \infty$. 
 
 \begin{defn}
 Let $q=p^e$. For $\delta\in \left\{0,\frac{1}{2}, 1\right\}$, define the distribution $\sD_q^\delta$ on non-negative integers 
 \[\sD_q^\delta(r)=\left(\prod_{j=0}^{\infty}(1+q^{-\delta-j})^{-1}\right)\left(\prod_{j=1}^r\frac{q^{1-\delta}}{q^j-1}\right).\]
 Here $e$ is even if $\delta=\frac{1}{2}$. 
 \end{defn}
One can also interpret these distributions directly as those arising from the intersections of maximal isotropic subspaces of a metabolic infinite-dimensional quadratic space (cf. Proposition \ref{maximal-isotropic} and \ref{infinite-quadratic}). 
 
 From the dictionary between form parameters and types of quadratic spaces (cf. \S \ref{statistics}, Table \ref{t1}), through the invariant $\delta$, one can readily read off the type of metabolic quadratic spaces where the distribution $\sD_q^\delta$ is induced. To indicate where they come from, we classify these distributions into orthogonal distributions $\sD^\Ort_q$ ($\delta=0$), symplectic distributions $\sD^\Sym_q$ ($\delta=1$) and unitary distributions $\sD^\Uni_q$ ($e$ even and $\delta=\frac{1}{2}$). 
 
 Unlike the orthogonal distributions, non-orthogonal distributions don't have equal parity. They have a bias to taking even values, and even a bias to taking $0$ value. As $q\ra \infty$, while $\sD^\Ort_q$ tends to assign probability $50\%$ to each of $0$ and $1$, a non-orthogonal distribution $\sD_q^\delta, \delta\neq 0$, tends to take $0$ with $100\%$ probability (cf. Proposition \ref{behavior-distribution}).  
 
If one takes $k=k_0\otimes k_0$ to be the split quadratic extension over a finite field $k_0$, with the involution $\sigma(x,y)=(y,x)$. Any nondegenerate hermitian form on a locally compact ``$k$-space"  $V$ identifies $V$ with a metabolic quadratic $k$-space $\BH(V_0)$ where $V_0$ is a $k_0$-space. Such quadratic spaces are called split unitary. Maximal isotropic subspaces in a split unitary space $\BH(V_0)$ are essentially subspaces of $V_0$. Thus they behave very differently from the uniform behavior of maximal isotropic subspaces in quadratic spaces over a finite field and also result in interesting distributions. This is discussed in Section \ref{split-unitary}. 
 
 No wonder almost all these quadratic distributions arise elsewhere as the distribution of (co-)ranks of random matrices of certain type over finite fields \cite{Gerth86, FG15} and as the distribution of ranks of fixed subspaces of a random element in classical groups over finite fields \cite{Fulman97, FS16}.
\subsection{Quadratic structures on cohomologies and Selmer groups}\label{intro-selmer}
Let $F$ be a global field of characteristic $\ell$. Let $A$ be an abelian variety over $F$ with a symmetric isogeny $\lambda:A\ra A^\vee$ where $A^\vee$ denotes the dual abelian variety. The endomorphism algebra $\End^0(A)=\End_F(A)\otimes \BQ$ is a semi-simple algebra. The symmetric isogeny $\lambda$ induces a Rosati involution $\dag$ on $\End^0(A)$: $\phi^\dag=\lambda^{-1}\phi^\vee\lambda$. Let $K$ be a number field contained in $\End^0(A)$ and suppose the order $\CO= K\cap \End_F(A)$ is stable under the Rosati involution $\dag$. Necessarily the Rosati involution stabilizes $K$ and restricts to an involution on $K$. If $\lambda$ is a polarization, the Rosati involution is positive and then $K$ is either a totally real field or a CM field with the Rosati involution restricting to the complex conjugation on $K$. Let $K_0\subset K$ be the subfield fixed by $\dag$ and put $\CO_0=\CO\cap K_0$. 

Let $p$ be a rational prime relatively prime to both the conductor of $\CO$ and $\deg(\lambda)$. Any prime ideal of $\CO$ above $p$ is invertible.  Let $\fp$ be such a prime ideal. Then $\fp_0=\fp\cap \CO_0$ is also an invertible prime ideal of $\CO_0$. We distinguish the following cases:
\begin{itemize}
\item[$(\Ort)$] $\dag|_K=1$, i.e.  $K=K_0$,
\item[($\Sym$)] $\dag|_K\neq1 $ and $\fp$ is ramified in $K/K_0$, 
\item[$(\Uni)$] $\dag|_K\neq1$ and $\fp$ is inert in $K/K_0$, or
\item[$(\SU)$] $\dag|_K\neq1$ and $\fp_0$ is split in $K/K_0$.
\end{itemize}
In first three cases, the Rosati involution $\dag$ stabilizes the prime ideal $\fp$ and induces an involution on the residue field $k=\CO/\fp$. Through duality theorems, the $\CO$-multiplication and the symmetric isogeny $\lambda$ furnish the ad\'elic $k$-space $\RH^1(\BA_F,A[\fp])$ with a metabolic quadratic $k$-structure. The Kummer image $\CL=\prod_{v} A(F_v)/\fp A(F_v)$ is a compact open subspace of $\RH^1(\BA_F,A[\fp])$. For the natural map $\RH^1(F,A[\fp])\ra \RH^1(\BA_F,A[\fp])$, denote its image by $W$ and kernel by $\Sha^1(F,A[\fp])$, then $$W\cong \RH^1(F,A[\fp])/\Sha^1(F,A[\fp]).$$ 
Define the $\fp$-Selmer group
\[\Sel_\fp(A):=\Ker\left(\RH^1(F,A[\fp])\ra \prod_{v}\RH^1(F_v,A)\right).\]

Using duality theorems, one can prove
\begin{thm}\label{main-introduction}
If $p=2$, assume the quadratic conditions as in Theorem \ref{intersection}. 
\begin{itemize}
\item[(1)] $(\RH^1(\BA_F,A[\fp]),h,q)$ is a metabolic orthogonal, symplectic resp. unitary $k$-space in the $(\Ort)$, $(\Sym)$ resp. $(\Uni)$ case accordingly. 
\item[(2)] Both $\CL$ and $W$ are maximal isotropic subspaces in $\RH^1(\BA_F,A[\fp])$ and 
\[\Sel_\fp(A)/\Sha^1(F,A[\fp])\cong \CL\cap W.\]
\end{itemize}
\end{thm}

Such quadratic structures are uniquely determined by the $\CO$-multiplication and the symmetric isogeny $\lambda$ except for the $(\Ort)$ case with $p=2$ whence the quadratic structures also depend on an auxiliary choice of a line bundle inducing the isogeny $\lambda$ (cf. Section \ref{global-cohomology}). For the extreme cases of abelian varieties with $\frac{2\dim A}{[K:\BQ]}=1$ or $2$, the quadratic conditions for $p=2$ are almost automatic (cf. Proposition \ref{global-criterion}, \ref{additional} and Theorem \ref{Rational-condition}, \ref{intersection}) and $\Sha^1(F,A[\fp])=0$ (cf. Remark \ref{vanishing-sha}). As a consequence, we have
\begin{coro}\label{RM+CM}
Assume $\ell\neq p$ and one of the following holds:
\begin{itemize}
\item[(1)] $\dag|_K=1$, $[K:\BQ]=\dim A$ and if $p=2$, assume $2$ is prime to the discriminant of $\CO$;
\item[(2)] $\dag|_K\neq 1$, $[K:\BQ]=2\dim A$ and $\fp$ is stable under the Rosati involution.
\end{itemize}
Then $(\RH^1(\BA_F,A[\fp]),h,q)$ is a metabolic quadratic $k$-space and 
\[\Sel_\fp(A)\cong \CL\cap W.\]
\end{coro}
\begin{remark}
The case $(\SU)$ is discussed in Section \ref{split-unitary} and \ref{split-case}. In the $(\SU)$ case, $\RH^1(\BA_F,A[\fp_0])$ acquires the split unitary $k$-structure and $\Sel_{\fp_0}(A)/\Sha^1(F,A[\fp_0])$ is also an intersection of two maximal isotropic subspaces where $k=\CO/\fp_0$.
\end{remark}

\subsection{Arithmetic families}
In view of Theorem \ref{main-introduction}, if a family of abelian varieties with their ad\'elic cohomological space $\RH^1(\BA_F, A[\fp])$ identified as an infinite-dimensional metabolic quadratic space $V$ and, as $A$ runs over this family, the Kummer image (and the global cohomology) of $A$ sits properly in $V$, one may model the Selmer groups $\Sel_\fp(A)$ as intersections of maximal isotropic subspaces and expect the distribution of Selmer ranks of this family to  be described in terms of quadratic distributions.
\subsubsection{RM abelian varieties}
Abelian varieties with real multiplication are those resemble elliptic curves in many aspects. 
Let $K$ be a totally real field of degree $g$ with $\CO$ its ring of integers and denote $\CO^+$ the set of totally positive elements of $\CO$. 
\begin{defn}
An abelian variety $A$ over a base scheme $S$ of relative dimension $g$  has real multiplication by $\CO$ if there is an embedding of rings $\iota: \CO\ra \End_S(A)$ such that 
\[A\otimes_\CO \Hom_\CO^\sym(A,A^\vee)\cong A^\vee,\quad  (x,\mu)\mapsto \mu(x).\]
\end{defn}
The isomorphism classes of the triple $(A,\iota,\lambda)$, with auxiliary level structures, where $(A,\iota)$ is an abelian variety with RM by $\CO$ and $\lambda$ is an $\CO$-linear principal polarization on $A$, are parametrized by certain Hilbert modular schemes (cf. \cite{Rap78, DP94}). 

Let $F$ be an algebraic number field. Let $(A,\iota)$ be an abelian variety over $F$ with RM by $\CO$ and an $\CO$-linear principal polarization $\lambda$. Fixing the isomorphism class of $(A,\iota)$, the set of the isomorphism class of the triple $(A,\iota,\mu)$, with $\mu$ an $\CO$-linear principal polarization, is an $(\CO^{\times})^{+}/(\CO^\times)^2$-torsor. For each such a torsor, we choose one representative $(A,\iota,\lambda)$ and denote the set of such triples by $\RM_\CO(F)$. That is, $\RM_\CO(F)$ contains exactly one triple $(A,\iota,\lambda)$ for each isomorphism class of abelian varieties over $F$ with RM by $\CO$, if possible. 

Let $(A, \iota,\lambda)$ be an element in $\RM_\CO(F)$. Let $p$ be a rational prime and let $\fp$ be a prime ideal of $\CO$ above $p$ with residue field $k$. If $p\neq 2$, through duality theorems, the real multiplication $\iota$ and the principal polarization $\lambda$  give rise to a unique orthogonal $k$-space $\RH^1(\BA_F,A[\fp])$ with $\Sel_\fp(A)$ as the intersection of two maximal isotropic subspaces. If $p=2$ and the discriminant $d_K$ is odd, with an auxiliary choice of a rational line bundle, the same holds true (see Corollary \ref{RM+CM} and Theorem  \ref{intersection}). 

If $|\RM_\CO(F)|=\infty$, we expect this family to be ``large" and, by Shafarevich's conjecture, proved by Faltings \cite{Faltings86}, we can order elements in $\RM_\CO(F)$ by their conductors. There seem no obvious obstructions to make the following conjecture, generalizing Conjecture \ref{PRC}.
\begin{conj}
Let $K$ be a totally real field with $\CO$ its ring of integers. Let $p$ be a rational prime and $\fp$ a prime ideal of $\CO$ above $p$ with residue field $\BF_q$. Let $F$ be an algebraic number field such that $|\RM_\CO(F)|=\infty$. If $p$ is odd, then, as $(A,\iota, \lambda)$ varies in $\RM_\CO(F)$, ordered by conductors, 
\[\Prob(\Sel_\fp(A)=d)=\sD_q^\Ort(d).\]
Moreover, if $2\nmid d_K$, the same statement holds also for $p=2$.  
\end{conj}

In \cite{KMR14}, Klagsbrun, Mazur and Rubin developed their method in the setting of two dimensional self-dual Galois $\BF_p$-representations. Thus their methods and results, under analogous conditions on the residual Galois representations, apply to the $\fp$-Selmer groups of $p$-th twist families of principally polarized abelian varieties with RM by $\CO$ if  $\CO/\fp=\BF_p$.

\subsubsection{CM abelian varieties}
Let $K$ be a CM field with $\CO_K$ its ring of integers. An abelian variety $A$ over a global field $F$ has complex multiplication by $K$ if there is an embedding $\iota: K\ra \End^0(A)$ and $[K:\BQ]=2\dim A$.  Suppose $\CO=K\cap \End_F(A)$ is stable under the complex conjugation. As in Section \ref{intro-selmer}, given a symmetric isogeny $\lambda$ with the Rosati involution restricting to the complex conjugation on $K$, a prime $p\nmid [\CO_K:\CO]\deg(\lambda)$ and a prime ideal $\fp$ of $\CO$ above $p$, we are in the $(\Sym)$, $(\Uni)$ or $(\SU)$ case according to $\fp_0$ is ramified, inert or split in $\CO$.  Thus if $\fp_0$ is ramified or inert, then $\RH^1(\BA_F, A[\fp])$ is a metabolic symplectic or unitary $k$-space; Otherwise $\RH^1(\BA_F,A[\fp_0])$ is a split unitary $k$-space. In \cite{Shu-CM} and its sequel, we describe the distribution of Selmer ranks in certain families of  $p$-th twists of CM abelian varieties over number fields in terms of the symplectic distribution $\sD_p^\Sym$, the unitary distribution $\sD_{p^2}^\Uni$ and the uniform distribution $\sU_p^t$ (cf. Definition \ref{uniform-distribution0}) respectively.

\subsection{Shafarevich-Tate groups}
Let $K$ be a number field with an order $\CO$ of odd discriminant. Let $A$ be an abelian variety over a global field $F$ of characteristic $\ell\neq 2$. Suppose there is an embedding $K\hookrightarrow \End^0(A)$ and $\CO=K\cap \End^0(A)$. Theorem \ref{Rational-condition} states that, for $A$, being of $\GL_2$-type, i.e. $[K:\BQ]=\dim A$ or the presence of a nontrivial involution $\dag$ on $K$ breaks the obstruction of a $\dag$-sesquilinear symmetric isogeny $\lambda$ being induced from a rational symmetric line bundle.  This ensures the pairing on the Shafarevich-Tate group $\Sha(A)$, which is induced from the Cassels-Tate pairing through $\lambda$, to be alternating (cf. \cite[Theorem 3.3]{Tateduality} and \cite[Corollary 7]{PS99}) and thus in turn has an effect on the linear structure of the Shafarevich-Tate group, especially on the $2$-primary component. 

Let the notations be as in Section \ref{intro-selmer}. There is a symmetric isogeny $\lambda:A\ra A^\vee$ and the order $\CO=K\cap \End^0(A)$ is stable under the Rosati involution $\dag$ induced by $\lambda$. Let $p\nmid [\CO_K:\CO]\deg(\lambda)$ be a prime. Let $\Sha(A)_{/\div}$ denote the quotient of $\Sha(A)$ by its divisible part. 
\begin{thm}
\begin{itemize}
\item[(1)] If $p\neq 2$, then $\Sha(A)_{/\div}[p^\infty]\cong M\bigoplus M$ for some finite $\CO_{0}$-module $M$.
\item[(2)] Assume $p=2$, $\ell\neq 2$ and $2$ is prime to the discriminant of $\CO$. Suppose either $\dag|_K=1$ and $[K:\BQ]=\dim A$ or $\dag|_K\neq 1$. Then $\Sha(A)_{/\div}[2^\infty]\cong M\bigoplus M$ for some finite $\CO_{0}$-module $M$.
\end{itemize}
\end{thm}
The $p\neq 2$ part of the theorem is a direct consequence of a result of M. Flach \cite{Flach90}. As an immediate consequence, we have
\begin{coro}\label{square-sha}
Let $A$ be an abelian variety over a global field $F$ of characteristic $\ell\neq 2$ with an embedding of a number field $K\hookrightarrow \End_F^0(A)$ and $\CO_K=K\cap \End_F(A)$. Suppose there is a symmetric isogeny $\lambda:A\ra A^\vee$ of degree one such that $\CO_K$ is stable under the induced Rosati involution $\dag$. Suppose the discriminant $d_K$ is odd and either $\dag|_K=1$ and $[K:\BQ]=\dim A$ or $\dag|_K\neq 1$.  If $\Sha(A)$ is finite, then $\Sha(A)\cong M\bigoplus M$ for some finite $\CO_{K_0}$-module $M$. 
\end{coro}

The linear structure of the $\fp$-torsion  $\Sha(A)_{/\div}[\fp]$ is also discussed in the final section and examples of abelian varieties with RM/CM and Jacobian varieties of superelliptic curves are provided there.

\subsection{Organizations and conventions}
Section \ref{QS} is devoted to the basic theory of quadratic spaces over finite fields, proves a lifting theorem for adjoint quadratic structures and discusses some combinatorial aspects of quadratic spaces. Section \ref{WPAV} discusses Weil pairings induced from symmetric homomorphisms and their quadratic refinements. There is a strong impact of the multiplication on $A$ to the existence of the quadratic refinements. Section \ref{QC} discusses, focusing on $p=2$, the criterion of local Tate pairings to be even so that local cohomologies acquire quadratic structures, and constructs quadratic maps for the $(\Ort)$ case. In Section \ref{QSS}, we equip the local and ad\'elic cohomologies with quadratic structures and realize the Selmer groups as intersections of maximal isotropic subspaces. In the final section, we discuss the structure of Shafarevich-Tate groups under the impact of the multiplication.

Throughout the sections \ref{WPAV}-\ref{QSS}, we will use sheaves of abelian groups over the big fppf site over $F$ and fppf cohomologies. Usually commutative group schemes will be viewed as sheaves.  Points on schemes are generally scheme-valued points. The flat cohomology of a sheaf represented by a smooth commutative group scheme can be interpreted by the \'etale cohomology, or even the Galois cohomology (cf. \cite[\S 11]{Groth68}).

\subsection*{Acknowledgements} 
The debt to the works \cite{PR11,PR12} of B. Poonen and E. Rains is evident. We thank Ye Tian and Lishuang Lin for helpful communications.

\section{Quadratic spaces and maximal isotropic subspaces}\label{QS}
\subsection{Hermitian category of locally compact $k$-spaces}\label{hermitian-spaces}
Let $k$ be a finite field of characteristic $p$ and cardinality $q$. Let $\sigma$ be an involution on $k$ having $k_0$ as its fixed field. If $\sigma=1$, then $k=k_0$. If $\sigma\neq 1$, $k$ is a quadratic extension of $k_0$ and $\sigma$ is the nontrivial element in the Galois group $\Gal(k/k_0)$. 

We endow the finite field $k$ with the discrete topology. A topological $k$-space is a Hausdorff abelian group $V$ together with a continuous $k$-linear structure on $V$ (cf. \cite[Chap. I-\S 1]{Bourbaki-EVT}). A finite dimensional topological $k$-space is always discrete. A topological $k$-space is called linearly topologized if it has a local topological base at $0$ consisting of open $k$-linear subspaces. A theorem of Dantzig states that any locally compact $k$-space is linearly topologized (cf. Lemma \ref{TDLC}). 

For topological $k$-spaces $V$ and $W$, denote by $\Hom_k(V,W)$ the group of continuous $k$-linear maps which is endowed with the compact-open topology and a $k$-linear structure by setting  $(af)(v)=f(\sigma(a)v)$. Denote the dual space by $V^*=\Hom_k(V,k)$.

Let $k'$ be a subfield of $k$. Any topological $k$-space $V$ can be viewed naturally as a $k'$-space. For any topological $k'$-space $W$,  $\Hom_{k'}(V,k')$ also admits a $k$-linear structure by $(af)(v)=f(\sigma(a)v)$. For a topological $k'$-space $W$, we endow $k\otimes_{k'}W\cong k'^{\oplus^{[k:k']}}\otimes_{k'}W$ with the product topology which is independent of the identification $k\cong k'^{\oplus^{[k:k']}}$ we choose. There is a continuous map $$\Tr_{k/k'}\otimes_{k'} 1_W: k\otimes_{k'} W\ra W, \quad a\otimes w\mapsto \Tr_{k/k'}(a)w.$$ It is known that $V^*$ is isomorphic to the Pontryagin dual $V'=\Hom_{\BF_p}(V,\BF_p)$ (cf. \cite[Theorem p.332]{Flood79}) and the following proposition gives an explicit construction of this isomorphism. 
\begin{prop}\label{duality}
Let $V$ be a topological $k$-space and $W$ a linearly topologized $k'$-space.
\begin{itemize}
\item[(1)] Composing $\Tr_{k/k'}\otimes_{k'} 1_W$ induces an isomorphism $\varphi_{V,W}: \Hom_k(V, k\otimes_{k'} W) \ra \Hom_{k'}(V,W)$ of topological $k$-spaces.
\item[(2)] If $V$ is locally compact, then the canonical map $\can_V: V\ra {V^*}^*$ defined by $\can_V(v)(f)=f(v)$ for $v\in V,f\in V^*$, is an isomorphism of locally compact $k$-spaces. 
\end{itemize}
\end{prop}

\begin{proof}
\noindent (1) The map $\varphi_{V,W}$ is clearly $k$-linear. The injectivity follows from the fact that the trace pairing $(x,y)=\Tr_{k/k'}(xy)$ is perfect and $k$ is free over $k'$. Next we prove $\varphi_{V,W}$ is surjective following \cite[Lemma A.3]{MR07}. Fix an $k'$-basis $\{\alpha_1,\cdots,\alpha_n\}$ of $k$ and let $\{\alpha_1^*,\cdots,\alpha_n^*\}$ be the dual basis with respect to the trace pairing, i.e.  $\Tr_{k/k'}(\alpha_i\alpha_j^*)=\delta_{ij}$. For $f\in \Hom_{k'}(V,W)$, define
\begin{equation}\label{inverse}
\wh{f}(v)=\sum_{i=1}^n\alpha_i^*\otimes f(\alpha_i v).
\end{equation}
One can directly verify $\wh{f}\in \Hom_k(V,k\otimes_{k'}W)$ and $\varphi_{V,W}(\wh{f})=f$.

Next we prove $\varphi_{V,W}$ is a homeomorphism. Fix an identification $k\otimes_{k'} W\cong W^{\oplus^n}$ where $n=[k:k']$. For any compact subset $K\subset V$ and any open $k'$-linear subspace $U\subset W$, define 
\[W_k(K,U^{\oplus^n})=\{f\in \Hom_k(V,k\otimes_{k'}W): f(K)\subset U^{\oplus^n}\}\]
and $W_{k'}(K,U)\subset \Hom_{k'}(V,W)$ is similarly defined. As $K$ runs through all compact subsets of $V$ and $U$ runs through all open $k'$-linear subspaces of $W$, $W_k(K,U^{\oplus^n})$ resp. $W_{k'}(K,U)$ constitutes a local base of the zero map in $\Hom_k(V,k\otimes_{k'}W)$ resp. $\Hom_{k'}(V,W)$. Clearly $\varphi_{V,W}(W_k(K,U^{\oplus^n}))\subset W_{k'}(K,U)$ and thus $\varphi_{V,W}$ is continuous. It follows from (\ref{inverse}) that $$\varphi_{V,W}^{-1}\left(W_{k'}(\cup_{i=1}^n\alpha_i K,U)\right)\subset W_{k}(K,U^{\oplus^n})$$ and thus $\varphi_{V,W}^{-1}$ is continuous. So $\varphi_{V,W}$ is a homeomorphism.

\noindent(2) Taking $W=k'=\BF_p$, the second statement follows directly from the Pontryagin duality for locally compact abelian groups. 
\end{proof}

Denote by $\LC(k)$ the category of locally compact $k$-vector spaces. For $V,W\in \LC(k)$, $V\oplus W$ endowed with the product topology defines a direct sum in $\LC(k)$. It is not hard to see $\LC(k)$ is an additive category. Every object $V\in \LC(k)$ is associated with a dual space $V^*$ and every morphism $f:W\ra V$ is associated with a dual morphism $f^*:V^*\ra W^*$ by $f^*(v^*)=v^*\circ f$. The induced map 
\[\Hom_k(W,V)\lra \Hom_k(V^*,W^*),\quad f\mapsto f^*\]
is $\sigma$-sesqlinear. Thus we obtain a contravariant $\sigma$-sesqlinear functor $*:\LC(k)\ra\LC(k)$. 
\begin{coro}
$(\LC(k), *)$ is a hermitian catergory in the sense of \cite[\S 7.2]{Sch85}.
\end{coro}
\begin{proof}
The duality functor $*$ respects direct sums and, by Proposition \ref{duality}.(3), the canonical maps $\can_V$ gives a natural equivalence between the functor $**$ and the identity functor of $\LC(k)$.
\end{proof}
As a result, we can freely use results from Chapter 7 of \cite{Sch85}.
\begin{defn}
Let $V$ be a locally compact $k$-space. 
\begin{itemize}
\item[(1)] A $\sigma$-sesquilinear form $h$ on $V$ is a continuous map $h: V\times V\ra k$ satisfying 
\begin{eqnarray*}
&h( x+y,z)=h( x,z)+h( y,z), &h( x,y+z)=h( x,y)+h(x,z),\\
 &h(ax,y)=ah(x,y), &h(x,ay)=\sigma(a)h(x,y),
 \end{eqnarray*}
 for all $x,y\in V$ and $a\in k$.

\item[(2)] The transpose $h^*$ of a sesquilinear form $h$ is defined by 
\[h^*(x,y)=\sigma(h(y,x)).\]
For $\lambda=\pm 1$, a sesquilinear form $h$ is called $\lambda$-hermitian if $h=\lambda h^*$.  If $s$ is a sesquilinear form, then the form $h=s+\lambda s^*$ is $\lambda$-hermitian and such forms are called even $\lambda$-hermitian. A $\lambda$-hermitian space $(V,h)$ is a locally compact $k$-space $V$ endowed with a $\lambda$-hermitian form. Moreover,  a hermitian space $(V,h)$ is called even if $h$ is even.
\item[(3)] A sesquilinear form $h$ is called nondegenerate if the homomorphism $V\ra V^*=\Hom_k(V,k)$ defined by $y\mapsto (x\mapsto h(x,y))$ is an isomorphism. A hermitian space $(V,h)$ is called nondegenerate if $h$ is nondegenerate.
\end{itemize}
\end{defn}

If $\sigma=1$, a $\sigma$-sesquilinear form is also called bilinear. A bilinear form $h$ satisfying $h(x,x)=0$ for all $x\in V$ is called alternating. By a subspace of $V$, we mean a subspace in the category $\LC(k)$, or equivalently, a closed $k$-subspace (cf. \cite[Proposition I-6]{RV98}).

\begin{lem}\label{TDLC}
Any $V\in \LC(k)$ has a local topological base at $0$ consisting of compact open $k$-subspaces.
\end{lem}
\begin{proof}
As a locally compact abelian group of exponent $p$, $V$ is totally disconnected. Then Dantzig's theorem \cite[p. 411]{Dantzig36} states there is a local base at $0$ consisting of compact open subgroups. For any $W$ in this local base, $\bigcap_{a\in k^\times}aW$ is a compact open $k$-subspace contained in $W$.
\end{proof}
\begin{prop}\label{even-alternating}
Suppose $\sigma=1$ and $h$ is a bilinear form on $V$. Then $h$ is an even $(-1)$-hermitian form if and only if $h$ is alternating.
\end{prop}
\begin{proof}
If $h$ is an even $(-1)$-hermitian form, it is clear from the definition that $h(x,x)=0$. Conversely suppose $h$ is alternating, necessarily $h$ is $(-1)$-hermitian. In the following we show $h$ is even. 

First we reduce to the case where either $V$ is compact or discrete. By the previous lemma, there exists a compact open subspace $W\subset V$ and therefore $V/W$ is discrete.  Let $W/V\ra V$ be an arbitrary (automatically continuous) $k$-linear section to the quotient map $V\ra V/W$ whose existence is guaranteed by Zorn's lemma.
It is not hard to see the induced map $W\oplus V/W\ra V$ is an isomorphism of locally compact $k$-spaces. Now we may assume $V=W\oplus U$ with $W$ a compact $k$-space and $U$ a discrete $k$-space. Then both $h_W=h|_{W\times W}$ and $h_U=h|_{U\times U}$ are alternating. If there exist continuous bilinear forms $s_W$ and $s_U$ on $W$ and $U$ respectively such that $h_W=s_W-s_W^*$ and $h_U=s_U-s_U^*$, then 
\[s=\left(\begin{matrix}s_W&h|_{W\times U}\\0&s_U\end{matrix}\right)\]
is a continuous bilinear form on $V$ and $h=s-s^*$.

Assume $V$ is discrete and we are not bothered with topologies. Let $\CS$ be the set of pairs $(W,s_W)$ where $W\subset V$ is a subspace and $s_W$ is a  bilinear form on $W$ such that $h_W=s_W-s_W^*$. By considering finite dimensional subspaces, we see $\CS$ is non-empty (see \cite[I-(2.4)]{Knus91}). The set $\CS$ is equipped with a partial order: $(W,s_W)<(W',s_{W'})$ if $W\subset W'$ and $s_{W'}$ restricts to $s_W$. For any totally ordered chain $\{(W_i,s_i): i\in I\}$ in $\CS$, $s_i$ extends to a bilinear form $s_W$ on $W=\bigcup_{i\in I} W_i$ and $h|_W=s_W-s_W^*$.  Thus each totally ordered chain has an upper bound in $\CS$. By Zorn's lemma, there exists a maximal pair, say $(W,s_W)$, in $\CS$. We claim $W=V$. If not, let $v\notin W$, form the direct sum $U=W\bigoplus kv$ and define $s_U$ by 
\[s_U(w+u,w'+u')=s_W(w,w')+h(w,u'),\quad w,w'\in W,u,u'\in kv.\]
It is direct to check $s_U$ is a bilinear form on $U$ and $h_U=s_U-s_U^*$, because $h$ is alternating. Then $(U,s_U)$ belongs to $\CS$, which conflicts with the maximality of $(W,s_W)$. Thus $W=V$ and $h=s_W-s_W^*$ is even.

Finally we assume $V$ is compact. By the continuity of $h$ and compactness of $V$, there exists an open neighborhood $W$ of $0$ such that $h(W,V)=h(V,W)=0$. By Lemma \ref{TDLC}, we may assume $W$ is a compact open subspace. Then $h$ descends to an alternating form on the discrete $k$-space $V/W$ and we are done.

\end{proof}

\subsection{Quadratic $k$-spaces}\label{quadratic-spaces}
Fix $\lambda=\pm 1$. According to Bak \cite[\S 1]{Bak81}, a form parameter in $k$ is an additive subgroup $\Lambda$ of $k$ with $\Lambda^{\min}\subset \Lambda\subset \Lambda^{\max}$, where
\[\Lambda^{\min}=\{a-\lambda \sigma(a): a\in k\},\quad \Lambda^{\max}=\{a\in k: a=-\lambda \sigma(a)\},\]
and $a\sigma(a)x\in \Lambda$ for all $x\in \Lambda$ and $a\in k$. For each $V\in \LC(k)$, let $\Sesq(V)=\Hom_k(V,V^*)$ be the group of sesquilinear forms on $V$, and define
\[\Lambda(V)=\{s\in \Sesq(V):s=-\lambda s^*, s(v,v)\in \Lambda\}.\]
For a class $[s]\in \Sesq(V)/\Lambda(V)$, one can associates with $[s]$ a well-defined even $\lambda$-hermitian form $h=s+\lambda s^*$ and a map, called {\em the quadratic map} associated to $[s]$,
\[q:V\ra k/\Lambda,\quad q(v)=s(v,v)\mod \Lambda.\] 
The class $[s]$ is uniquely determined by the associated even hermitian form $h$ and the quadratic form $q$. 
\begin{defn}[Cf. {\cite[\S 1]{Bak81}, \cite[7-3.3, 7-3.5]{Sch85}}]
Fix $\lambda=\pm 1$ and the form parameter $\Lambda$. For $V\in \LC(k)$, $[s]\in \Sesq(V)/\Lambda(V)$, let $h$ and $q$ be the even $\lambda$-hermitian form and quadratic map associated to $[s]$ respectively.  The triple $(V, h,q)$ is called the $(\lambda,\Lambda)$-quadratic space  over $k$ associated to the class $[s]$.
\end{defn}

\begin{remark}\label{quadratic-evenhermitian2}
Let $(V,h)$ be an even $\lambda$-hermitian space. Any choice of a sesquilinear form $s$ satisfying $h=s+\lambda s^*$ and a form parameter $\Lambda$ will make $(V,h)$ into a $(\lambda,\Lambda)$-quadratic space with quadratic map $q(v)=s(v,v)\mod \Lambda$. Moreover, if $\Lambda=\Lambda^{\max}$, the class $[s]$ of any $s$ satisfying $h=s+\lambda s^*$, and hence the quadratic map $q$, are uniquely determined by $h$. 
\end{remark}

\begin{lem}\label{auto-quadratic}
Suppose either $p\neq 2$ or $\sigma\neq 1$. The following hold:
\begin{itemize}
\item[(1)] $\Lambda^{\min}=\Lambda^{\max}$.
\item[(2)] Any $\lambda$-hermitian space $(V,h)$ is even and determines a unique quadratic space $(V,h,q)$.
\end{itemize}
\end{lem}
\begin{proof}
This is \cite[Lemma 7-6.5]{Sch85}.
\end{proof}

\begin{defn}
If $\sigma=1$, we call a nondegenerate $(1,\Lambda^{\min})$-quadratic space an orthogonal space, and a nondegenerate $(-1,\Lambda^{\max})$-quadratic space a symplectic space. If $\sigma\neq 1$, a nondegenerate quadratic spaces is called a unitary space. 
\end{defn}

\begin{remark}
If $p=2$ and $\sigma=1$, for $\lambda=\pm 1$, by definition, we see $\Lambda^{\min}=0$, $\Lambda^{\max}=k$ and a form parameter is either $\Lambda^{\min}$ or $\Lambda^{\max}$. Thus, together with Lemma \ref{auto-quadratic}, we conclude that orthogonal, symplectic and unitary spaces exhaust all nondegenerate quadratic spaces.
\end{remark}

\begin{remark}\label{classical-orthogonal}
Suppose $\sigma=1$ and $\lambda=1$. Let $q:V\ra k$ be a (continuous) quadratic form in the usual sense, i.e.  $h(w,v)=q(w+v)-q(w)-q(v)$ is $k$-bilinear and $q(av)=a^2q(v)$ for any $a\in k$. One can show there exists a bilinear form $s$ such that 
\[q=s(v,v)\text{ and } h=s+s^*.\]
Thus orthogonal spaces are exactly nondegenerate quadratic spaces in the usual sense. The situation of finite dimensional spaces are treated in \cite[Lemma 7-3.2]{Sch85}). The general case can be proved using the construction technique in  {\em loc.\ cit.} and Zorn's lemma in the same manner as in the proof of Proposition \ref{even-alternating}.
\end{remark}

Let $(V,h,q)$ be a quadratic $k$-space. A vector $v\in V$ is called isotropic if $q(v)=0$ (necessarily $h(v,v)=0$). If $X\subset V$, we denote by $X^\perp$ the orthogonal complement of $X$ under the hermitian form $h$. 
\begin{defn}
Let $(V,h,q)$ be a quadratic $k$-space.
\begin{itemize}
\item[(1)] A  subspace $X\subset V$ is called isotropic if $q(X)=0$ and $X\subset X^\perp$, or equivalently $s|_{X\times X}\in \Lambda(X)$.  An isotropic subspace $X$ is called maximal if $X=X^\perp$.
\item[(2)] The quadratic space $V$ is called metabolic if $h$ is nondegenerate and $V$ contains a maximal isotropic subspace as a direct summand in the category $\LC(k)$.
\end{itemize}
\end{defn}

\begin{lem}\label{orth-iso}
Let $V$ be a nondegenerate quadratic $k$-space. If $p\neq 2$ or $V$ is symplectic or unitary, then a subspace $X\subset V$ is maximal isotropic if and only if $X=X^\perp$.
\end{lem}
\begin{proof}
Suppose the quadratic space $V$ is associated to $[s]\in \Sesq(V)/\Lambda(V)$. For $x\in V$, $h(x,x)=0$ is equivalent to $s(x,x)\in \Lambda^{\max}$ while $q(x)=0$ is equivalent to $s(x,x)\in \Lambda$. By Lemma \ref{auto-quadratic} and the assumption for the quadratic spaces,  the form parameter $\Lambda=\Lambda^{\max}$ and the lemma follows.
\end{proof}

For any locally compact $k$-space $V$, the hyperbolic space 
\[\BH(V)=\left(V\oplus V^*,\left[\left(\begin{matrix}0&1_{V^*}\\0&0\end{matrix}\right)\right]\right)\]
is metabolic with $V$ a maximal isotropic subspace. 

\begin{prop}\label{metabolic}
Let $V$ be a nondegenerate quadratic $k$-space. 
\begin{itemize} 
\item[(1)] If $V$ contains a maximal isotropic subspace, then it contains a compact open maximal isotropic subspace.
\item[(2)] If $V$ contains a maximal isotropic subspace, then $V$ is metabolic.
\item[(3)] If $V$ is metabolic and $X$ is an open maximal isotropic subspace, then $V$ is isomorphic to the hyperbolic space $\BH(X)$. 
\item[(4)] Open maximal isotropic subspaces are compact.
\end{itemize}
\end{prop}
\begin{proof}
(1) The case of orthogonal spaces is explained in \cite[Remark 2.15]{PR12} and the case of other quadratic spaces can be discussed analogously. 

\noindent (2)-(3) If $V$ contains a maximal isotropic subspace, then it contains a (compact) open maximal isotropic subspace. On the other hand, for any open maximal isotropic subspace $X$, $V\cong X\oplus V/X$ as locally compact $k$-spaces.  By a result of Knebusch \cite[Lemma 7-3.7]{Sch85}, $V\cong \BH(X)$ as quadratic spaces. 

\noindent (4) If $X$ is an open maximal isotropic subspaces, then $V/X$ is discrete and $X=X^\perp\cong (V/X)^*$ is compact by Pontryagin duality.

\end{proof}

\subsection{Adjoint quadratic structures}
\begin{defn}
\begin{itemize}  
\item[(1)] Let $V$ be a  topological $k$-space and $W$ a topological $\BF_p$-space. A continuous $\BF_p$-bilinear map $h: V\times V\ra W$ is called  $\sigma$-adjoint  $($or simply adjoint if $\sigma=1$$)$ if $h(ax,y)=h(x,\sigma(a)y)$ for all $a\in k$ and $x,y\in V$. 
\item[(2)] Suppose $\sigma=1$ and $k'\subset k$ is a subfield. For topological $k$-space $V$ and $k'$-space $W$, a continuous map $q:V\ra W$ is called  $k'$-quadratic if $h(x,y)=q(x+y)-q(x)-q(y)$ is $k'$-bilinear and $q(ax)=a^2q(x)$ for all $a\in k'$ and $x\in V$. Moreover, a $k'$-quadratic map $q$ is adjoint if $h$ is adjoint. 
\end{itemize}
\end{defn}

Our aim in this subsection is to prove
\begin{thm}\label{lifting-thm}
Let $V\in \LC(k)$ and  $(V,h,q)$ be a nondegenerate quadratic $\BF_p$-space with $h$ $\sigma$-adjoint.
\begin{itemize}
\item[(1)] Suppose $\sigma\neq 1$. There is a unique unitary $k$-space $(V,h',q')$ such that $\Tr_{k/\BF_p}h'=h$. Moreover, if $p\neq 2$ or $(V,h,q)$ is symplectic, then $\Tr_{k/\BF_p}q'=q$.

\item[(2)] If $\sigma=1$, then there is a unique nondegenerate quadratic $k$-space $(V,h',q')$ of the same type as $(V,h,q)$ such that $\Tr_{k/\BF_p}h'=h$ and $\Tr_{k/\BF_p}q'=q$.
\item[(3)] If $p=2$, $\sigma\neq 1$ and $(V,h,q)$ is orthogonal, assume additionally that $\Tr_{k/\BF_2}q'=q$. Let $X$ be a $k$-subspace of $V$. Then $X$ is maximal isotropic for $(V,h,q)$ if and only if  it is maximal isotropic for $(V,h',q')$. 
\end{itemize}
\end{thm}

\begin{proof}
(1)-(2) The only involution on $\BF_p$ is the identity map, and therefore $(V,h,q)$ is either orthogonal or symplectic. If $h'$ is  the unique $\sigma$-sesquilinear form satisfying $\Tr_{k/\BF_p}h'=h$, guaranteed by Proposition \ref{adjoint-form}, and $h$ is $\lambda$-hermitian, then $h'$  is  nondegenerate and $\lambda$-hermitian because the bijection of Proposition \ref{adjoint-form} preserves nondegeneracy and transpose.

If $p\neq 2$ or $\sigma=1$, there exists a $\sigma$-adjoint $\BF_p$-bilinear form $s$ such that $(V,h,q)$ is the quadratic space associated to the class of $s$ with corresponding form parameter. Indeed, if $p\neq 2$, just let $s=h/2$; If $p=2$ and $\sigma=1$, forgetting the form parameter, $q$ is an adjoint $\BF_2$-quadratic form. The existence of the required $s$ is guaranteed by Corollary \ref{adjoint-quadratic}.(2). By Proposition \ref{adjoint-form}, there exists a unique $\sigma$-sesquilinear form $s'$ such that $\Tr_{k/\BF_p}s'=s$. The class of $s'$ with the indicated form parameter is uniquely determined and let $(V,h',q')$ be the corresponding quadratic $k$-space. It follows that $h'$ is the unique hermitian form satisfying $\Tr_{k/\BF_p}h'=h$ and $\Tr_{k/\BF_p}q'=q$. 

If $p=2$ and $\sigma\neq 1$,  let $h'$ be the unique hermitian form, guaranteed by Proposition \ref{adjoint-form}, satisfying $\Tr_{k/\BF_p}h'=h$. Then $h'$ defines a unique unitary $k$-structure $(V,h',q')$ by Lemma \ref{auto-quadratic}. Moreover, if $(V,h,q)$ is symplectic, $\Tr_{k/\BF_p}q'=q=0$.

\noindent (3) First we show $X^{\perp,h'}=X^{\perp,h}$. Clearly we have the inclusion $X^{\perp,h'}\subset X^{\perp,h}$. For any $x\in X$, $y\in X^{\perp,h}$ and $a\in k$,
\[\Tr_{k/\BF_p}ah'(x,y)=\Tr_{k/\BF_p}h'(ax,y)=h(ax,y)=0.\]
By the nondegeneracy of the trace form, we have $h'(x,y)=0$. Thus $X^{\perp,h}\subset X^{\perp,h'}$ and hence the equality holds.  Thus we can write $X^\perp$ without referring to the hermitian forms. 

Suppose $p\neq 2$ or $(V,h,q)$ is symplectic. If $(V,h,q)$ is symplectic, then  $(V,h',q')$ is either symplectic or unitary. By Lemma \ref{orth-iso}, $X$ is maximal isotropic for $(V,h,q)$ resp. $(V,h',q')$ if and only if $X=X^{\perp}$. 

Suppose $p=2$ and $(V,h,q)$ is orthogonal. It suffices to prove that $q|_X=0$ if and only if $q'|_X=0$. By assumption, in any case, we have $\Tr_{k/\BF_2}q'=q$. Clearly $q'|_X=0$ implies $q|_X=0$. Conversely, suppose $q|_X=0$. For $b\in k$ and any $x\in X$,
\begin{equation}\label{trace-form}
\Tr_{k/\BF_p}( bb^\sigma q'(x))=\Tr_{k/\BF_p}q'(bx)=q(bx)=0.
\end{equation}
Let $\Lambda$ and $\Lambda'$ be the from parameter of $(V,h,q)$ and $(V,h',q')$ respectively. If $\sigma=1$, i.e.  $(V,h',q')$  is orthogonal, then $\Lambda =\Lambda'=0$. As $b$ runs through $k$, $b^2=bb^\sigma$ runs through $k$. By the nondegeneracy of the trace form and (\ref{trace-form}), we have $q'(x)=0$. 

If $\sigma\neq 1$, i.e.  $(V,h',q')$  is unitary, then $\Lambda'=\{a-\sigma(a):a\in k\}$. As $b$ runs through $k$, $bb^\sigma$ runs through $k_0$. The annihilator of $k_0$ under the trace form is exactly $\Lambda'$. By (\ref{trace-form}),  $q'(x)=0 \mod \Lambda'$ as desired.
\end{proof}
\begin{remark}
For Assertion (1), the unitary $k$-structure on $V$ is determined by the hermitian form $h'$, and hence by $h$. But if $p=2$ and $(V,h,q)$ is orthogonal, different quadratic forms $q$ may give rise to the same bilinear form $h$. So in this case, $q$ may not be determined by $h'$. For example of finite dimensional spaces, if $(V,h,q)$ is associated to the class of $s$ and $s$ is represented by a matrix, then arbitrarily modifying the diagonal entries of $s$ gives different quadratic forms $q$ resulting the same bilinear form $h$. This phenomenon also accounts for the additional assumption of Assertion (3).
\end{remark}

Let's proceed to develop results that we need to prove Theorem \ref{lifting-thm}. First we introduce the tensor product of topological $k$-spaces. For linearly topologized $k$-spaces $V$ and $W$, define the tensor product of $V$ and $W$ to be $V\otimes_k W$ endowed with the strongest topology such that $V\otimes_k W$ is a linearly topologized $k$-space and the bilinear map $V\times W\ra V\otimes W$ is continuous. We list the following basic facts whose proof we omit.
\begin{itemize}
\item[(1)]  $V\otimes_k W$ has the following universal property: for any (abstract) bilinear map $h: V\times W\ra Z$ where $Z$ is a linearly topologized $k$-space, the induced linear map $\overline{h}: V\otimes_kW\ra Z$ is continuous if and only if $h$ is continuous. 
\item[(2)] A topological local  base at $0$ of $V\otimes_k W$ can be given by the following subspaces
\[\sum_{i\in I}X_i\otimes Y_i\]
where $X_i$ and $Y_i$ are any indexed open subspaces of $V$ and $W$ respectively satisfying $\bigcup_{i\in I}X_i=V$ and $\bigcup_{i\in I}Y_i=W$.  Here the sum means finite sums. 
\item[(3)] The tensor product topology on $k^{\oplus^n}\otimes_k V\cong V^{\oplus^n}$ is equivalent to the product topology. In particular, $k\otimes_k V\cong V$ as $k$-spaces.
\end{itemize}

\begin{prop}\label{adjoint-form}
Let $V$ be a topological $k$-space and $W$ a linearly topologized $\BF_p$-space. Composing $\Tr_{k/\BF_p}\otimes_{\BF_p} 1_W$ gives a bijection from the set of continuous $\sigma$-sesquilinear maps $V\times V\ra  k\otimes_{\BF_p}W$ to the set of continuous $\sigma$-adjoint $\BF_p$-bilinear maps $V\times V\ra W$. Moreover, if $W=\BF_p$, this bijection preserves nondegeneracy and transpose.  
\end{prop}
\begin{proof}
It follows from Proposition \ref{duality}.(2) that composing $\Tr_{k/\BF_p}\otimes_{\BF_p} 1_W$ induces an isomorphism
\[\Hom_k(V,\Hom_k(V,k\otimes_{\BF_p}W))\cong \Hom_k(V,\Hom_{\BF_p}(V,W)).\]
The left-hand side is the set of $\sigma$-sesquilinear maps $V\times V\ra k\otimes_{\BF_p}W$ while the right-hand side is the set of $\sigma$-adjoint maps $V\times V\ra W$.
\end{proof}

Thanks to Lemma \ref{auto-quadratic}, the above proposition (with $W=\BF_p$) enable us to prove Theorem \ref{lifting-thm} when $p\neq 2$ or $\sigma\neq 1$. In the following, we assume $p=2$ and $\sigma=1$. Let $V$ be a linearly topologized $k$-space. Inspired by the construction of \cite[\S 2.1]{PR11}, let $\mathtt{U}V=V\oplus (V\otimes_k V)$ be the topological group with the identity element $0$ and the multiplication as
\[(m+t)(m'+t')=(m+m')+((m\otimes m')+t+t'),\quad m,m'\in V,t,t'\in V\otimes_kV.\]
There is a central extension of topological groups
\[0\ra V\otimes_k V\ra \mathtt{U}V\ra V\ra 0.\]
The surjection $\mathtt{U}V\ra V$ admits a continuous section $\iota:V\ra \mathtt{U}V$ sending $m\mapsto m$. If $m,m'\in V$, then
\begin{equation}\label{multi}
\iota(m)\iota(m')\iota(m+m')^{-1}=m\otimes m' \in V\otimes V \subset \mathtt{U}V.
\end{equation}

The closed subgroup $H$ generated by all commutators $[m+t,m'+t']=(m\otimes m'-m'\otimes m)$ and all squares $(m+t)^2=(m\otimes m)$ is a normal subgroup of $\mathtt{U}V$.  The Frobenius map $x\mapsto x^2$ is an automorphism of the field $k$ and denote its inverse by $x\mapsto x^{1/2}$. It is {\em fortunate} that the quotient group  $\mathtt{U}V/H$ admits a topological $k$-linear structure through the scalar multiplication:
\begin{equation}\label{scalar}
a(m+t)=a^{1/2}m+at.
\end{equation}
The topological $k$-space $\mathtt{U}V/H$ classifies adjoint quadratic maps on $V$ as follows.
\begin{prop}
Let $V$ be a linearly topologized $k$-space. The induced map $\iota:V\ra \mathtt{U}V/H$ has the following universal property: for any subfield $k'\subset k$ and any continuous adjoint and $k'$-quadratic map $q:V\ra W$ where $W$ is a linearly topologized  $k'$-space, there exists a unique continuous $k'$-linear map $\varphi: \mathtt{U}V/H\ra W$ such that $q=\varphi \iota$. 
\end{prop}
\begin{proof}
Direct calculation using (\ref{multi}) and (\ref{scalar}). The continuity of $\varphi$ follows from the continuity of $q$ through the universal property of tensor product. 
\end{proof}
\begin{coro}\label{adjoint-quadratic}
Let $k'$ be a subfield of $k$ and suppose $V\in \LC(k)$ and $W\in \LC(k')$. 
\begin{itemize}
\item[(1)] 
If $q:V\ra W$ is a continuous adjoint and $k'$-quadratic map, then there exists a unique continuous $k$-quadratic map $q':V\ra k\otimes_{k'}W$ such that $q=(\Tr_{k/k'}\otimes_{k'} 1_W)q'$.
\item[(2)] If $q:V\ra k'$ is a continuous adjoint $k'$-quadratic form, then there exists a continuous adjoint $k'$-bilinear form $s$ such that $q(x)=s(x,x)$ for all $x\in V$.
\end{itemize}
\end{coro}
\begin{proof}
(1) By the previous proposition, $q$ corresponds to a unique map in $\Hom_{k'}(\mathtt{U}V/H,W)$, and hence a unique map $\varphi\in \Hom_k(\mathtt{U}V/H,k\otimes_{k'}W)$ via Proposition \ref{duality}.(1). Then $q'=\varphi \iota:V\ra k\otimes_{k'} W$ is the unique $k$-quadratic map as required.  

\noindent (2) Taking $W=k'$ in (1), we obtain a $k$-quadratic form $q'$ satisfying $q=\Tr_{k/k'}q'$. By Remark \ref{classical-orthogonal} there is a $k$-bilinear form $s'$ such that $q'(x)=s'(x,x)$ and then $s=\Tr_{k/k'}s'$ is the required adjoint $k'$-bilinear form.
\end{proof}

\subsection{Some combinatorics on quadratic $k$-spaces}\label{statistics}
Let $(V,h,q)$ be a nondegenerate $(\lambda,\Lambda)$-quadratic $k$-space. Denote by $\CI_V$ the set of maximal isotropic subspaces in $V$.

\begin{lem}\label{dim2}
If $V$ is a $2$-dimensional metabolic quadratic $k$-space. Then $|\CI_V|=1+|\Lambda|$.
\end{lem}
\begin{proof}
Suppose $V$ is a $(\lambda,\Lambda)$-quadratic $k$-space. Since $V$  is metabolic, $V$ is isomorphic to the hyperbolic plane $\left(k^2, \left[\matrixx{0}{1}{0}{0}\right]\right)$.  Then
\[h=\matrixx{0}{1}{\lambda }{0},\quad q=x\sigma(y) \mod \Lambda.\]
Any line through the origin contains a unique representative point of the form $(1:0)$ or $(a:1)$ with  $a\in k$. The line is maximal isotropic if and only if the representative is isotropic, i.e.  the points $(1:0)$ and $(a:1)$ with $a\in \Lambda$. Thus there are $1+|\Lambda|$ isotropic lines through the origin. 
\end{proof}

Let $X\subset V$ be an isotropic subspace. The nondegenerate $\lambda$-hermitian form $h$ descends to a nondegenerate $\lambda$-hermitian form $h_X$ on $X^\perp/X$. The quadratic structure on $V$ determines a unique nondegenerate $(\lambda,\Lambda)$-quadratic $k$-space $(X^\perp/X,h_X,q_X)$. Indeed, if $p\neq 2$ or $\sigma\neq 1$, the induced quadratic structure is determined by the hermitian form $h_X$ (cf. Lemma \ref{auto-quadratic}); If $p=2$ and $\sigma=1$, it is determined by the quadratic form $q_X$  induced from $q$ (cf. Remark \ref{classical-orthogonal}).
If $W\in \CI_V$, then $X+(X^\perp\cap W)\in \CI_V$ and $(X+(X^\perp\cap W))/X\in \CI_{X^\perp/X}$. Define 
\begin{equation}\label{projection}
\pi_{V,X^\perp/X}:\CI_V\ra \CI_{X^\perp/X},\quad W\mapsto ((W\cap X^\perp)+X)/X.
\end{equation}
\begin{prop}\label{maximal-isotropic}
Let $V$ be a $2n$-dimensional metabolic $(\lambda,\Lambda)$-quadratic $k$-space. 
\begin{itemize}
\item[(a)] All fibers of $\pi_{V,X^\perp/X}:\CI_V\ra \CI_{X^\perp/X}$ have size $\prod_{i=1}^{\dim X}(|\Lambda|q^{n-i} +1)$.
\item[(b)] We have $|\CI_V|=\prod_{j=0}^{n-1}(|\Lambda|q^{j}+1)$.
\item[(c)] Let $W$ be a fixed maximal isotropic subspace of $V$. Let $X_{q,n}^\Lambda$ be the random variable $\dim (Z\cap W)$ where $Z$ is chosen uniformly at random from $\CI_V$. Then $X_{q,n}^\Lambda$ is a sum of independent Bernoulli random variables $B^\Lambda_{q,1},\cdots,B^\Lambda_{q,n}$ where $B^\Lambda_{q,i}$ is $1$ with probability $1/(1+|\Lambda|q^{i-1})$ and $0$ otherwise.
\item[(d)] For $0\leq d\leq n$, let $\alpha_{q,n,d}^\Lambda:=\Prob(X_{q,n}^\Lambda=d)$, and let $\alpha_{q,d}^\Lambda:=\lim_{n\ra \infty} \alpha_{q,n,d}^\Lambda$. Then
\[\sum_{d\geq 0}\alpha_{q,n,d}^\Lambda z^d=\prod_{i=0}^{n-1}\frac{1+|\Lambda|^{-1}q^{-i}z}{1+|\Lambda|^{-1}q^{-i}},\]
\[\sum_{d\geq 0}\alpha_{q,d}^\Lambda z^d=\prod_{i=0}^{\infty}\frac{1+|\Lambda|^{-1}q^{-i}z}{1+|\Lambda|^{-1}q^{-i}},\]
\item[(e)] For $0\leq d\leq n$, we have
\[\alpha_{q,n,d}^\Lambda=\prod_{j=0}^{n-1}(1+|\Lambda|^{-1}q^{-j})^{-1}\prod_{j=1}^d\frac{|\Lambda|^{-1}q}{q^j-1}\prod_{j=0}^{d-1}(1-q^{j-n}),\]
\[\alpha_{q,d}^\Lambda=\alpha_q^\Lambda\prod_{j=1}^d\frac{|\Lambda|^{-1}q}{q^j-1},\]
where 
\[\alpha_q^\Lambda=\alpha^\Lambda_{q,0}=\prod_{j=0}^{\infty}(1+|\Lambda|^{-1}q^{-j})^{-1}.\]
\end{itemize}
\end{prop}
\begin{proof}
The proof goes verbatim as in \cite[Proposition 2.6]{PR12} with necessary remedies for general metabolic spaces provided by Lemma \ref{dim2}.
\end{proof}

By Proposition \ref{maximal-isotropic}.(d), the function $d\mapsto \alpha^\Lambda_{q,d}$ is a distribution on the states $\BZ_{\geq 0}$. As in \cite[\S 2.3]{PR12}, to interpret the value $\alpha^\Lambda_{q,d}$ as a probability, not only as a limit of probabilities, one needs to consider infinite-dimensional quadratic spaces.  Let $V$ be a second-countable infinite-dimensional metabolic $(\lambda,\Lambda)$-quadratic $k$-space. Let $\CX_V$ be the set of compact open isotropic subspaces of $V$. The set $\CX_V$ is a countable directed poset under inclusion. For $X\in \CX_V$, the projection maps $\pi_{V,X^\perp/X}$ of (\ref{projection}) form an inverse system and induce a bijection
\[\CI_V\ra \lim_{\begin{subarray}{c}\lla\\X\in \CX_V\end{subarray}}\CI_{X^\perp/X}.\]
Equipping $\CI_V$ with the inverse topology, there exists a unique Borel measure $\mu$ on $\CI_V$ which pushes forward to the uniform probability measure on the finite set $\CI_{X^\perp/X}$ under $\pi_{V,X^\perp/X}$.  Mimicking the proof of \cite[Proposition 2.19]{PR12}, we have
\begin{prop}\label{infinite-quadratic}
Fix a compact open maximal isotropic subpace $W\in \CI_V$. If $Z\in \CI_V$ is distributed according to the measure $\mu$, then 
 \[\Prob(\dim_k(Z\cap W)=d)=\alpha_{q,d}^\Lambda.\]
\end{prop}

In view of the previous two propositions, the probability distribution of the dimension of the intersection of a fixed maximal isotropic subspace with a random one in a metabolic quadratic space $V$ over $k=\BF_q$ is determined by the cardinality of its form parameter $\Lambda$. If we denote $|\Lambda|=q^{\delta(\Lambda)}$, then we may denote $\alpha^\Lambda_{q,d}$ as $\alpha^{\delta(\Lambda)}_{q,d}$. Then we have the following dictionary:
\begin{center}
\tablefirsthead{
\hline  
$V$ & Form parameter $\Lambda$ & $\delta(\Lambda)$\\ 
\hline }

\bottomcaption{Type of spaces v.s. cardinality of form parameters}
\begin{supertabular}{|c|c|c|}\label{t1}
Orthogonal spaces&$\Lambda^{\min}=0$&0\\
\hline 
Symplectic spaces&$\Lambda^{\max}=k$&1\\
\hline 
Unitary spaces&$\Lambda^{\min}=\Lambda^{\max}=k_0(1-2\alpha)$& $\frac{1}{2}$\\
\hline 
\end{supertabular}
\end{center}
The unitary form parameter given here is only for the case $\lambda=1$ which is harmless (cf. \cite[Lemma 7-6.6]{Sch85}), and $\alpha$ is an element satisfying $\alpha+\sigma(\alpha)=1$.

\begin{defn}
Let $q=p^e$, $\delta\in \{0,\frac{1}{2},1\}$ and if $\delta=\frac{1}{2}$, we assume $e$ is even. Define $\sD_{q}^\delta:\BZ_{\geq 0}\ra \BR$ to be the probability distribution given as $\sD_q^\delta(d)=\alpha_{q,d}^\delta$ for all $d\in \BZ_{\geq 0}$.
\end{defn}
From Table \ref{t1}, to indicate where they come from, we classify these distributions into orthogonal distributions $\sD^\Ort_q$ ($\delta=0$), symplectic distributions $\sD^\Sym_q$ ($\delta=1$) and unitary distributions $\sD^\Uni_q$ ($e$ even and $\delta=1/2$).

Given a random variable $X$, let $\BE(X)$ be its expectation. If $m\in \BZ_{\geq 0}$, then $\BE(X^m)$ is its $m$-th moments. Let $X_{q}^\delta$ be the $\BZ_{\geq 0}$-valued random variable which has $\sD^\delta_q$ as its probability distribution.

\begin{prop}\label{behavior-distribution}
\begin{itemize}
\item[(1)] We have
\[\BE(X_q^\delta)=\sum_{i=0}^\infty \frac{1}{1+q^{\delta+i}}.\]
\item[(2)] For $m\geq 1$, 
\[\BE\left(\left(q^{X_{q}^\delta}\right)^m\right)=\prod_{i=1}^m (1+q^{-\delta+i}).\]
In particular, $\BE(q^{X_q^\delta})=1+q^{1-\delta}$.
\item[(3)] If we set $P(\delta)=\prod_{i\geq 0}\frac{1-q^{-\delta-i}}{1+q^{-\delta-i}}$, then 
\[\Prob\left(X^\delta_q \text{ is even}\right)=\frac{1+P(\delta)}{2}.\]
In particular, $\Prob\left(X^\delta_q \text{ is even}\right)=1/2$ if and only if $\delta=0$, i.e.  $X_q^\delta$ obeys  an orthogonal distribution.
 \item[(4)] We have 
\[\Prob\left(X^\delta_q =0 \right)>1-\frac{q^{1-\delta}}{q-1}.\]
If $\delta\neq 0$, then 
\[\lim_{q\ra \infty}\Prob\left(X^\delta_q =0 \right)=1. \]
\end{itemize}
\end{prop}
\begin{proof}
For (1), taking the derivative of the first identity of Proposition \ref{maximal-isotropic}.(d) and taking values at $z=1$, we have
\[\sum_{d=0}^n \alpha^\Lambda_{q,n,d} d =\sum_{i=0}^n \frac{|\Lambda|^{-1}q^{-i}}{1+|\Lambda|^{-1}q^{-i}}.\]
Then (1) follows as $n\ra \infty$. For (2)-(3), see the proof of \cite[Proposition 2.22]{PR12} and just take $z=q^m$ and $z=-1$ in Proposition \ref{maximal-isotropic}.(d). 

As for (4), taking $m=1$ in (2), 
\begin{eqnarray*}
1+q^{1-\delta}&=&\alpha_{q,0}^\delta+q \alpha_{q,1}+q^2\alpha_{q,2}+\cdots\\
&>& \alpha_{q,0}^\delta+q(1-\alpha^\delta_{q,0}).
\end{eqnarray*}
The inequality for $\alpha_{q,0}^\delta$ follows as desired.
\end{proof}
If the random variable $X_q^\delta$ is orthogonal, then it has equal parity. Otherwise, $X_q^\delta$ has a bias to taking even values, and even a bias to taking $0$ value. This bias of taking $0$ value tends to $100\%$ as $q\ra \infty$.
If $X_q^\delta\neq X_2^1$ is symplectic or unitary, then $\Prob\left(X^\delta_q =0 \right)>1/2$.

\subsection{Split unitary spaces}\label{split-unitary}
Let $k_0$ be a finite field of cardinality $q_0$ and $k=k_0\oplus k_0$ with the involution $\sigma(x,y)=(y,x)$. All definitions in Section \ref{hermitian-spaces}-\ref{quadratic-spaces} apply in this situation, and the quadratic structures are straight to describe.  Let $V$ be a locally compact $k$-space (by which we mean a $k$-module). Suppose $\lambda=\pm 1$ and $h:V\times V\ra k$ is a nondegenerate $\lambda$-hermitian form. By multiplying $(1,-1)\in k$ if necessary, we may assume $\lambda=1$. The elements $\alpha=(1,0)$ and $\sigma(\alpha)=(0,1)$ are idempotents in $k$ satisfying $\alpha\sigma(\alpha)=0$ and $\alpha+\sigma(\alpha)=1$. If we denote $V_0=\alpha V$, then $V_0$ is a $(k_0=\alpha k)$-space and the nondegenerate form $h$ identifies $\sigma(\alpha)V$ with the dual space $V_0^*=\Hom_{k_0}(V_0, k_0)$. Thus $V=V_0\bigoplus V_0^*$ as a $k$-space. As in the unitary case, the presence of the nontrivial involution forces the uniqueness of the form parameter: $\Lambda^{\min}=\Lambda^{\max}=k_0(1,-1)$. If we set $s=\alpha h=\left(\begin{matrix}0&1_{V_0^*}\\0&0\end{matrix}\right)$, then $h=s+s^*$ is even. For any $\sigma$-sesquilinear form $s'$ satisfying $h=s'+s'^*$, the class $[s']=[s]$ is uniquely determined by $h$. Thus $(V,h)$ determines a unique metabolic quadratic $k$-space $(\BH(V_0)=V_0\bigoplus V_0^*,[s])$. We call these quadratic $k$-spaces {\em split unitary $k$-spaces}. In particular, in any split unitary $k$-space, a $k$-subspace $X$ is maximal isotropic if and only if $X=X^\perp$ (cf. Remark \ref{quadratic-evenhermitian2}, Lemma \ref{auto-quadratic} and {\ref{orth-iso}).

Viewing $V$ as a $k_0$-space via the diagonal embedding $k_0\hookrightarrow k$, $V$ also acquires a metabolic orthogonal $k_0$-space structure $(\BH(V_0),[s_0])$ where $s_0=\Tr_{k/k_0}s$. Since these two quadratic structures are trace compatible, a $k$-subspace in $V$ is maximal isotropic in $(\BH(V_0),[s])$ if and only if it is maximal isotropic in 
$(\BH(V_0),[s_0])$ (cf. Theoreom \ref{lifting-thm}).
\begin{prop}
Any maximal isotropic $k$-subsapce of $(\BH(V_0),[s])$ is of the form
\[\wt{W}=W\bigoplus W^\circ\]
where $W$ runs through $k_0$-subspaces of $V_0$ and $W^\circ=W^\perp\bigcap V_0^*$.
\end{prop}
\begin{proof}
Since $\wt{W}^\perp=W^\perp\bigcap {W^\circ}^\perp=(V_0\bigoplus W^\circ)\bigcap (W\bigoplus V_0^*)=\wt{W}$, $\wt{W}$ is maximal isotropic. Suppose $Z$ is a maximal isotropic subspace. Since $\alpha+\sigma(\alpha)=1$, $Z=\alpha Z\bigoplus \sigma(\alpha)Z$. It suffices to prove $\sigma(\alpha)Z=(\alpha Z)^\perp\bigcap V_0^*$. Clearly $\sigma(\alpha)Z\subset (\alpha Z)^\perp\bigcap V_0^*$ and we prove the converse. Let $f$ be an arbitrary element of $(\alpha Z)^\perp\bigcap V_0^*$. For any $z\in Z$,
\[h(z,f)=h(z, \sigma(\alpha)f)=h(\alpha z, f)=0.\]
Thus $f\in Z^\perp=Z$. Since $f\in V_0^*$, $f=\sigma(\alpha) f \in \sigma(\alpha) Z$ as desired. 
\end{proof}
Thus the combinatorics of maximal isotropic $k$-subspaces in a split unitary space $\BH(V_0)$ is essentially the same as that of $k_0$-subspaces in $V_0$. In particular, If $V_0$ has finite dimension $n$, then the set $\CI_{\BH(V_0)}$ of maximal isotropic subspaces of the split unitary space $\BH(V_0)$ has the cardinality as the number of subspaces of $V_0$, i.e.  $|\CI_{\BH(V_0)}|=\sum_{i=0}^n\qbinomial{n}{i}_{q_0}$ where $\qbinomial{n}{i}_{q_0}$ denote the $q_0$-binomial coefficients. Fix a maximal isotropic subspace $\wt{W}=W\bigoplus W^\circ\in \CI_{\BH(V_0)}$ with $\dim_{k_0} W=d_W$. A standard linear argument shows that the probability of $\dim_{k_0} \wt{W}\bigcap \wt{Z}=r$ as $\wt{Z}$ runs uniformly through $\CI_{\BH(V_0)}$ is given by
\[\left. \left(\sum^n_{\begin{subarray}{c}d=0\\
d-d_W\equiv n-r \mod 2 \end{subarray}}q_0^{\frac{(n-r)^2-(d-d_W)^2}{4}}\qbinomial{d_W}{\frac{d+d_W+r-n}{2}}_{q_0}\qbinomial{n-d_W}{\frac{d-d_W-r+n}{2}}_{q_0}\right)\right/|\CI_{\BH(V_0)}|,\]
which depends on the dimension of $W$.

Thus the maximal isotropic subspaces in a split unitary space behave very differently from those in a quadratic space over a finite field and various limit probability distributions can be obtained as the dimension goes to infinity. For example, the well known ``uniform distribution" for co-rank of random matrices can be recoverd as in Proposition \ref{uniform-distribution}. 

\begin{defn}\label{uniform-distribution0}
Let $k_0=\BF_{q_0}$ and $t\geq 0$ be an integer. Define the distribution $\sU_{q_0}^t$ on non-negative integers
\[\sU_{q_0}^t(r)=q_0^{-r(t+r)}\frac{\prod_{j=r+1}^\infty (1-q_0^{-j})}{\prod_{j=1}^{t+r}(1-q_0^{-j})}.\]
\end{defn}
It is known that, as $s\ra \infty$, $P(s,s+t,r)$ converges to $\sU_{q_0}^t(r)$ where $P(s,s+t,r)$ denotes the probability of a uniformly random $s\times(s+t)$ matrix over $k_0$ having co-rank $r$ (cf. \cite[\S1]{Gerth86, FG15}}).

Fix integers $m_0$ and $m_1$ of the same parity. For large $n\equiv m_0\mod 2$, let $V_0$ be a $k_0$-space of dimension $n$ and let $W$ be a fixed subspace of $V_0$ of dimension $\frac{n+m_0}{2}$. Since $\GL_n(k_0)$ acts transitively on the set of subspaces of $V_0$ of a given dimension.  The probability of $\dim_{k_0} W\cap Z=r$, as $Z$ runs uniformly through subspaces of dimension $\frac{n+m_1}{2}$, depends only on the dimensions of $W$ and $Z$ and we denote it as $P_{m_0,m_1,n,q_0}(r)$.

Let $V=\BH(V_0)$ be the split unitary $k$-space of $k_0$-dimension $2n$. A maximal isotropic subspace $\wt{W}=W\bigoplus W^\circ$ is said of type $m$ if $\dim_{k_0}W=\frac{n+m}{2}$ while $\dim_{k_0}W^\circ=\frac{n-m}{2}$.

\begin{prop} \label{uniform-distribution}
Fix integers $m_0$ and $m_1$ of the same parity and let the notations be as above.
\begin{itemize} 
\item[(1)] For $r\geq \max\left(0,\frac{m_0+m_1}{2}\right)$, we have \[\lim_{n\ra \infty} P_{m_0,m_1,n,q_0}(r)=\sU_{q_0}^{\left|\frac{(m_0+m_1)}{2}\right|}\left(r-\max\left(0,\frac{m_0+m_1}{2}\right)\right).\]
\item[(2)] Let $\wt{W}\subset \BH(V_0)$ be a fixed maximal isotropic subspace of type $m_0$. For $r\geq |\frac{m_0+m_1}{2}|$, the probability of $\dim_{k_0} \wt{W}\cap \wt{Z}=r$, as $\wt{Z}$ runs uniformly through maximal isotropic subspaces of type $m_1$, converges to $\sU_{q_0}^{\left|\frac{(m_0+m_1)}{2}\right|}\left(\frac{r+(m_0+m_1)/2}{2}-\max\left(0,\frac{m_0+m_1}{2}\right)\right)$ if $r\equiv \frac{m_0+m_1}{2}\mod 2$ and $0$ otherwise.  
\end{itemize}
\end{prop}
\begin{proof}
First from the dimension formula, we have
\[\max\left(0,\frac{m_0+m_1}{2}\right)\leq \dim_{k_0} W\cap Z\leq \max\left(\frac{n+m_0}{2},\frac{n+m_1}{2}\right).\]
As is noted that the probability $P_{m_0,m_1,n,q_0}(r)$ depends only on the dimensions of $W$ and $Z$, it is also the probability of $\dim_{k_0} W\cap Z=r$, as $W$ runs through subspaces of dimension $\frac{n+m_0}{2}$ and $Z$ runs through subspaces of dimension $\frac{n+m_1}{2}$. For an $n\times (\frac{n+m_0}{2})$ matrix $A$ and an $n\times (\frac{n+m_1}{2})$ matrix $B$ over $k_0$, the matrix $(A,B)\in \mathrm{M}_{n\times \left(n+\frac{m_0+m_1}{2}\right)}(k_0)$ is called admissible if both $A$ and $B$ are of full rank (Note $n$ is large) and denote the set of admissible matrices by $\mathrm{M}^{\mathrm{adm}}$. There is a surjection with fibers of the same cardinality
\[\mathrm{M}^{\mathrm{adm}}\ra \left\{(W,Z): \begin{aligned}&\text{pairs of subspaces of $V_0$ of}\\
&\text{dimension $\left(\frac{n+m_0}{2},\frac{n+m_1}{2}\right)$}\end{aligned}\right\}\]
\[(A,B)\mapsto (\mathrm{Span}(A),\mathrm{Span}(B)),\]
where $\mathrm{Span}(A)$ denotes the subspace generated by the column vectors of $A$. Suppose $(A,B)$ maps to $(W,Z)$. If $\frac{m_0+m_1}{2}\geq 0$, then $\corank(A,B)=\dim_{k_0}(W\bigcap Z)-\frac{m_0+m_1}{2}$; Otherwise $\corank(A,B)=\dim_{k_0}(W\bigcap Z)$. Thus if we denote $t=\max\left(0,\frac{m_0+m_1}{2}\right)$, then
\begin{eqnarray*}
P_{m_0,m_1,n,q_0}(r)&=&\frac{|\{X\in \mathrm{M}^{\mathrm{adm}}: \corank(X)=r-t\}|}{|\mathrm{M}^{\mathrm{adm}}|}\\
&=&\frac{|\mathrm{M}|}{| \mathrm{M}^{\mathrm{adm}}|}\left(\frac{|\{X\in \mathrm{M}: \corank(X)=r-t\}|}{|\mathrm{M}|}-\frac{|\{X\notin \mathrm{M}^{\mathrm{adm}}: \corank(X)=r-t\}|}{|\mathrm{M}|}\right).
\end{eqnarray*}
Here $\mathrm{M}$ simply denotes $\mathrm{M}_{n\times \left(n+\frac{m_0+m_1}{2}\right)}(k_0)$.  As $n\ra \infty$, $\frac{n+m_0}{2}\sim \frac{n}{2}$, and thus, for a random $A$, the probability of $A$ having full rank, i.e.  co-rank zero, tends to $1$; The same holds for a random $B$. Thus as $n\ra \infty$, $\frac{| \mathrm{M}^{\mathrm{adm}}|}{|M|}\ra 1$ and $P_{m_0,m_1,n,q_0}(r)\ra\sU_{q_0}^{\left|\frac{(m_0+m_1)}{2}\right|}(r-t)$. 

The second assertion follows from (1) and the following identity
\[\dim_{k_0} \wt{W}\bigcap \wt{Z}=-\frac{m_0+m_1}{2}+2\dim_{k_0} W\bigcap Z.\]
\end{proof}

\begin{defn}\label{type-density}
For the parity $\epsilon\in \{0,1\}$, define the distribution $\sV^{\epsilon}_{q_0}$ on $\BZ$ as
\[\sV^\epsilon_{q_0}(n)=q_0^{-n(n+\epsilon)}\left(\prod_{i=1}^{\infty}(1-q_0^{-2i})(1+q_0^{-(2i-1+\epsilon)})(1+q_0^{-(2i-1-\epsilon)})\right)^{-1}. \]
\end{defn}
The fact that $\sV_{q_0}^\epsilon$ is indeed a probability distribution follows from the following identity
\begin{equation}\label{Jacobi}
\sum_{-\infty}^{\infty} q_0^{-n(n+\epsilon)}=\prod_{i=1}^{\infty}(1-q_0^{-2i})(1+q_0^{-(2i-1+\epsilon)})(1+q_0^{-(2i-1-\epsilon)}),
\end{equation}
which is obtained by substituting $q=q_0^{-1}$ and $z=q_0^{-\epsilon}$ in Jacobi's triple product formula \cite[Theorem 11.1]{KC01}.
\begin{prop}\label{mi-density}
Fix $m\in \BZ$ of parity $\epsilon$. Let $V_0$ be a $k_0$-space of dimension $n$ with $n$ large of parity $\epsilon$. The proportion of subspaces of dimension $\frac{n+m}{2}$ in $V_0$, or equivalently, the proportion of maximal isotropic subspaces of type $m$ in $\BH(V_0)$, converges to $\sV^\epsilon_{q_0}\left(\frac{m-\epsilon}{2}\right)$ as $n\ra \infty$.
\end{prop}
\begin{proof}
Let $G(n,q_0)$ be the Galois number of $V_0$, i.e.  the number of subspaces of $V_0$. It suffices to compute the limit of the ratio $\qbinomial{n}{\frac{n+m}{2}}_{q_0}/G(n,q_0)$ as $n\ra \infty.$ Set $n=2k-\epsilon$. We have the symmetry $\qbinomial{n}{i}_{q_0}=\qbinomial{n}{n-i}_{q_0}$ and for any $d\geq 0$, as $k\ra \infty$,
\[\left.\qbinomial{n}{k+d}_{q_0}\right/\qbinomial{n}{k}_{q_0}\ra q_0^{-d(d+\epsilon)}.\]
Therefore the proposition follows from Jacobi's identity (\ref{Jacobi}) and Tannery's theorem \cite[\S 49]{Brom64} as desired. 
\end{proof}

We prove the limit of the distribution of the dimension of the intersection of a fixed maximal isotropic subspace with a uniformly random one in $\BH(V_0)$ of dimension $2n$ exists as $n\ra \infty$.
\begin{prop} 
Fix an integer $m$ of parity $\epsilon$ and an integer $r\geq 0$. Let $n$ be a large integer of parity $\epsilon$ and let $\BH(V_0)$ be a split unitary $k$-space of $k_0$-dimension $2n$ with a fixed maximal isotropic subspace $\wt{W}$ of type $m$. As $n\ra \infty$, the probability of $\dim_{k_0} \wt{W}\cap \wt{Z}=r$, as $\wt{Z}$ runs uniformly through $\CI_{\BH(V_0)}$, converges to, 
\begin{equation}\label{su-dist}
q_0^{-\left(r+\frac{m+\epsilon}{2}\right)\left(r+\frac{m+\epsilon}{2}+\epsilon\right)}\left(\sum_{k=0}^{r} q_0^{2(k+\epsilon)(r-k)+(2k+\epsilon)(m+\epsilon)}\qbinomial{r}{k}_{q_0}\right)\prod_{r+1}^\infty(1-q_0^{-j})\sV_{q_0}^\epsilon(0).
\end{equation}
\end{prop}
\begin{proof}
Fix $r\geq 0$. For $\wt{Z}$ of type $m_1$, in order that $\dim_{k_0} \wt{W}\bigcap \wt{Z}=r$, necessarily $r\geq |\frac{m+m_1}{2}|$, i.e.  $-2r-m\leq m_1\leq 2r-m$. Moreover, by Proposition \ref{uniform-distribution}, $r\equiv \frac{m+m_1}{2} \mod 2$. We have $m_1=-2r-m+4k$ with $0\leq k\leq r$. Thus
\[\Prob(\dim \wt{W}\bigcap \wt{Z}=r)=\sum_{k=0}^r \Prob(\type(\wt{Z})=m_1)A_{m,m_1,n}(r),\]
where $A_{m,m_1,n}(r)$ denotes the probability of $\dim \wt{W}\bigcap \wt{Z}=r$ with $\wt{Z}$ of type $m_1$. As $n\ra\infty$, by Proposition \ref{uniform-distribution} and \ref{mi-density}, take the limits and we obtain
\begin{eqnarray*}
\lim\limits_{n\ra \infty}\Prob(\dim \wt{W}\bigcap \wt{Z}=r)&=&\sum_{0\leq k\leq r/2} \sV^\epsilon_{q_0}\left(-r+2k-\frac{m+\epsilon}{2}\right)\sU_{q_0}^{r-2k}(k)\\
&+&\sum_{r/2<k\leq r} \sV^\epsilon_{q_0}\left(-r+2k-\frac{m+\epsilon}{2}\right)\sU_{q_0}^{2k-r}(r-k).
\end{eqnarray*}
Taking into account of the explicit expressions of $\sV_{q_0}^\epsilon$ and $\sU^t_{q_0}$, the desired formula (\ref{su-dist}) follows.
\end{proof}
Accordingly we make the following definition.
\begin{defn}
For any integer $m$, define the distribution $\sD_{q_0}^{\SU,m}$ on non-negative integers such that $\sD_{q_0}^{\SU,m}(r)$ is given by the expression $(\ref{su-dist})$.
\end{defn}
Flipping the roles of $V_0$ and $V_0^*$ in $\BH(V_0)$, we see $\sD_{q_0}^{\SU,m}=\sD_{q_0}^{\SU,-m}$. It is not obvious at first glance that $\sD_{q_0}^{\SU,m}$ is a distribution, i.e.  the total mass of $\sD_{q_0}^{\SU,m}$ is equal to $1$. Interestingly, it is equivalent to a type of Rogers-Ramanujan identities (cf. \cite[\S 7.3]{And76}) whose proof is communicated to me by Lishuang Lin. 
\begin{prop}
Let $m$ be an integer of parity $\epsilon$. For $|t|<1$, we have the following Rogers-Ramanujan type identity
\begin{equation}\label{t-identity}
\sum_{k=0}^\infty\sum_{i=0}^\infty \frac{t^{k^2+i^2-ki+\epsilon(k-i)-(k-i)(m+\epsilon)}}{(t,t)_k(t,t)_i}=(-t,t)_\infty(-t^{1+\epsilon},t^2)_\infty(-t^{1-\epsilon},t^2)_\infty
\end{equation}
where $(a,t)_n=(1-a)(1-at)\cdots(1-at^{n-1})$ denotes the $t$-Pochhammer symbol. 
\end{prop}
\begin{proof}
Rearrage the exponent of $t$ as 
\begin{eqnarray*}
&&\frac{i(i-\epsilon)}{2}+\frac{i(m+\epsilon)}{2}+\frac{k(k+\epsilon)}{2}-\frac{k(m+\epsilon)}{2}+\\
&&\frac{(k-i)(k-i+\epsilon)}{2}-\frac{(k-i)(m+\epsilon)}{2}+\left(\frac{m+\epsilon}{2}\right)^2-\frac{\epsilon(m+\epsilon)}{2}.
\end{eqnarray*}
Define
\[F_1(z)=\sum_{i=0}^\infty \frac{t^{\frac{i(i-\epsilon)}{2}}}{(t,t)_i}z^i,\quad F_2(z)=\sum_{k=0}^\infty \frac{t^{\frac{k(k+\epsilon)}{2}}}{(t,t)_k}z^{-k}\text{ and } F_3(z)=\sum_{n=-\infty}^{\infty} t^{\frac{n(n+\epsilon)}{2}-n(m+\epsilon)}z^n.\]
Thus the LHS of (\ref{t-identity}) is the constant coefficient of $t^{\left(\frac{m+\epsilon}{2}\right)^2-\frac{\epsilon(m+\epsilon)}{2}}(F_1F_2F_3)(t^{\frac{m+\epsilon}{2}}z)$, or equivalently, the constant coefficient of $t^{\left(\frac{m+\epsilon}{2}\right)^2-\frac{\epsilon(m+\epsilon)}{2}}(F_1F_2F_3)(z)$.

By Euler's identity \cite[(9.3), p30]{KC01},
\[F_1(z)=\left(-t^{\frac{1-\epsilon}{2}}z,t\right)_\infty \text{ and }F_2(z)=\left(-t^{\frac{1+\epsilon}{2}}z^{-1},t\right)_\infty.\]
Also by Jacobi's triple product identity \cite[Theorem 11.1]{KC01}, 
\[\sum_{\ell=-\infty}^\infty t^{\frac{\ell(\ell-\epsilon)}{2}}z^\ell=(t,t)_\infty\left(-t^{\frac{1-\epsilon}{2}}z,t\right)_\infty \left(-t^{\frac{1+\epsilon}{2}}z^{-1},t\right)_\infty.\]
Thus
\begin{eqnarray*}
(t,t)_\infty(F_1F_2F_3)(z)&=&\left(\sum_{\ell=-\infty}^\infty t^{\frac{\ell(\ell-\epsilon)}{2}}z^\ell\right)\left(\sum_{n=-\infty}^{\infty} t^{\frac{n(n+\epsilon)}{2}-n(m+\epsilon)}z^n\right),
\end{eqnarray*}
and the constant coefficient of $t^{\left(\frac{m+\epsilon}{2}\right)^2-\frac{\epsilon(m+\epsilon)}{2}}(F_1F_2F_3)(z)$ is given as
\begin{eqnarray*}
(t,t)^{-1}_\infty \sum_{n=-\infty}^\infty  t^{n(n+\epsilon)-n(m+\epsilon)+\left(\frac{m+\epsilon}{2}\right)^2-\frac{\epsilon(m+\epsilon)}{2}} &=& \frac{\sum_{n=-\infty}^\infty  t^{\left(n+\frac{m+\epsilon}{2}\right)\left(n+\epsilon+\frac{m+\epsilon}{2}\right)}}{(t,t)_\infty}\\
&=&(-t,t)_\infty(-t^{1+\epsilon},t^2)_\infty(-t^{1-\epsilon},t^2)_\infty
\end{eqnarray*}
 as desired by Jacobi's identity (\ref{Jacobi}).
\end{proof}

It is ready to verify that $\sD_{q_0}^{\SU,m}$ has total mass $1$. Substituting $t=q_0^{-1}$, interchanging the summation and substituting $i=r-k$,
\begin{eqnarray*}
\sum_{r=0}^\infty \sD_{q_0}^{\SU,m}(r)&=&\sum_{r=0}^\infty \sum_{k=0}^{r} t^{\left(r+\frac{m+\epsilon}{2}\right)\left(r+\frac{m+\epsilon}{2}+\epsilon\right)-2(k+\epsilon)(r-k)-(2k+\epsilon)(m+\epsilon)-k(r-k)}\qbinomial{r}{k}_{t} \frac{(t,t)_\infty}{(t,t)_r}\sV^\epsilon(0)\\
&=&(t,t)_\infty \sV^\epsilon(0) \sum_{k=0}^{\infty} \sum_{i=0}^\infty  \frac{t^{k^2+i^2-ki+\epsilon(k-i)-(k-i)(m+\epsilon)}}{(t,t)_k(t,t)_i}=1
\end{eqnarray*}
as desired by invoking Definition \ref{type-density} and the identity (\ref{t-identity}).

\section{Weil pairings and quadratic refinements}\label{WPAV}
Let $F$ be a global field of characteristic $\ell$.  Let $A/F$ be an abelian variety. The endomorphism algebra $\End^0(A)=\End_F(A)\otimes \BQ$ is a semi-simple algebra.  If $A$ is simple, by Albert classification (see \cite[Theorem 2, p.201]{Mumford-AV}), $\End^0(A)$ is a central simple division algebra with center either a totally real or CM number field.  Let $K\subset \End^0(A)$ be a number field endowed with an involution $\dag$. Suppose the order $\CO= K\cap \End_F(A)$ is stable under $\dag$. Let $K_0\subset K$ be the subfield fixed by $\dag$ and put $\CO_0=\CO\cap K_0$.

Let $A^\vee$ be the dual abelian variety. For any $\phi\in \End_F(A)$, denote by $\phi^\vee$ its dual map. The assignment $a\mapsto a^\vee$ gives rise to an embedding $\CO\hookrightarrow \End_F(A^\vee)$.  Let $\Hom_\dag^\sym(A,A^\vee)$ denote group of $\dag$-sesquilinear symmetric homomorphisms $\mu$, i.e. $\mu a^\dag=a^\vee\mu$ for all $a\in \CO$. If $\mu$ is a symmetric isogeny, then $\dag$-sesquilinearity exactly means  the Rosati involution associated to $\mu$ stabilizes $\CO$ and restricts to the involution $\dag$ on $K$. 

Let $p$ be a rational prime number not dividing  $[\CO_K:\CO]$. Then any prime ideal of $\CO$ above $p$ is invertible.  Let $\fp$ be such a prime ideal and then $\fp_0=\fp\cap \CO_0$ is also an invertible prime ideal of $\CO_0$. Let $\CO_{0,\fp_0}$ be the completion of $\CO_0$ with respect to the prime ideal $\fp_0$ and $\CO_{\fp_0}=\CO\otimes_{\CO_{0}}\CO_{0,\fp_0}$. Then $K_{\fp_0}=\BQ\otimes_\BZ\CO_{\fp_0}$ is a separable  quadratic extension of ${K}_{0,\fp_0}=\BQ\otimes_{\BZ}\CO_{0,\fp_0}$ and the involution $\dag$ induces the nontrivial automorphism of the extension $K_{\fp_0}/K_{0,\fp_0}$.  Both $\CO_{\fp_0}$ and $\CO_{0,\fp_0}$ are rings of integers of $K_{\fp_0}$ and $K_{0,\fp_0}$ respectively.

\subsection{$p$-adic Tate groups and extended Weil pairings}\label{local-tate} 
For a nonzero integer $m$, the multiplication map $m:A\ra A$ is finite, faithfully flat. Thus the kernel $A[m]=\Ker(m)$ is a finite flat group scheme over $F$.  For $m\geq n\geq 0$, the morphisms $p^{m-n}:A[p^m]\ra A[p^n]$ form an inverse system with affine transition maps. The inverse limit $T_pA =\varprojlim A[p^n]$, the $p$-adic Tate group scheme, exists in the category of $F$-schemes. The canonical morphism $T_pA\ra A[p^n]$ is faithfully flat with kernel $p^nT_pA$ (cf. \cite[\S 8.2-8.3]{EGA-4-3}). As an inverse limit of finite ableian $p$-power torsion groups, $T_pA$ acquires a $\BZ_p$-module structure, and the inclusion $\CO\subset \End_F(A)$ furnishes $T_p A$ with a module structure over $\CO_p=\CO\otimes_\BZ \BZ_p$.

The Weil pairings $e_{p^n}:A[p^n]\times A^\vee[p^n]\lra \mu_{p^n}\subset \BG_m$
also form an inverse system of morphisms \cite[Proposition 11.21]{Edixhoven-AV}. By passage to limit, we get a natural pairing 
\[E_p:T_pA\times T_p A^\vee\lra \BZ_p(1):=\varprojlim \mu_{p^n}.\]
The pairing $E_p$ is $\BZ_p$-bilinear and nondegenerate, by which we mean, $E_p$ induces an isomorphism of $\BZ_p$-sheaves
\[T_pA^\vee\cong \bhom_{\BZ_p}(T_p A,\BZ_p(1)).\]  
Moreover, 
\begin{equation}\label{sesqui0}
E_p(a x,y)=E_p(x,a^\vee y) \text{ for any $a\in \CO$.}
\end{equation} 
Composing a symmetric homomorphism $\mu\in \Hom_\dag^\sym(A,A^\vee)$, we obtain an alternating pairing (see \cite[pp. 228-230]{Mumford-AV})
$$
E_p^\mu:T_p A\times T_p A\ra \BZ_p(1),\quad E_p^\mu(x,y)=E_p(x,\mu(y)).
$$

For prime ideals $\fP$ of $\CO_0$ above $p$, denote $T_\fP A=T_pA\otimes_{\CO_p}\CO_{\fP}$ the Serre tensor (see \cite[Theorem 7.2]{Conrad04}).  There is a decomposition
\[T_pA=\bigoplus_{\fP\mid p} T_{\fP} A\]
through which $E_p^\mu$ factors.  Denote by
\[E^\mu_{\fp_0}: T_{\fp_0} A\times T_{\fp_0} A\ra \BZ_p(1)\]
the $\fp_0$-component of $E^\mu_p$. By the $\dag$-sesquilinearity of $\mu$ and (\ref{sesqui0}), $E^\mu_{\fp_0}$ is $\dag$-adjoint, i.e.  for all $a\in \CO_{\fp_0}$, 
\[E^\mu_{\fp_0}(ax,y)=E^\mu_{\fp_0}(x,a^\dag y).\] 
Fix a generator $d\in \partial^{-1}_{\CO_{\fp_0}/\BZ_p}$ of the inverse different. By Proposition \ref{sesquilinear-adjoint}, there exists a unique $\dag$-sesquilinear pairing
\begin{equation}\label{Weil2}
\Theta^\mu_{\fp_0}: T_{\fp_0} A\times T_{\fp_0} A\ra \CO_{\fp_0}(1):=\CO_{\fp_0}\otimes_{\BZ_p}\BZ_p(1)
\end{equation}
which is characterized by the following property: for any $x,y \in T_{\fp_0}A \text{ and }a\in \CO_{\fp_0},$
\begin{equation}\label{characterization}
E^\mu_{\fp_0}(ax,y)=E^\mu_{\fp_0}(x,a^\dag y)=\Tr_{K_{\fp_0}/\BQ_p}(da \Theta_{\fp_0}^\mu(x,y) ).
\end{equation}

If $\dag\neq 1$, i.e.  $K/K_0$ is a quadratic extension, by the classification of local quadratic extensions \cite[\S2-3]{Casselman-quadratic}, we can require the fixed generator $d\in \partial^{-1}_{\CO_{\fp_0}/\BZ_p}$ satisfying
\begin{equation}\label{gen-dif}
d^{\dag-1}:=d^\dag/d=\left\{
\begin{aligned}
1,\quad & \text{if $\fp_0$ is inert or split in $K/K_0$};\\
-1,\quad & \text{if $\fp_0$ is ramified  in $K/K_0$}.
\end{aligned}
\right.
\end{equation}
In the following we fix such a generator $d\in \partial^{-1}_{\CO_{\fp_0}/\BZ_p}$ once and for all.

\begin{lem}\label{symmetric-effect}
For  $a\in \CO_0$ and $\mu\in \Hom_\dag^\sym(A,A^\vee)$, the pairing $(x,y)\mapsto E_p^\mu(x,ay))$ is alternating. 
\end{lem}
\begin{proof}
For any $x\in T_p A$,
\[E_p^\mu(x,ax)=E_p^\mu(a^\dag x,x)=E_p^\mu(ax,x)=E_p^\mu(x,ax)^{-1}.\]
Since  $E_p^\mu$ takes values in the torsion-free group $\BZ_p(1)$, $E_p^\lambda(x,ax)=1$ as desired.
\end{proof}

\begin{prop}\label{extended-Weil}
Let $\mu\in \Hom_\dag^\sym(A,A^\vee)$ be a symmetric homomorphism.
\begin{itemize}
\item[(1)] If $\dag=1$,  $\Theta_{\fp_0}^\mu$ is $\CO_{\fp_0}$-bilinear and alternating. 
\item[(2)] If $\dag\neq 1$,   $\Theta_{\fp_0}^\mu$ is $\dag$-sesquilinear and $(-d^{(\dag-1)})$-hermitian, i.e.  for any $x,y\in T_{\fp_0}A$, 
\[\Theta_{\fp_0}^\mu(x,y)^\dag=-d^{(\dag-1)}\Theta_{\fp_0}^\mu(y,x).\]
\end{itemize}
Moreover, if $\mu$ is an isogeny of degree prime to $p$, then $\Theta_{\fp_0}^\mu$ is nondegenerate. 
\end{prop}
\begin{proof}
By the construction, the pairing $E_{\fp_0}^\mu$ is $\dag$-adjoint and alternating, in particular, $(-1)$-hermitian. By Proposition \ref{sesquilinear-adjoint}, the corresponding paring $\Theta_{\fp_0}^\mu$ is $\dag$-sesquilinear and $(-d^{(\dag-1)})$-hermitian. By Lemma \ref{symmetric-effect}, for $x\in T_{\fp_0}A$ and all $a\in \CO_{0,\fp_0}$,
\begin{equation*}
\Tr_{K_{\fp_0}/\BQ_p}(da\Theta_{\fp_0}^\mu(x,x))=E_{\fp_0}^\mu(x,ax)=0.
\end{equation*}
If $K=K_0$, then $\CO_{\fp_0}=\CO_{0,\fp_0}$ and by the nondegeneracy of the trace form we conclude that $\Theta_{\fp_0}^\lambda(x,x)$=0 for all $x\in T_{\fp_0}A$.

If $\mu$ is an isogeny of degree prime to $p$, then $E_{\fp_0}^\mu$ is nondegenerate, and by Proposition \ref{sesquilinear-adjoint}, $\Theta_{\fp_0}^\mu$ is nondegenerate.
\end{proof}

From now on, till the end of this section, we assume the prime $\fp$ is stable under $\dag$, i.e.  either that $\dag=1$, or that $\dag\neq 1$ and $\fp_0$ is inert or ramified in $K/K_0$. In this case, $\fp$ is the unique prime ideal of $\CO$ above $\fp_0$. We fix a local generator $\varpi\in \fp$ so that $\varpi$ generates $\fp/\fp^2$. If $\dag\neq 1$ and $\fp$ is inert in $K/K_0$, we require additionally $\varpi^\dag=\varpi$. Let $k=\CO/\fp$ and $k_0=\CO_0/\fp_0$ be the residue fields. If $\dag=1$, or that $\dag\neq 1$ and $\fp_0$ is ramified in $K/K_0$, then $k=k_0$ and the involution $\dag$ induces the trivial action on $k$; Otherwise, $\dag$ induces the nontrivial involution on $k$ with fixed field $k_0$.

Denote by $e$ the ramification index of $\fp$ over $p$. Note $p/\varpi$ is $\fp$-adic integral. The morphism
\[T_{\fp} A=\varprojlim_n A[\fp^{en}]\lra A[\fp], \quad (\cdots, x_2,x_1)\mapsto (p/\varpi)x_1\]
induces an identification $T_{\fp} A/\fp T_{\fp} A\cong A[\fp]$ and we fix this identification throughout. Reducing (\ref{Weil2}) modulo $\fp$ yields a pairing 
\[\theta_\fp^\mu: A[\fp]\times A[\fp]\ra k(1):=k\otimes_{\BZ_p} \BZ_p(1).\]
As an immediate consequence of Proposition \ref{extended-Weil}, one has
\begin{coro}\label{mod-Weil}
Let $\mu\in \Hom_\dag^\sym(A,A^\vee)$ be a symmetric homomorphism.
\begin{itemize}
\item[(1)] If $\dag=1$,  $\theta_\fp^\mu$ is $k$-bilinear and alternating. 
\item[(2)] If $\dag\neq 1$ and $\fp$ is ramified in $K/K_0$, $\theta_\fp^\mu$ is $k$-bilinear and symmetric. 
\item[(3)] If $\dag\neq 1$ and $\fp$ is inert in $K/K_0$, $\theta_\fp^\mu$ is $\dag$-sesquilinear and $(-1)$-hermitian, i.e.  for any $x,y\in A[\fp]$, 
\[\theta_\fp^\mu(x,y)^\dag=-\theta_\fp^\mu(y,x).\]
\end{itemize}
Moreover, if $\mu$ is an isogeny of degree prime to $p$, then $\theta_{\fp}^\mu$ is nondegenerate. 
\end{coro}

\subsection{Induced Weil pairings on torsions}\label{IWP}
One may refer to Appendix \ref{Hom-construction} for the following constructions and notations. The injection $i_\fp: \fp\ra \CO$ and the identity map $1_A$  induces an isogeny, i.e.  the ``$\fp$-multiplication" map,
 \[h(i_\fp,1_A): A=\shom{\CO}{A}\lra \shom{\fp}{A}\]
with kernel $\shom{\CO/\fp}{A}$. The dual isogeny is
\[h(j_\fp,1_{A^\vee}):\shom{\fp^{-1}}{A^\vee}\lra \shom{\CO}{A^\vee}=A^\vee\]
with kernel $\shom{\fp^{-1}/\CO}{A^\vee}$, where $j_\fp:\CO\ra \fp^{-1}$ is the inclusion (cf. Proposition \ref{dual}). Associating sections with their images at $1$, we identify $\shom{\CO/\fp}{A}$ with $A[\fp]$. Fix $t_0\in \fp^{-1}$ once and for all such that $t_0\varpi\equiv 1\mod \fp$, and necessarily $t_0$ generates $\fp^{-1}/\CO$. Similarly,  associating sections with their images at the generator $t_0$ of $\fp^{-1}/\CO$, $\shom{\fp^{-1}/\CO}{A^\vee}$  is identified with $A^\vee[\fp]$ (This identification is independent of the choice of such a $t_0$, but depends on the choice of $\varpi$). The Cartier duality between the kernels of the dual isogenies $h(i_\fp,1_A)$ and $h(j_\fp,1_{A^\vee})$ induces the Weil pairing
\begin{equation}\label{Weil-p}
e_\fp: A[\fp]\times A^\vee[\fp]\lra \mu_p.
\end{equation}

For any $\mu\in \Hom^\sym_\dag(A,A^\vee)$, denote 
\[e_\fp^\mu: A[\fp]\times A[\fp]\lra \mu_p,\quad e_\fp^\mu(x,y)=e_\fp(x,\mu(y)).\]
If $\mu=0$, then $e_\fp^\mu=1$ identically. For $\mu\neq 0$, take $X=(\Ker(\mu))^0$ and $Y=\Im(\mu)$. The quotient $A/X$ is again an abelian variety and denote $\pi:A\ra A/X$ the quotient map (cf. \cite[\S 9.5]{Pol03}). The symmetric property of $\mu$ identifies $Y$ with $(A/X)^\vee$ under $\pi^\vee$ and forces the induced map $\mu':A/X\ra Y$ to be a symmetric isogeny of abelian varieties. Both $A/X$ and $Y$ have induced $\CO$-actions so that $\mu'$ is $\dag$-sesquilinear.  For an isogeny $f:A\ra B$ of abelian varieties, denote by
\[e_f:A[f]\times B^\vee[f^\vee]\ra \BG_m\]
the Weil pairing associated to $f$. 
\begin{prop}\label{muwp}
Let $\mu\in \Hom_\dag^\sym(A,A^\vee)$ be a non-zero homomorphism. 
\begin{itemize}
\item[(1)] For any $x,y\in A[\fp]$,
\[e_\fp^\mu(x,y)=e_{\mu' \varpi}(\pi(x),\pi(y))\]
where $\mu':(A/X)\ra Y$ is the induced isogeny.
\item[(2)]
\begin{itemize}
\item[(a)] If $\dag=1$,  $e_\fp^\mu$ is alternating. 
\item[(b)] If $\dag\neq 1$ and $\fp$ is ramified in $K/K_0$, $e_\fp^\mu$ is symmetric. 
\item[(c)] If $\dag\neq 1$ and $\fp$ is inert in $K/K_0$, $e_\fp^\mu$ is alternating. 
\end{itemize}
\end{itemize}
\end{prop}
\begin{proof}
(1) First we reduce to the case where $\mu$ is an isogeny. Denote the $\fp$-multiplication map as
\[g=h(i_\fp,1_A):A=\shom{\CO}{A}\lra B=\shom{\fp}{A}.\]
We have the following commutative diagram
\[\xymatrix{0\ar[r]&A[g]\ar[d]\ar[r]&A\ar[d]^{\pi}\ar[r]^g&B\ar[d]^{\pi'}\ar[r]&0\\
0\ar[r]&(A/X)[g']\ar[r]&A/X\ar[r]^{g'}&B/g(X)\ar[r]&0,}\]
with vertical quotient maps and $B/g(X)=\shom{\fp}{A/X}$. Via evaluation of sections at the generator $t_0\in \fp^{-1}/\CO$ as previously, we identify  the map $\pi'^\vee:(B/g(X))^\vee[g'^\vee]\ra B^\vee[g^\vee]$ as $\pi^\vee:(A/X)^\vee[\fp]\ra A^\vee[\fp]$.  Applying \cite[Proposition 11.20.(ii)]{Edixhoven-AV} to $h=\pi'g=g'\pi$, a diagram chasing via $\mu=\mu'\pi$ shows that for any $x,y\in A[\fp]$, $e_g^\mu(x,y)=e_{g'}^{\mu'}(\pi(x),\pi(y))$. Thus taking $\mu'$ in place of $\mu$,  without loss of generality, we may assume $\mu$ is an isogeny.

Let $m^\dag_\varpi: \CO\ra \CO$ be the $\dag$-sesquilinear map that $m^\dag_\varpi(a)=\varpi a^\dag$, and denote
\[f=\mu\varpi=h(m^\dag_\varpi, \mu):A=\shom{\CO}{A}\lra A^\vee=\shom{\CO}{A^\vee},\] 
\[k=h(m^\dag_\varpi,\mu):\shom{\fp}{A}\ra A^\vee=\shom{\CO}{A^\vee}.\] We have the following commutative diagram
\begin{equation}\label{com-1}
\xymatrix{
0\ar[r]&A[\fp]\ar[r]\ar[d]^{i}&\shom{\CO}{A}\ar[r]^{g}\ar@{=}[d]&\shom{\fp}{A}\ar[r]\ar[d]^{k}&0\\
0\ar[r]&A[f]\ar[r]&\shom{\CO}{A}\ar[r]^{f}&\shom{\CO}{A^\vee}\ar[r]&0,
}
\end{equation}
where $i$ denotes the inclusion map. By Corollary \ref{sym-cri}, the dual map $f^\vee=\mu \varpi^\dag$ and 
\[k^\vee=h(m^\dag_{\varpi^\dag}, \mu):A\ra B^\vee=\shom{\fp^{-1}}{A^\vee}.\]
Under the identification $B^\vee[g^\vee]=A^\vee[\fp]$, the induced map $k^\vee:A[f^\vee]\ra B^\vee[g^\vee]=A^\vee[\fp]$ is given by $x\mapsto \mu((\varpi t_0)^\dag x)$. Since $t_0\varpi\equiv 1\mod \fp$, for $x\in A[\fp^\dag]=A[\fp]\subset A[f^\vee]$, $k^\vee(x)=\mu(x)$. The Weil pairing $e_g$ is denoted as $e_\fp$ in (\ref{Weil-p}). By \cite[Proposition 11.21]{Edixhoven-AV}, for $x,y\in A[\fp]$,
\[e_\fp^\mu(x,y)=e_\fp(x,\mu(y))=e_g(x,k^\vee(y))=e_f(x,y).\]

\noindent (2) As is known from the proof of (1) that the pairing $e_\fp^\mu$ factors through the quotient map $\pi: A\ra A/X$, replacing $A$ by $A/X$, we may assume that $\mu$ is an isogeny. If $\dag=1$ or $\dag\neq 1$ and $\fp$ is inert in $K/K_0$, the local uniformizer $\varpi$ is chosen so that $\varpi^\dag=\varpi$ and thus $f^\vee=f$. As the commutator map of a theta group, $e_f$ is alternating (cf. \cite[Proposition 11.2 and 11.20]{Edixhoven-AV} and \cite[Proposition 10.3]{Pol03}). If $\dag\neq 1$ and $\fp$ is ramified in $K/K_0$, then $\varpi^\dag+\varpi=a\in \fp_0$, and then $f^\vee+f=\mu a$. Let $x,y\in A[\fp]\subset A[f]\cap A[f^\vee]$ be arbitrary elements. If $a=0$, then $f^\vee=-f$,  one has
\begin{equation}\label{inverse1}
e_f(x,y)e_{f^\vee}(x,y)=1.
\end{equation}
Suppose $a\neq 0$. Then $\mu a$ is an isogeny. By unraveling the definition of Weil pairings explicitly (cf. \cite[11.12]{Edixhoven-AV}), we have
\begin{equation}\label{inverse2}
1=e_{\mu a}(x,y)=e_{f+f^\vee}(x,y)=e_f(x,y)e_{f^\vee}(x,y).
\end{equation}
The first equality holds because $\ord_\fp(a)\geq 2$ and $A[\fp]$ is killed by $\fp$. To see this, one just replaces the isogeny $f$ by $\mu\circ a$ in the proof of (1), sees that $k^\vee$ vanishes on $A[\fp]$ and thus $e_{\mu a}(x,y)=1$. Combining (\ref{inverse1})-(\ref{inverse2}) with \cite[Proposition 11.17.(i)]{Edixhoven-AV}, $e_f(x,y)=e_{f}(y,x)$ as desired. 
\end{proof}

The pairing $\theta_\fp^\mu$ should be no surprisingly extended from the induced Weil pairing $e^\mu_\fp$. We explain as follows and first present a lemma.

\begin{lem}\label{trace-reduction}
Let $L$ be the maximal unramified extension of $\BQ_p$ contained in $K_\fp$, $\partial_{\CO_\fp/\BZ_p}=\fp^n$ and $p\CO_\fp=\fp^e$.
\begin{itemize}
\item[(1)] $\Tr_{K_\fp/L}(\varpi^{-n+e-1}\CO_\fp^\times)\subset \CO_{L}^\times$. 
\item[(2)] For $a\in \varpi^{-n+e-1}\CO_\fp^\times$, denote
\[\alpha=\left(\Tr_{K_\fp/L}(a) \mod p  \right)\in k^\times.\]
Then for any $x\in \CO_\fp$, 
\[\Tr_{K_\fp/\BQ_p}(ax) \mod p=\Tr_{k/\BF_p}(\alpha \overline{x})\]
where $\overline{x}$ denotes the reduction of $x$ in $k$.
\end{itemize}
\end{lem}
\begin{proof}
In order to prove (1), it suffices to prove that, for any unit $u\in \CO_\fp^\times$,  $\Tr_{K_\fp/L}\left(\frac{u}{\varpi^{n+1}}\right)$ has $p$-adic valuation $-1$. Since $L/\BQ_p$ is unramified, $\partial_{\CO_\fp/\CO_{L}}=\partial_{\CO_\fp/\BZ_p}$. By  \cite[Proposition 4, \S VIII]{Weil-ANT},
\[\Tr_{K_\fp/L}(\fp^{-n-1})=p^{-1}\CO_{L},\quad \Tr_{K_\fp/L}(\fp^{-n})=\CO_{L}.\]
Thus there exists a unit $u_0\in \CO_\fp^\times$ such that $\Tr_{K_\fp/L}\left(\frac{u_0}{\varpi^{n+1}}\right)$ has $p$-adic valuation $-1$. Let  $u\in \CO_\fp^\times$ be an arbitrary unit and suppose 
\[u\equiv a\mod \fp,\quad u_0\equiv a_0\mod \fp\]
where $a,a_0\in \CO_{L}^\times$. Then
\[u=aa_0^{-1}u_0+\varpi t, \quad t\in \CO_\fp,\]
and 
\[\Tr_{K_\fp/L}\left(\frac{u}{\varpi^{n+1}}\right)=aa_0^{-1}\Tr_{K_\fp/L}\left(\frac{u_0}{\varpi^{n+1}}\right)+\Tr_{K_\fp/L}\left(\frac{t}{\varpi^{n}}\right),\]
which has $p$-adic valuation $-1$. This proves (1).

Let $L'/\BQ_p$ be the Galois closure of $K_\fp/\BQ_p$ and $G=\Gal(L'/\BQ_p)$. Let $H_1 \text{ resp. }H_2\subset G$ be the subgroup fixing $K_\fp$ resp. $L$. Let $x_0\in \CO_{L}$ be a lifting of $\overline{x}$ under reduction modulo $\fp$. Then, noting $\Tr_{K_\fp/\BQ_p}((x-x_0)a)\in p\BZ_p$, one has 
\begin{eqnarray*}
\Tr_{K_\fp/\BQ_p}(ax)\mod p&=&\sum_{\sigma\in H_1\backslash G} x_0^\sigma a^\sigma \mod p\\
&=&\sum_{\sigma\in H_2\backslash G} x_0^\sigma \left(\sum_{\tau\in H_1\backslash H_2} a^\tau\right)^\sigma\mod p\\
&=&\sum_{\sigma\in \Gal(L/\BQ_p)} x_0^\sigma \left(\Tr_{K_\fp/L}(a)\right)^\sigma\mod p\\
&=& \Tr_{k/\BF_p} (\alpha \overline{x}).
\end{eqnarray*}
as desired.
\end{proof}
\begin{prop}\label{compare-weil}
Denote
\[\alpha_0=\left(\Tr_{K_\fp/L}(dp/\varpi^\dag)\mod p\right) \in k^\times.\]
For $x,y\in A[\fp]$, we have
\[e_\fp^\mu(x,y)=\Tr_{k/\BF_p}(\alpha_0 \theta_\fp^\mu(x,y)).\]
\end{prop}
\begin{proof}
Let $h:\shom{\fp}{A}\ra\shom{p\CO}{A}$ be the isogney induced by the inclusion $p\CO\ra \fp$ and the identity $1_A$. The composition 
\[h(i_{p\CO},1_A):A=\shom{\CO}{A}\xrightarrow{h(i_\fp,1_A)}\shom{\fp}{A}\xrightarrow{h} \shom{p\CO}{A}\]
is essentially the $p$-multiplication morphism on $A$, whose kernel $\shom{p^{-1}\CO/\CO}{A}$, by associating a section with its image at $p^{-1}$, is identified with $A[p]$. The Cartier duality induces the Weil pairing
\[e_p:A[p]\times A^\vee[p]\lra \mu_p.\]
By Proposition \ref{dual}, the dual map $h^\vee:\shom{p^{-1}\CO}{A^\vee}\ra\shom{\fp^{-1}}{A^\vee}$ is induced by the inclusion $\fp^{-1}\ra p^{-1}\CO$ and the identity $1_{A^\vee}$. The restricted morphism $h^\vee:A^\vee[p]\ra A^\vee[\fp]$ is the multiplication by $pt_0$. Moreover, when restricted to $A^\vee[\fp^e]\subset A^\vee[p]$, it is even given by the multiplication by $p/\varpi$. It follows from \cite[Proposition 11.21]{Edixhoven-AV} that for any $x\in A[\fp]$ and $y\in A^\vee[\fp^e]$ we have
\[e_\fp(x,(p/\varpi)y)=e_\fp(x,h^\vee(y))=e_p(x,y).\]
Let $x,y\in A[\fp]$ be arbitrary elements. Since the map $A[\fp^m]\ra A[\fp^{m-1}]$ induced by the local uniformizer $\varpi$ is surjective, one can choose $u\in A[\fp^e]$ satisfying $(p/\varpi^\dag)u=y$. Let 
\begin{eqnarray*}
&\wt{x}=(\cdots, x_2,x), \quad \wt{u}=(\cdots, u_2,u)\text{ and }\\
&\wt{x}'=(\cdots, \varpi x_3, \varpi x_2) \in T_\fp A=\lim\limits_{\begin{subarray}{c}\lla\\ n\end{subarray}} A[\fp^{e n}].
\end{eqnarray*}
Then  $(p/\varpi)\wt{x}'=\wt{x}$ and
\[e_\fp^\mu(x,y)=e_p^\mu(x,u)=E_p^\mu(\wt{x},\wt{u})\mod p=E_{\fp_0}^\mu(\wt{x},\wt{u})\mod p.\]
From (\ref{characterization}), 
\begin{eqnarray}\label{weil-trace}
e_\fp^\mu(x,y)&=&\Tr_{K_{\fp_0}/\BQ_p} (d\Theta_{\fp_0}^\mu(\wt{x},\wt{u}))\mod p\nonumber\\
&=&\Tr_{K_{\fp_0}/\BQ_p} ((dp/\varpi^\dag)\Theta_{\fp_0}^\mu(\wt{x}',(\varpi/\varpi^\dag)\wt{u}))\mod p.
\end{eqnarray}
Let $L\subset K_\fp$ be the maximal unramified extension of $\BQ_p$, and let $b\in \CO_L$ be an integer satisfying $$b \equiv \Theta_{\fp_0}^\mu(\wt{x}',(\varpi/\varpi^\dag)\wt{u}) \mod \fp.$$ 
Under the identification $T_\fp A/\fp T_\fp A\cong A[\fp]$, $\wt{x}'$ resp. $(\varpi/\varpi^\dag)\wt{u}$ is mapped to $x$ resp. $y$. Thus we have 
\[(b\mod \fp)=\theta_\fp^\mu(x,y)\in k.\]
It follows from (\ref{weil-trace}) and Lemma \ref{trace-reduction} that
\[e_\fp^\mu(x,y)=\Tr_{k/\BF_p}(\alpha_0 \theta_\fp^\mu(x,y))\]
as desired.
\end{proof}

\subsection{Quadratic refinements of of induced Weil pairings}
 In this subsection, we assume $p=2$.  Let $\sU$ be the functor, constructed in \cite[\S 2]{PR11}, from the category of abelian sheaves to the category of sheaves of groups over the big fppf site over $F$. The abelian sheaf $(\sU A[\fp])^\ab$ classifies quadratic maps from $A[\fp]$ to abelian sheaves (cf. \cite[Proposition 2.3]{PR11} and its sheaf analogue), and sits into an exact sequence
 \begin{equation}\label{symmetric-exact}
 0\ra S^2A[\fp]\lra (\sU A[\fp])^\ab\lra A[\fp]\ra 0.
 \end{equation}
  Applying $\bhom(\cdot, \BG_m)$ to the exact sequence and using ${\bf Ext^1}(A[\fp],\BG_m)=0$ (\cite[Theorem 1]{Waterhouse71}), we obtain 
\begin{equation}\label{Cartier-dual}
0\ra\bhom(A[\fp],\BG_m)\ra\bhom((\sU A[\fp])^\ab,\BG_m)\ra \bhom(S^2A[\fp],\BG_m)\ra 0.
\end{equation}
Let $\bhom_\dag^\sym(A,A^\vee)$ be the sheaf of $\dag$-sesquilinear symmetric homomorphisms. Since $p=2$, by Proposition \ref{muwp}, for any section $ \mu\in \bhom_\dag^\sym(A,A^\vee)$, the pairing $e_\fp^\mu$ is symmetric. Pulling back through $\bhom_\dag^\sym(A,A^\vee)\ra \bhom(S^2 A[\fp],\BG_m), \mu \mapsto e_\fp^\mu$, gives rise to the following commutative diagram
\begin{equation}\label{pullback}
\xymatrix@C=9pt{
0\ar[r]&A^\vee[\fp]\ar[r]\ar[d]^{e_\fp}_{\wr}&\fX\ar[r]\ar[d]^{\fq}& \bhom_\dag^\sym({A,A^\vee})\ar[r]\ar[d]^{e_\fp^\bullet}&0\\
0\ar[r]&\bhom(A[\fp],\BG_m)\ar[r]&\bhom((\sU A[\fp])^\ab,\BG_m)\ar[r]&\bhom(S^2A[\fp],\BG_m)\ar[r]&0.
}
\end{equation}
The left vertical isomorphism is the Cartier duality. The cohomology of the top row of the above diagram gives a connecting homomorphism
\[\Hom_\dag^\sym(A,A^\vee)\ra \RH^1(F,A^\vee[\fp]),\quad \mu\mapsto c_\mu.\]
The above construction works with any $F$-scheme as the base. 
\begin{remark}
If we take $\CO=\BZ$, $\fp=(2)$, then $\bhom_\dag^\sym(A,A^\vee)=\bhom^\sym(A,A^\vee)$ and the pull-back $\fX$ is the symmetric Picard scheme $\Pic_{A/F}^\sym$. Then $c_\mu$ is the obstruction to find a quadratic refinment of the (alternating) pairing $e_2^\mu$ over $F$, or equivalently, a rational symmetric line bundle $L$ with $\varphi_L=\mu$ (cf. \cite[Proposition 3.6]{PR11} and \cite[Proposition 13.1]{Pol03}).
\end{remark}

\begin{defn}
Let $S$ be an $F$-scheme and $\mu\in \Hom_\dag^\sym(A_S,A_S^\vee)$. A morphism $q:A[\fp]_S\ra {\BG_m}_S$ of schemes is  called a quadratic refinement of the pairing $e_\fp^\mu$ if $q(x+y)=q(x)q(y)e_\fp^\mu(x,y)$ for all $x,y\in A[\fp]_S$.
\end{defn}
\begin{prop}\label{self-cup}
Let $S$ be an $F$-scheme, $\mu\in \Hom_\dag^\sym(A_S,A_S^\vee)$. 
\begin{itemize}
\item[(1)] The pairing $e_\fp^\mu$ has a quadratic refinement over $S$ if and only if $c_{\mu}=0$.
\item[(2)]For any  $x\in \RH^1(S,A[\fp])$, 
\[x\mathop{\cup}_{e_\fp^\mu}x=x\mathop{\cup}_{e_\fp}c_{\mu}\]
in $\RH^2(S,\BG_m)$.
\end{itemize}
\end{prop}
\begin{proof}
(1) Change the base to $S$, and consider the derived diagram from (\ref{pullback}):
\[\xymatrix@C=10pt{&\Hom_\dag^\sym(A_S,A^\vee_S)\ar[r]\ar[d]&\RH^1(S,A^\vee[\fp])\ar[d]^{\wr}\\
\Hom((\sU A[\fp])^\ab_S,\BG_{m,S})\ar[r]& \Hom(S^2 A[\fp]_S,\BG_{m,S})\ar[r]^\delta&\RH^1(S,\bhom(A[\fp],\BG_m)).}\]
Then $c_{\mu}=0$ just means $e^{\mu}_\fp\in \ker(\delta)$, i.e.  it comes from a quadratic refinement $A[\fp]_S\ra \BG_{m,S}$.

\noindent (2) The statement follows by applying \cite[Corollary 2.16.(a),(d)]{PR11} with $M=A[\fp]$, $N=\BG_m$ and $\beta=e_\fp^\mu$  (cf. \cite[Theorem 3.9]{PR11}).

\end{proof}
We have the following vanishing criterion for $c_\mu$ generalizing \cite[Proposition 3.12.(a)]{PR11}.
\begin{prop}\label{vanishing-criterion}
Suppose $\ell\neq 2$. Let $E$ be a field extension of $F$ with $E^s$ as its separable closure, and $\mu\in \Hom_\dag^\sym(A_E,A_E^\vee)$. If the image $G$ of $\Gal(E^s/E)\ra \Aut(A[\fp](E^s))$ is cyclic, then $c_\mu=0$.
\end{prop}
\begin{proof}
For simplicity, we drop the subscript $E$ which is used to indicate the base. View the sequence (\ref{Cartier-dual}) as an exact sequence of finite $\BZ[G]$-modules, and the Pontryagin dual exact sequence is (\ref{symmetric-exact}). Since $(\sU A[\fp])^\ab\ra A[\fp]$ has a section as $G$-sets (cf. \cite[\S 2]{PR11}), the sequence 
\[0\ra (S^2A[\fp])^G\lra ((\sU A[\fp])^\ab)^G\lra (A[\fp])^G\ra 0\]
is exact. By Lemma \ref{invariants-dual} and counting cardinality, 
\[0\ra(\bhom(A[\fp],\BG_m))^G\ra(\bhom((\sU A[\fp])^\ab,\BG_m))^G\ra (\bhom(S^2A[\fp],\BG_m))^G\ra 0\]
is also exact, which implies the connecting map $(\bhom(S^2A[\fp],\BG_m))^G\ra \RH^1(E, A^\vee[\fp])$ is zero.  
\end{proof}
Let $M$ be a finite abelian group and $g:M\ra M$ is a group homomorphism.  Denote $M'=\Hom(M,\BC^\times)$ the Pontryagin dual of $M$ and denote $g':M'\ra M'$ the homomorphism induced by composing $g$. 
\begin{lem}\label{invariants-dual}
$|\Ker(g-1)|=|\Ker(g'-1)|.$
\end{lem}
\begin{proof}
 Consider the perfect pairing $M\times M'\ra \BC^\times, (m,f)\mapsto f(m)$. Then $g'$ is the adjoint map of $g$ under this pairing. Thus $\Ker(g'-1)=((g-1)M)^\perp\cong M/(g-1)M$.
From the exact sequence
\[0\ra \Ker(g-1)\ra M\xrightarrow{g-1}M\ra M/(g-1)M\ra 0,\]
we have $|\Ker(g-1)|=|M/(g-1)M|=\Ker(g'-1)$ as desired.  
\end{proof}

The multiplication of $K$ on $A$ has a strong impact on the vanishing of $c_\mu$ as follows.
\begin{prop}\label{global-criterion}
Suppose $\ell\neq 2$ and $\mu\in \Hom_\dag^\sym(A,A^\vee)$. If one of the following conditions holds:
\begin{itemize}
\item[(1)] $\dag=1$ and $[K:\BQ]=\dim A$,
\item[(2)] $\dag\neq 1$ and $[K:\BQ]=2\dim A$,
\item[(3)]$\dag\neq 1$ and the prime $\fp$ is inert in $K/K_0$,
\item[(4)] $\dag= 1$ and there exists a quadratic extension $K'/K$ contained in  $\End^0(A)$ together with a nontrivial involution $\dag'$ on $\CO'=K'\cap \End_F(A)$ extending $\dag$ such that $\mu$ is $\dag'$-sesquilinear, $2\nmid [\CO_{K'}:\CO']$ and $\fp$ is either inert or split in $K'/K$,
\end{itemize}
then $c_\mu=0$.
\end{prop}
\begin{proof}
Since $\ell \neq 2$, $A[\fp]$ is a finite \'etale group scheme over $F$ and we view it as a $G_F$-module.  

\noindent (1) Suppose $\dag=1$ and $[K:\BQ]=\dim A$. Since $p=2$ is prime to conductor of $\CO$, by Proposition \cite[Proposition 2.2.1-2]{Ribet76b}, $A[\fp]$ is a $k$-vector space of dimension $2$.  If $\theta_\fp^\mu$ is trivial, by Proposition \ref{compare-weil}, $e_\fp^\mu$ is trivial and $c_\mu=0$. Assume $\theta_\fp^\mu$ is not trivial. We may choose a basis $e_1,e_2$ of $A[\fp]$ such that $\theta_\fp^\mu(e_1,e_2)=1\otimes -1\in k(1)$. Then for any $a,b,c,d\in k$,
\[\theta_\fp^\mu(ae_1+be_2,ce_1+de_2)=(ad+bc)\otimes -1.\]
Let $q:A[\fp]\ra k(1)$ be the map defined on all $F^s$-points of $A[\fp]$ by $$q(ae_1+be_2)=ab\otimes -1.$$
As a morphism of group schemes, $\theta_{\fp}^\mu$ is Galois equivariant. Then $q(ae_1+be_2)=\theta_\fp^\mu(ae_1,be_2)$ is Galois equivariant. It is readily verified that $q$ is a quadratic refinement of $\theta_\fp^\mu$. Let $\alpha_0\in k$ be given as in Proposition \ref{compare-weil}. Then $q'=\tr_{k/\BF_p}(\alpha_0 q)$ is a quadratic refinement of $e^\mu_\fp$ over $F$.

\noindent (2) In this case $A[\fp]$ is a $k$-vector space of dimension $1$. Since $G_F$ lands in $\Aut_\CO(A[\fp])\cong k^\times$, the image is cyclic. It follows from Proposition \ref{vanishing-criterion} that $c_\mu=0$.

\noindent (3)-(4) The situation (3) is contained in (4), and (4) follows from the linear aspect of the presence of the nontrivial involution. We consider the situation where $A$ has multiplication by $K'$ and then $K'_0=K$, $\fp'_0=\fp$. If we denote $k'=\CO'/\fp_0'\CO'$ and $k'_0=\CO/\fp$, then $\dag'$ induces a nontrivial involution of $k'$ with $k'_0$ as its fixed field. There exists an element $\alpha\in k'$ such that $\alpha+\alpha^\dag=1$. By Proposition \ref{extended-Weil}, reducing pairing $\Theta_{\fp_0'}^\mu$ modulo $\fp_0'$, we get a $(-1)$-hermitian pairing 
\[\theta_{\fp_0'}^\mu:A[\fp_0']\times A[\fp_0']\lra k'(1).\]
In particular, $\theta_{\fp_0'}^\mu=\alpha \theta_{\fp_0'}^\mu+ \left(\alpha \theta_{\fp_0'}^\mu\right)^*$.
If $\fp_0'$ is inert in $K'/K'_0$, take $\alpha_0\in {k'}^\times$ to be the element of  Proposition \ref{compare-weil} which actually lies in ${k_0'}^\times$; But in the split case, take $\alpha_0\in {k'_0}^\times$ of  same proposition applied to the situation where $A$ has multiplication by $K=K_0$ and $\dag=1$.  For any $x\in A[\fp]$, take $q(x)=\Tr_{k'/\BF_p}(\alpha_0\alpha \theta_{\fp_0'}^\mu(x,x))$. Since $\theta_{\fp'_0}^\mu$ is Galois equivariant, so is $q$. It is direct to verify that $q$ is a quadratic refinement of $e_{\fp}^\mu$ over $F$.
\end{proof}

\begin{thm}\label{Rational-condition}
Suppose $\ell\neq 2$ and $2$ is prime to the discriminant of $\CO$. Suppose either
\begin{itemize} 
\item[(1)] $\dag=1$ and $[K:\BQ]=\dim A$, \text{ or }
\item[(2)]  $\dag\neq 1$.
\end{itemize}
For any $\mu\in \Hom^\sym_\dag(A,A^\vee)$, there exists a rational symmetric line bundle $L$ on $A$ such that $\varphi_L=\mu$. 
\end{thm}
\begin{proof}
We decompose $2\CO=\fp_1\cdots\fp_r\fq_1\fq_1^\dag\cdots \fq_s\fq_s^\dag$ into a product of distinct normal prime ideals of $\CO$ with $\fp_i=\fp_i^\dag$.  If $\dag=1$, there are no such primes $\fq_i$'s. Through the decomposition 
\[A[2]=\left(\bigoplus_i A[\fp_i]\right)\bigoplus \left(\bigoplus_j A[\fq_j]\bigoplus A[\fq_j^\dag]\right),\] 
$e_2^\mu=\sum_i e_{\fp_i}^\mu+\sum_j (e_{\fq_j}^\mu+e^\mu_{\fq_j^\dag})$. By Proposition \ref{global-criterion} and Proposition \ref{self-cup}, $e_{\fp_i}^\mu$ has a quadratic refinement over $F$. Since $p=2$, by \cite[Proposition 11.21.(1)]{Edixhoven-AV}, $e^\mu_{\fq_j}(x,y)=e^\mu_{\fq_j^\dag}(y,x)$. Then $e^\mu_{\fq_j}+e^\mu_{\fq_j^\dag}$ has the quadratic refinement $(x,y)\mapsto e_{\fq_j}^\mu(x,y)$  (This gives another proof of Proposition \ref{global-criterion}.(4) in the split case). Thus $e_2^\mu$ has a quadratic refinement over $F$ and the theorem follows by \cite[Proposition 3.6]{PR11}.
\end{proof}
\begin{remark}
The above theorem states that, for an abelian variety $A$ with an embedding $K\hookrightarrow \End^0(A)$ and $\CO=K\cap \End^0(A)$ of odd discriminant, being of $\GL_2$-type, i.e. $[K:\BQ]=\dim A$ or the presence of a nontrivial involution $\dag$ on $K$ breaks the obstruction of a $\dag$-sesquilinear symmetric homomorphism being induced from a rational symmetric line bundle. This will have an effect on the linear structure of the Shafarevich-Tate groups as discussed in the final section.
\end{remark}

\section{Quadratic abelian varieties}\label{QC}

From now on we fix a symmetric isogeny $\lambda:A\ra A^\vee$ and denote by $\dag$ the associated Rosati involution on $\End^0(A)$, i.e. $\phi^\dag=\lambda^{-1}\phi^\vee\lambda$. Let $K$ be a number field contained in $\End^0(A)$ and suppose the order $\CO= K\cap \End_F(A)$ is stable under $\dag$. Thus the Rosati involution restricts an involution of $K$. If $\lambda$ is a polarization, then the Rosati involution is positive, and by \cite[Proposition 1.39]{Milne-CMAV}, $K$ is either a totally real field or a CM field, and the involution $\dag$ restricts to the unique complex conjugation on $K$. 

Let $p$ be a rational prime not diving $[\CO_K:\CO]\deg(\lambda)$. Let $\fp$ be an invertible prime ideal of $\CO$ above $p$ and $\fp_0=\fp\cap \CO_0$. In this section we assume $\fp^\dag=\fp$.
\subsection{Quadratic abelian varieties}
Let $\Omega$ be the set of places of $F$. For $v\in \Omega$, let $F_v$ denote the completion of $F$ at $v$. Fix the separable closures $F^s$ and $F_v^s$ of $F$ and $F_v$ respectively. Denote the corresponding Galois group as $G_F=\Gal(F^s/F)$ and $G_v=\Gal(F_v^s/F_v)$. For each $v\in \Omega$, the local Tate cup-product pairing 
\begin{equation}\label{local-Tate}
\bigcup_{e_\fp^\lambda}:\RH^1(F_v,A[\fp])\times\RH^1(F_v,A[\fp])\lra \RH^2(F_v, \BG_m)
\end{equation}
is $\BF_p$-bilinear. In this subsection, we discuss when it is an even $\BF_p$-bilinear form so that one can put quadratic structures on $\RH^1(F_v,A[\fp])$ (cf. Remark \ref{quadratic-evenhermitian2}). For $p\neq 2$, the local pairing $\bigcup_{e_\fp^\lambda}$ is always even by Lemma \ref{auto-quadratic}.  As a result, we restrict to the case  $p=2$. 

For each $v\in \Omega$, let $c_{\lambda,v}\in \RH^1(F_v,A^\vee[\fp])$ be the localization of $c_\lambda$ at $v$.
\begin{defn}
The pair  $(A,\lambda)$, with additional data $K$ and $\fp$, is called quadratic at $v$ if $2\nmid \deg(\lambda)$ and $c_{\lambda,v}=0$.
\end{defn}

\begin{prop}\label{quadratic-alternating}
Suppose $2\nmid \deg(\lambda)$ and let $v$ be a place of $F$. The following are equivalent:
\begin{itemize}
\item[(1)] The pair  $(A,\lambda)$ is quadratic at $v$.
\item[(2)] The pairing $e_\fp^\lambda$ has a quadratic refinement over $F_v$.
\item[(3)] The local Tate pairing $(\ref{local-Tate})$ is alternating. 
\end{itemize}
\end{prop}
\begin{proof}
The equivalence of (1) and (2) is just Proposition \ref{self-cup}.(1) with $S=\Spec F_v$. Applying Proposition \ref{self-cup}.(2), for any $x\in \RH^1(F_v,A[\fp])$, $$x\mathop{\cup}_{e_\fp^\lambda}x=x\mathop{\cup}_{e_\fp}c_{\lambda,v}.$$ 
Since $2\nmid \deg(\lambda)$, by local Tate duality, $\cup_{e_\fp}^\lambda$ is nondegenerate. Thus  $\cup_{e_\fp^\lambda}$ is alternating if and only if $c_{\lambda,v}=0$ and the equivalence of $(1)$ and $(3)$  follows.
\end{proof}

Thus $(A,\lambda)$ is quadratic at $v$ is equivalent to that the local Tate pairing $\cup_{e_\fp^\lambda}$ is an even $\BF_p$-bilinear form (cf. Proposition \ref{even-alternating}). For an example of non-quadratic $(A,\lambda)$, see \cite[Example 3.20.(a)]{PR11}. We have similar local vanishing criterions for $c_{\lambda,v}$ as in \cite[Proposition 3.12]{PR11}.
\begin{prop}\label{local-criterion}
The pair $(A,\lambda)$ is quadratic at $v$, if $2\nmid \deg(\lambda)$ and one of the following conditions holds:
\begin{itemize}
\item[(1)] $\ell\neq 2$ and the image of $G_v$ in $\Aut(E[\fp])$ is cyclic.
\item[(2)] $v=\BR$ or $\BC$.
\item[(3)] The residue characteristic of $F_v$ is not $2$ and $A$ has good reduction at $v$.
\end{itemize}
\end{prop}
\begin{proof}
All follows from Proposition \ref{vanishing-criterion}.  As for (3), by the criterion of N\'eron-Ogg-Shafarevich and its generalization by Serre and Tate \cite[Theorem 1]{Serre-Tate68}, $A[\fp]$ is unramified and the image of $G_v$ is generated by the image of the Frobenius map.
\end{proof}

We also have the following everywhere results.
\begin{prop}\label{global-criterion1}
The pair $(A,\lambda)$ is quadratic everywhere, if $2\nmid \deg(\lambda)$ and one of the following conditions holds:
\begin{itemize}
\item[(1)] $\ell\neq 2$, $\dag|_K=1$ and $[K:\BQ]=\dim A$.
\item[(2)] $\ell\neq 2$, $\dag|_K\neq 1$ and $[K:\BQ]=2\dim A$.
\item[(3)]$\dag|_K\neq 1$ and the prime $\fp$ is inert in $K/K_0$.
\item[(4)] $\dag|_K= 1$ and there exists a quadratic extension $K'/K$ contained in  $\End^0(A)$ such that $\CO'=K'\cap \End_F(A)$ is stable under $\dag$, $\dag|_{K'}\neq 1$, $2\nmid [\CO_{K'}:\CO']$ and $\fp$ is either inert or split in $K'/K$.
\end{itemize}
\end{prop}
\begin{proof}
If $\ell\neq 2$, the proposition follows by applying Proposition \ref{global-criterion} to $\mu=\lambda$.  The case (3) is again contained in (4) and in order to prove (4), we work with local Tate pairings rather than Weil pairings. Let $\alpha$ and $\alpha_0$ be as chosen in the proof of Proposition \ref{global-criterion}.(4). For any $v\in \Omega$,  the local Tate pairing $h_v$ of (\ref{extended-pairing}) and (\ref{extended-pairing2}) is $1$-hermitian and thus $h_v=\alpha h_v +(\alpha h_v)^*$ is automatically even. If we set $s=\Tr_{k'/\BF_p}(\alpha\alpha_0 h_v)$, then $\cup_{e_{\fp}^\lambda}=s+s^*$ is even and therefore alternating. By Proposition \ref{quadratic-alternating}, $(A,\lambda)$ is quadratic at $v$.

\end{proof}

\subsection{Quadratic maps arising from line bundles}
In this and the next subsection, the goal is to provide quadratic maps for the local Tate pairings $\cup_{e_\fp^\lambda}$ in the ($\Ort$) case of Section \ref{QSS}. In both subsections, we assume $\dag|_K=1$ and $p=2$. Let $\Pic_{A/F}$ be the relative Picard group scheme. For any $F$-scheme $S$, $\Pic_{A/F}(S)=\Pic(A_S)/\mathrm{pr}_S^*\Pic(S)$ where $\Pic(A_S)$  denotes the group of isomorphism classes of line bundles on $A_S=A\times_F S$. Given a line bundle $L$ on $A$, the theorem of square implies the map
\[\varphi_L:A\lra \Pic_{A/F},\quad x\mapsto t_x^*L\otimes L^{-1}\]
is a homomorphism with image contained in $A^\vee=\Pic_{A/F}^0$. Thus we may view $\varphi_L$ as a homomorphism $A\ra A^\vee$ (cf. \cite[p.59, Corollary 4 and p.131, Corollary 5]{Mumford-AV}). Moreover, $\varphi_L$ is a symmetric homomorphism.

As explained in Section \ref{IWP} or Appendix \ref{Hom-construction}, the $\fp$-multiplication map
\[g=h(i_\fp,1_A):A=\shom{\CO}{A}\lra \shom{\fp}{A}\]
sits in the short exact sequence
\[0\ra A[\fp]\lra A\xlongrightarrow{g}\shom{\fp}{A}\ra 0.\]
Over any $F$-scheme $S$,  the cohomology induces a connecting map, i.e.  the Kummer $\fp$-descent map over $S$
\begin{equation}\label{kummer-map}
\delta_g:\Hom_\CO(\fp,A(S))\lra \RH^1(S,A[\fp]).
\end{equation}

Recall $m^\dag_\varpi: \CO\ra \CO$ is the map that $m^\dag_\varpi(a)=\varpi a^\dag=\varpi a$. Denote
\[f=\lambda \varpi=h(m^\dag_\varpi, \lambda):A=\shom{\CO}{A}\lra A^\vee=\shom{\CO}{A^\vee}.\] 
By Corollary \ref{sym-cri}, the isogeny $f=\lambda \varpi$ is symmetric. By \cite[Proposition 11.2]{Edixhoven-AV}, there exists a finite extension $E$ of $F$ and a nondegenerate line bundle $L$ on $A_E$ such that $\varphi_L=f$.

Let $S$ be an $F$-scheme and suppose there is a line bundle $L$ on $A_S$ such that $\varphi_L=f$. Let $\CG(L)$ be the theta group over $S$ associated to the line bundle $L$ (cf. \cite[p. 225, Theorem 1]{Mumford-AV}, \cite[pp. 44-46]{Mumford-Tata3}). There is a central extension over $S$
\[1\ra \BG_m\lra \CG(L)\lra A[f]\ra1\]
with the Weil pairing $e_f$ as its commutator map.  Taking cohomology gives rise to a connecting map $q_{f}:\RH^1(S,A[f])\ra \RH^2(S,\BG_m)$, and denote 
\[q_{\fp}=-q_{f}\circ i:\RH^1(S,A[\fp])\xrightarrow{i}\RH^1(S,A[f])\xrightarrow{-q_{f}} \RH^2(S,\BG_m),\]
where $i:A[\fp]\ra A[f]$ denotes the inclusion map.
\begin{prop}\label{RM-quadratic-map}
Let $S$ be an $F$-scheme. Let $L$ be a line bundle on $A_S$ such that  $f=\varphi_L$ and $q_\fp$ is the map constructed from $L$ as above.
\item[(1)] For any $x,y\in \RH^1(S,A[\fp])$,
\[x\cup_{e_\fp^\lambda}y=q_\fp(x+y)-q_\fp(x)-q_\fp(y).\]
\item[(2)] Let $\delta_g$ be the Kummer $\fp$-descent map over $S$ in $(\ref{kummer-map})$. Then $q_\fp \delta_g=0$.
\end{prop}
\begin{proof}
(1) Applying Proposition \ref{muwp}.(1) with $\mu=\lambda$, we obtain that for $x,y\in A[\fp]$, $e_\fp^\lambda(x,y)=e_f(x,y).$ Thus we have the following commutative diagram
\[\xymatrix@C=6pt{
\RH^1(S, A[\fp])&\times \ar[d]^{i\times i}&\RH^1(S,A[\fp])\ar[rr]^-{\cup_{e_\fp^\lambda}}&&\RH^2(S, \BG_m) \ar@{=}[d]\\
\RH^1(S,A[f])&\times&\RH^1(S, A[f])\ar[rr]^-{\mathop{\cup}_{e_f}}&&\RH^2(S,\BG_m).
}
\]
 Applying \cite[Proposition 2.17]{PR11} to $f$, the above diagram shows that for $x,y\in \RH^1(S,A[\fp])$
\[x\mathop{\cup}_{e_\fp^\lambda}y=q_{\fp}(x+y)-q_{\fp}(x)-q_{\fp}(y).\]

\noindent (2) The commutative diagram (\ref{com-1}) (with $\mu=\lambda$) induces a commutative diagram
\[\xymatrix{
\Hom_\CO(\fp,A(S))\ar[d]^{k}\ar[r]^{\delta_g}&\RH^1(S,A[\fp])\ar[d]^i\\
A^\vee(S)\ar[r]^{\delta_f}&\RH^1(S,A[f])
}\]
where $\delta_g$ and $\delta_f$ are the connecting morphisms. 
By \cite[Proposition 4.9 and Remark 4.10]{PR12} , we have $q_f \delta_f=0$, and thus $q_\fp\delta_g=0$ as desired.
\end{proof}

\begin{remark}\label{local-quadratic}
\begin{itemize}
\item[(1)] Let $v\in \Omega$ be a place of $F$ and $S=\Spec F_v$.  A result \cite[Lemma 1, p. 1119]{PS99} of Poonen-Stoll tells that there always exists a line bundle $L$ over $F_v$ such that $\varphi_L=f$. The quadratic map $q_\fp$ generally takes values in $4$-torsions of $\RH^2(F_v,\BG_m)\cong \BQ/\BZ$ and it takes $2$-torsions as values if and only if $c_{\lambda,v}=0$, i.e. $(A,\lambda)$ is quadratic at $v$.

\item[(2)] The construction of Proposition \ref{RM-quadratic-map} also works for $p\neq 2$. But one needs to take a symmetric $L$, i.e.  $[-1]^*L\cong L$, to ensure that $q_\fp(ax)=a^2q_\fp(x)$ for all $a\in \BF_p$ (cf. the proof of \cite[Corollary 4.7]{PR12}). If this is the case, no wonder that, $q_\fp(x)=\frac{1}{2}(x\cup_{e^\lambda_\fp} x)$.
\end{itemize}
\end{remark}

\subsection{A global condition}\label{global-hyp}
In order to construct local quadratic maps that vanish on Kummer images and interplay with global  cohomologies (for the $(\Ort)$ case of the next section), we impose the following global condition:

\[
\left\{\begin{aligned}
&\text{The symmetric isogeny $\lambda=\varphi_M$ for some rational symmetric} &&\\  
&\text{ line bundle $M$ on $A/F$ and $2$ is unramified in $\CO$.}&&
 \end{aligned}\right.
 \leqno{\mathrm{(\bf{L}):}}
\]
A line bundle $L$ is symmetric if $[-1]^*L\cong L$. The symmetric Picard group scheme, denoted by $\Pic^\sym_{A/F}$, is the closed subgroup scheme of $\Pic_{A/F}$ that parametrizes symmetric line bundles on $A$.
\begin{prop}\label{rational-line-bundle}
Supposethe hypothesis $({\bf L})$ holds. There exists a unique group homomorphism $L: \CO\ra \Pic^\sym_{A/F}(F)$  such that $L(1)=M$, $\varphi_{L(a)}=\lambda a$ and $b^*L(a)=L(ab^2)$ for any $a,b\in \CO$.
\end{prop}
\begin{proof}
This is proved by Polishchuk (cf. \cite[Proposition 1.5.1]{Pol-CM}) and we describe the construction as follows.
Since $2$ is unramifed in $\CO$, $\CO/2\CO$ is a direct sum of finite fields. The Frobenius map $x\mapsto x^2$ is an automorphism of $\CO/2\CO$. For any $a\in \CO$, $a=x^2+2y$ for some $x,y\in \CO$ with $x$ uniquely chosen modulo $2$. Let $\sP$ be the normalized Poincar\'e line bundle on $A\times A^\vee$ and take $$L(a)=x^*M\otimes (y, \lambda)^*\sP.$$ It is straightly verified that the map $L$ is  well-defined with required property.
\end{proof}
\begin{remark}
The assumption that $2$ is unramified in $\CO$ can be refined to any weaker condition under which if one can show that the fppf-sheaf version of the extension \cite[(1.5.2)]{Pol-CM} is split. In fact,  \cite[Theorem 1.5.6]{Pol-CM} proves this sheaf extension is locally split.
\end{remark}

\begin{prop}\label{additional}
\begin{itemize}
\item[(1)] Suppose the hypothesis $({\bf L})$ holds. Then $c_\lambda=0$. In particular, $(A,\lambda)$ is quadratic everywhere. 
\item[(2)] Suppose $\ell\neq 2$ and $2$ is unramified in $\CO$. If $\dag|_K=1$ and $[K:\BQ]=\dim A$, then  the hypothesis $({\bf L})$ holds.
\end{itemize}
\end{prop}
\begin{proof}
(1) By \cite[Proposition 3.6]{PR11}, the existence of the rational symmetric line bundle $M$ with $\lambda=\varphi_M$ implies that $e_2^\lambda$ has a quadratic refinement over $F$. Since $2$ is unramified in $\CO$, taking $\varpi=2$ to be the local uniformizer for $\fp$, Proposition \ref{muwp}.(1) tells $e_\fp^\lambda=e_{2\lambda}|_{A[\fp]\times A[\fp]}=e_2^\lambda|_{A[\fp]\times A[\fp]}$. Thus the quadratic refinement of $e_2^\lambda$, restricted to $A[\fp]$, gives a quadratic refinement of $e_\fp^\lambda$ over $F$.  By Proposition \ref{self-cup}, $c_\lambda=0$. (2) follows by Theorem \ref{Rational-condition}.
\end{proof}

\section{Quadratic structures on cohomologies and Selmer groups}\label{QSS}
We distinguish the following cases:
\begin{itemize}
\item[$(\Ort)$] $\dag|_K=1$, i.e.  $K=K_0$,
\item[($\Sym$)] $\dag|_K\neq1 $ and $\fp$ is ramified in $K/K_0$, 
\item[$(\Uni)$] $\dag|_K\neq1$ and $\fp$ is inert in $K/K_0$, or
\item[$(\SU)$] $\dag|_K\neq1$ and $\fp_0$ is split in $K/K_0$.
\end{itemize}
We first focus on the $(\Ort)$, $(\Sym)$ and $(\Uni)$ cases where the prime ideal $\fp$ is stable under the Rosati involution $\dag$ and the split case $(\SU)$ is delayed into the last subsection.

\subsection{Local cohomologies}
\subsubsection{Quadratic structures on local cohomologies}
Let $v\in \Omega$ be a place of $F$. As explained in \cite[VI-4]{Shatz72} and \cite[III. 6.5-6.7]{Milne-ADT06}, for any group scheme $G$ of finite type over $F_v$, the cohomology group $\RH^r(F_v,G)$ admits a natural topology making it locally compact and $\sigma$-compact. Under this topology, maps arising from derived exact sequences and cup products are continuous. The group $\RH^2(F_v,\BF_p(1))$ is given the discrete topology. For an abelian variety $B/F_v$, $B(F_v)$ is compact, and thus $\RH^1(F_v, B)$ is discrete by duality (cf. \cite{Tate95}, \cite[I.3.4, I.3.7, III.7.8]{Milne-ADT06}). The group $\RH^1(F_v,A[\fp])$ is a locally compact $k$-space with the $k$-structure inherited from that of $A[\fp]$. Moreover, If $p\neq \ell$, then $\RH^1(F_v,A[\fp])$ is finite (cf. \cite[I.2.3, I.2.13(a)]{Milne-ADT06}).

The cup-product pairing with respect to the pairing $e_\fp^\lambda$ induces the local Tate pairing
\begin{equation*}\label{tate-pairing}
h_v'=\mathop{\bigcup}_{e_\fp^\lambda}:\RH^1(F_v,A[\fp])\times\RH^1(F_v,A[\fp])\lra \RH^2(F_v,\BF_p(1))\xhookrightarrow[]{\Inv_v} \BF_p
\end{equation*}
which is $\BF_p$-bilinear. The last inclusion is the Hasse invariant map  which is an  isomorphism for non-archimedean places. The same construction applied to $\theta_\fp^\lambda$ yields an extended local pairing 
\begin{equation}\label{extended-pairing}
h_v=\mathop{\bigcup}_{\theta_\fp^\lambda}:\RH^1(F_v,A[\fp])\times\RH^1(F_v,A[\fp])\lra \RH^2(F_v,k(1))\xhookrightarrow[]{k\otimes \Inv_v} k. 
\end{equation}
The Rosati involution $\dag$ on $\CO$ induces an involution on the residue field $k$ which we also denoted as $\dag$. The induced involution is trivial in the $(\Ort)$ or $(\Sym)$ case, and nontrivial in the $(\Uni)$ case. The pairing $h_v'$ is a $\dag$-sesquilinear form of $k$-spaces. 
\begin{prop}\label{general-linear}
Let $v\in \Omega$ be a place of $F$.  The pairing $h_v'$ is a $\dag$-adjoint $\BF_p$-bilinear form and $h_v$ is a hermitian $k$-form with respect to the involution $\dag$  satisfying $\Tr_{k/\BF_p}(\alpha_0 h_v)=h'_v$. Both $h_v'$ and $h_v$ are nondegenerate with symmetry as listed below:
\begin{center}
\begin{tabular}{c|c|c|c}
&$(\Ort)$&$(\Sym)$&$(\Uni)$\\
\hline
$h'_v$&symmetric $\BF_p$-bilinear&anti-symmetric $\BF_p$-bilinear&symmetric $\BF_p$-bilinear\\
\hline
$h_v$&symmetric $k$-bilinear&anti-symmetric $k$-bilinear&$1$-hermitian\\
\hline
\end{tabular}
\end{center}
\end{prop}
\begin{proof}
Since $p\nmid \deg(\lambda)$, both $e_\fp^\lambda$ and $\theta_\fp^\lambda$ are nondegenerate, so are $h_v'$ and $h_v$ by local Tate duality \cite{Tateduality, Milne-ADT06}. The $\dag$-adjointness of $h_v'$ follows from that of $e_\fp^\lambda$. The symmetry of $h_v'$ and $h_v$ follows from Proposition \ref{muwp}.(2), Corollary \ref{mod-Weil}  and  the dimension-shifting of cup products \cite[\S IV, Proposition 9]{CF-ANT}. The trace relation $\Tr_{k/\BF_p}(\alpha_0 h_v)=h'_v$ follows from Proposition \ref{compare-weil}.
\end{proof}

\begin{prop}\label{even-tate}
Let $v\in \Omega$ be a place of $F$. If $p=2$, suppose $(A,\lambda)$ is quadratic at $v$ in the $(\Ort)$ and $(\Sym)$ cases. Then the pairing $h_v'$ is an even $\BF_p$-form $($with the trivial involution on $\BF_p$$)$. 
\end{prop}
\begin{proof}
(1) {Case $p\neq 2$}. It follows from Lemma \ref{auto-quadratic} that $h_v'$ is even.
(2) Case $(\Ort)$ or $(\Sym)$ when $p=2$. By assumption $(A,\lambda)$ is quadratic at $v$. It follows from Proposition \ref{quadratic-alternating}, the pairing $h'_v$ is alternating, and thus even by Lemma \ref{even-alternating}.  (3) Case $(\Uni)$ when $p=2$. Let $\alpha\in k$ be an element that $\alpha+\alpha^\dag=1$ and set $s=\alpha \alpha_0 h_v$. Since $h_v$ is $1$-hermitian, $\alpha_0 h_v=s+s^*$. If we put $s'=\Tr_{k/\BF_2} s$ , then,  by the trace relation of Proposition \ref{general-linear},  $h'_v=s'+s'^*$ is even as desired.
\end{proof}

One may consult Proposition \ref{local-criterion} and \ref{global-criterion1} on cases where $(A,\lambda)$ is quadratic. In case of the above proposition, $(\RH^1(F_v,A[\fp]),h_v')$ is an even hermitian $\BF_p$-space. By Remark \ref{quadratic-evenhermitian2}, to give  a quadratic $\BF_p$-structure on $\RH^1(F_v,A[\fp])$, it suffices to specifying the form parameter and  the quadratic map and we do this as follows:

\noindent{\bf $\bullet$ $(\Ort)$ case.} The form parameter is taken to be $\Lambda=\Lambda^{\min}=0$. If $p\neq 2$, the quadratic structure is uniquely determined by $h_v'$ and $q_v'(x)=h_v'(x,x)/2$. If $p=2$, take $q_v'=q_\fp$ where $q_\fp$ is given by Proposition \ref{RM-quadratic-map} with $S=\Spec F_v$ and a fixed line bundle $L$ over $F_v$ such that $\varphi_L=\lambda \varpi$ (cf. Remark \ref{local-quadratic}.(1)).  This results in an orthogonal $\BF_p$-sapce $(\RH^1(F_v,A[\fp]),h_v',q_v')$.

\noindent {\bf $\bullet$ $(\Sym)$ case.} The form parameter is taken to be $\Lambda=\Lambda^{\max}=\BF_p$. Then the even hermitian space $(\RH^1(F_v,A[\fp]),h_v')$ uniquely determines a symplectic $\BF_p$-space $(\RH^1(F_v,A[\fp]),h_v',q_v')$ with $q_v'=0$.

\noindent {\bf $\bullet$ $(\Uni)$ case.} We take the form parameter $\Lambda=\Lambda^{\min}=0$. As in the proof of Proposition \ref{even-tate},  fix an element $\alpha\in k$ satisfying $\alpha+\alpha^\dag=1$. Then the even hermitian space $(\RH^1(F_v,A[\fp]),h_v')$ determines an orthogonal $\BF_p$-space $(\RH^1(F_v,A[\fp]),h_v',q_v')$ with $q_v'(x)=\Tr_{k/\BF_p}(\alpha \alpha_0 h_v(x,x))$.

To summarize, we have

\begin{prop}\label{adjoint-quadratic2}
Let $v\in \Omega$ be a place of $F$. If $p=2$, suppose $(A,\lambda)$ is quadratic at $v$ in the $(\Ort)$ and $(\Sym)$ cases. Then in the $(\Ort)$ or $(\Uni)$ case resp. $(\Sym)$ case, $(\RH^1(F_v,A[\fp]),h_v',q_v')$ is  a $\dag$-adjoint orthogonal resp. symplectic $\BF_p$-space.
\end{prop}

Invoking Theorem \ref{lifting-thm} and Proposition \ref{general-linear}, the $\dag$-adjoint quadratic $\BF_p$-space $(\RH^1(F_v,A[\fp]),h_v',q_v')$ determines a unique trace-compatible quadratic $k$-space $(\RH^1(F_v,A[\fp]),h_v,q_v)$.

\begin{thm}\label{quadratic-structure}
Assume the same hypothesis as in Proposition \ref{adjoint-quadratic2} and let $\alpha_0\in k^\times$ be as in Proposition \ref{compare-weil}. 
\begin{itemize}
\item[(1)] In the $(\Ort)$ case, $(\RH^1(F_v,A[\fp]), h_v,q_v)$ is an orthogonal $k$-space.
\item[(2)] In the $(\Sym)$ case,  $(\RH^1(F_v,A[\fp]), h_v,q_v)$ is a symplectic  $k$-space.
\item[(3)] In the $(\Uni)$ case, $(\RH^1(F_v,A[\fp]), h_v,q_v)$ is a unitary $k$-space.
\end{itemize}
In all cases, we have $\Tr_{k/\BF_p}(\alpha_0 h_v)=h'_v$ and $\Tr_{k/\BF_p}(\alpha_0 q_v)=q_v'$.
\end{thm}

\subsubsection{Local Kummer descents and metabolic structures}
The exact sequence 
\[0\ra A[\fp]\lra A\xrightarrow{h(i_\fp,1_A)} \shom{\fp}{A}\ra 0\]
gives rise to the local Kummer $\fp$-descent exact  sequence
\begin{equation}\label{kummer-sequence}
\Hom_\CO(\fp,A(F_v))\xrightarrow{\delta_{A, \fp,v}}\RH^1(F_v,A[\fp])\lra \RH^1(F_v, A)[\fp]\ra0.
\end{equation}
Denote the local Kummer image as
\[\CL_{A,\fp,v}:=\Im\left( \Hom_\CO(\fp,A(F_v))\xrightarrow{\delta_{A, \fp,v}}\RH^1(F_v,A[\fp])\right).\]
Since $\Hom_\CO(\fp,A(F_v))$ is compact and $\RH^1(F_v,A)$ is discrete, $\CL_{A,\fp,v}$ is compact open in $\RH^1(F_v,A[\fp])$.

Let $v$ be a finite place of $F$  and let $\CO_{F,v}$ be the valuation ring of $F_v$. Suppose $A$ has good reduction at $v$, i.e.  $A$ extends to an abelian scheme $\CA$ over $\CO_{F,v}$. Then the fppf cohomology group $\RH^1(\CO_{F,v},\CA[\fp])$ is an open subgroup of $\RH^1(F_v,A[\fp])$. Let $\CB$ be a finite set of places of $F$ containing all the archimedean places and the places where $A$ has bad reductions.
\begin{prop}\label{isotropic}
\begin{itemize}
\item[(1)] For any $v\in \Omega$, $\CL_{A,\fp,v}$ is a compact open subspace of $\RH^1(F_v,A[\fp])$ which is the orthogonal complement of itself with respect to both $h_v'$ and $h_v$.
\item[(2)] For $v\notin \CB$, $\CL_{A,\fp,v}=\RH^1(\CO_{F,v},\CA[\fp])$.
\end{itemize}
\end{prop}
\begin{proof}
(1) Since $h_v$ and $h_v'$ are trace-compatible, it suffices to prove  $\CL_{A,\fp,v}$ is self-orthogonal under $h_v'$ which is a consequence of local Tate dualities \cite{Tateduality,Tate95, Milne-ADT06}. For simplicity, we denote $g=h(i_\fp,1_A), B=\shom{\fp}{A}$ and thus $g^\vee=h(j_\fp,1_{A^\vee}), B^\vee=\shom{\fp^{-1}}{A^\vee}$. Let $A(F_v)_\bullet$ be $A(F_v)$ modulo its identity connected component which is nonzero only if $F_v=\BR$ or $\BC$. The exact sequence 
\[0\ra A[g]\ra A\xrightarrow{g}B\ra 0\]
gives rise to the following commutative diagram (possibly up to sign):
\begin{equation}\label{kummer-duality}
\xymatrix@C=9pt{0\ar[r]&B(F_v)_\bullet/g(A(F_v)_\bullet)\ar[d]_{\alpha_1}\ar[r]&\RH^1(F_v,A[g])\ar[r]\ar[d]_{\alpha_2}&\RH^1(F_v,A)[g]\ar[r]\ar[d]_{\alpha_3}&0\\
0\ar[r]&\left(\RH^1(F_v,B^\vee)[g^\vee]\right)'\ar[r]&\left(\RH^1(F_v,B^\vee[g^\vee])\right)'\ar[r]&\left(A^\vee(F_v)_\bullet/g^\vee(B^\vee(F_v)_\bullet)\right)'\ar[r]&0.}
\end{equation}
Here $(\cdot)'=\Hom_{\BF_p}(\cdot,\BF_p)$ denotes the Pontryagin dual. The maps $\alpha_1,\alpha_3$ are isomorphisms induced from Tate duality for abelian varieties over local fields  (\cite{Tate95}, \cite[I.3.4, I.3.7, III.7.8]{Milne-ADT06}). The map $\alpha_2$ is induced from the cup-product pairing with respect to the Weil pairing $e_\fp=e_g$ in (\ref{Weil-p}) and is an isomorphism by Tate dualities for finite group schemes (\cite{Tateduality} , \cite[I.2.3, I.2.13, III.6.10]{Milne-ADT06}). The bottom row is the Pontryagin dual of the Kummer descent exact sequence induced by the dual isogeny $g^\vee$. A detailed account for the commutativity of the diagram can be found in \cite[p21]{Skoronotes}.

Composing the isogeny $\lambda$, by Proposition \ref{dual-kummer}, the diagram (\ref{kummer-duality}) transforms into
\begin{equation*}
\xymatrix@C=9pt{0\ar[r]&B(F_v)_\bullet/g(A(F_v)_\bullet)\ar[d]_{}\ar[r]&\RH^1(F_v,A[g])\ar[r]\ar[d]_{\alpha_2^\lambda}&\RH^1(F_v,A)[g]\ar[r]\ar[d]_{}&0\\
0\ar[r]&\left(\RH^1(F_v,A)[g]\right)'\ar[r]&\left(\RH^1(F_v,A[g])\right)'\ar[r]&\left(B(F_v)_\bullet/g(A(F_v)_\bullet)\right)'\ar[r]&0,}
\end{equation*}
with vertical isomorphisms. The map $\alpha_2^\lambda$ coincides with the one induced by  $\cup_{e_\fp^\lambda}$ and therefore Assertion (1) follows.

\noindent (2) The second assertion is \cite[Lemma 1.2]{Ulmer91}.
\end{proof}

\begin{prop}\label{metabolic}
Let $v\in \Omega$ be a place of $F$. Under the same hypothesis as in Proposition \ref{adjoint-quadratic2}, $\CL_{A,\fp,v}$ is a compact open maximal isotropic subspace of the quadratic $k$-space $(\RH^1(F_v,A[\fp]),h_v,q_v)$ which is thus metabolic.
\end{prop}
\begin{proof}
By Theorem \ref{lifting-thm} and \ref{quadratic-structure} , it suffices to prove $\CL_{A,\fp,v}$ is maximal isotropic in $(\RH^1(A[\fp]),h_v',q_v')$.  If $p\neq 2$ or $(\RH^1(A[\fp]),h_v',q_v')$ is symplectic,  the form parameter is $\Lambda^{\max}$. Any subspace $W$ with $W^\perp=W$ is automatically maximal isotropic by Lemma \ref{orth-iso}. As for the $(\Ort)$ case with $p=2$, recall $q_v'=q_\fp$ of Proposition \ref{RM-quadratic-map} with $S=\Spec F_v$ and a fixed line bundle $L$ over $F_v$ such that $\varphi_L=\lambda \varpi$. Therefore $q'_v(\CL_{A,\fp,v})=0$ is just Proposition \ref{RM-quadratic-map}.(2). 

\end{proof}

\subsection{Ad\'elic cohomologies and Selmer groups}\label{global-cohomology}
For each $v\in \Omega$, $\RH^1(F_v,A[\fp])$ is $\sigma$-compact and hence second-countable by \cite[Corollary 2.12]{PR12}. Equipped with the restricted product topology, the restricted product 
\[\RH^1(\BA_F,A[\fp]):=\sideset{}{'}\prod_{v\in \Omega} \left(\RH^1(F_v,A[\fp]), \CL_{A,\fp,v}\right)\]
is a second countable locally compact $k$-space (cf. \cite[\S 13]{Cassels-ANT} and \cite[Example 2.18]{PR12}).

To give quadratic $k$-structures on $\RH^1(\BA_F,A[\fp])$, we need the local quadratic structures everywhere. To this end, by Theorem \ref{quadratic-structure}, if $p=2$ we need additional hypotheses: 
\begin{itemize}
\item In the ($\Sym$) case with $p=2$, assume $(A,\lambda)$ is quadratic everywhere;
\item  In the $(\Ort)$ case with $p=2$, assume the global hypothesis ({\bf L}) of Section \ref{global-hyp} holds. By Proposition \ref{additional}, under the hypothesis ({\bf L}), $(A,\lambda)$ is quadratic everywhere. Moreover, It follows from Proposition \ref{rational-line-bundle}, there exists a rational symmetric line bundle $L$ on $A$ such that $\varphi_L=\lambda \varpi$. We use the coherent system of local orthogonal structures on $\RH^1(F_v,A[\fp])$ that are constructed from this single $L$ via Proposition \ref{RM-quadratic-map} with $S=\Spec F_v$  for all $v\in \Omega$. 
\end{itemize}
By Proposition \ref{metabolic}, for any $x=(x_v)$ and  $y=(y_v)\in \RH^1(\BA_F,A[\fp])$, $h_v(x_v,y_v)=q_v(x_v)=0$ for all but finitely many $v\in \Omega$. Thus $h=\sum_{v\in \Omega} h_v$ and $q=\sum_{v\in \Omega} q_v$ are well-defined. It follows from Theorem \ref{quadratic-structure} that $h$ is an even hermitian form and $q$ is a quadratic map on $\RH^1(\BA_F,A[\fp])$, so that $(\RH^1(\BA_F,A[\fp]),h,q)$ is a nondegenerate quadratic $k$-space. The subspace $\CL=\prod_{v\in \Omega} \CL_{A,\fp,v}$ is compact open and maximal isotropic by Proposition \ref{metabolic}, and thus $(\RH^1(\BA_F,A[\fp]),h,q)$ is a metabolic quadratic $k$-space.

By Proposition \ref{isotropic}.(2), 
\[\RH^1(\BA_F,A[\fp])=\sideset{}{'}\prod_{v\in \Omega} \left(\RH^1(F_v,A[\fp]), \RH^1(\CO_{F,v},\CA[\fp])\right).\]
Since each element of $\RH^1(F,A[\fp])$ belongs to $\RH^1(\CO_{F,v},\CA[\fp])$ for all but finitely many $v\in \Omega\backslash \CB$, there is a well-defined map 
\begin{equation}\label{global}\RH^1(F,A[\fp])\ra \RH^1(\BA_F,A[\fp]).\end{equation} 
We denote its image by $W$ and kernel by $\Sha^1(F,A[\fp])$. Then $$W\cong \RH^1(F,A[\fp])/\Sha^1(F,A[\fp]).$$ 
Define the $\fp$-Selmer group
\[\Sel_\fp(A):=\Ker\left(\RH^1(F,A[\fp])\ra \prod_{v\in \Omega}\RH^1(F_v,A)\right).\]

\begin{thm}\label{intersection}
If $p=2$, suppose the hypothesis $({\bf L})$ holds in the $(\Ort)$ case and $(A,\lambda)$ is quadratic everywhere in the $(\Sym)$ case. 
\begin{itemize}
\item[(1)] $(\RH^1(\BA_F,A[\fp]),h,q)$ is a metabolic orthogonal, symplectic resp. unitary $k$-space in the $(\Ort)$, $(\Sym)$ resp. $(\Uni)$ case accordingly. 
\item[(2)] Both $\CL$ and $W$ are maximal isotropic subspaces in $\RH^1(\BA_F,A[\fp])$ and the map $(\ref{global})$ induces an isomorphism
\[\Sel_\fp(A)/\Sha^1(F,A[\fp])\cong \CL\cap W.\]
\end{itemize}

\end{thm}
\begin{proof}
The isomorphism follows from the definition of Selmer groups and the Kummer exact sequence (\ref{kummer-sequence}). To prove $W$ is maximal isotropic, by Theorem \ref{lifting-thm} and Theorem \ref{quadratic-structure},  it suffices to prove $W$ is so for the $\dag$-adjoint quadratic $\BF_p$-space $(\RH^1(\BA_F,A[\fp]),h',q')$ where $h'=\sum_{v\in \Omega}h_v'$ and $q'=\sum_{v\in \Omega}q_v'$. As in the proof of \cite[Theorem 4.14]{PR12}, using the pairing $e_\fp^\lambda$ instead of $e_\lambda$ there, by the global Poitou-Tate duality, $W$ is the orthogonal complement of itself with respect to the cup-product pairing $h'$. By Lemma \ref{orth-iso}, $W$ is maximal isotropic if $p\neq 2$ or $(\RH^1(\BA_F,A[\fp]),h',q')$ is symplectic. 

Suppose $p=2$ and $(\RH^1(\BA_F,A[\fp]),h',q')$ is orthogonal. The construction of the quadratic maps $q_v'$ in Proposition \ref{RM-quadratic-map} via the rational line bundle $L$ is functorial with respect to base extensions. So for any $w\in W$ and any $v\in \Omega$, $q_v'(w)$ can be evaluated through a global quadratic map $\RH^1(F,A[\fp])\ra \RH^2(F,\BG_m)$ at $x\in \RH^1(F,A[\fp])$ where $w$ is the image of $x$. By the reciprocity law for Brauer groups, we have an exact sequence
\[0\ra \RH^2(F,\BG_m)\lra \bigoplus_{v\in \Omega}\RH^2(F_v,\BG_m)\xrightarrow{\sum_{v\in \Omega} \Inv_v} \BQ/\BZ\ra0.\]
Thus $q'(W)=0$. 
\end{proof}
\begin{remark}\label{vanishing-sha}
The group $\Sha^1(F,A[\fp])$ is often zero by the vanishing criterion \cite[Proposition 3.3]{PR12} and see for examples of hyperelliptic curves there \cite[Proposition 3.4]{PR12}.  Suppose either $[K:\BQ]=\dim A$ or $2\dim A$ and $p\neq \ell$. Then $A[\fp]\cong k^i$  and $\Aut(A[\fp])\cong \GL_i(k)$ for $i=1$ or $2$. In either case the $p$-sylow subgroups of $\Aut(A[\fp])$ are cyclic and by \cite[Proposition 3.3.(d)]{PR12}, $\Sha^1(F,A[\fp])=0$.
\end{remark}

\subsection{The split case}\label{split-case}
Finally we discuss the case $(\SU)$: $\dag|_K\neq 1$ and $\fp_0$ splits in $K/K_0$. In this case $\fp_0\CO=\fp\fp^\dag$ and $\CO_{\fp_0}=\CO_{\fp}\oplus \CO_{\fp^\dag}$. The residue field $k_0=\CO_0/\fp_0$ embeds into $k=\CO/\fp_0\CO=k_0\oplus k_0$ diagonally and the induced involution $\dag$ on $k$ acts as $(x,y)^\dag=(y,x)$. By Proposition \ref{extended-Weil}, reducing the extended pairing $\Theta_{\fp_0}^\lambda$ modulo $\fp_0$, we get a nondegenerate $(-1)$-hermitian pairing 
\[\theta_{\fp_0}^\lambda:A[\fp_0]\times A[\fp_0]\lra k(1).\]
For each $v\in \Omega$, the local Tate-duality yields a nondegenerate $1$-hermitian pairing of $k$-spaces
\begin{equation}\label{extended-pairing2}
h_v:\RH^1(F_v,A[\fp_0])\times \RH^1(F_v,A[\fp_0])\lra k.
\end{equation}
We have $A[\fp_0]=A[\fp]\oplus A[\fp^\dag]$, and 
$$\RH^1(F_v,A[\fp_0])=\RH^1(F_v,A[\fp])\bigoplus \RH^1(F_v,A[\fp^\dag]).$$
 The nondegenerate hermitian form $h_v$ identifies $\RH^1(F_v,A[\fp^\dag])$ as the $k_0$-dual of $\RH^1(F_v,A[\fp])$ and furnishes $\RH^1(F_v, A[\fp_0])$ with the split unitary $k$-structure. Moreover, $h=\sum_{v\in \Omega}h_v$ determines $\RH^1(\BA_F, A[\fp_0])$ as a split unitary $k$-space. 
 
Instead, if we take $K'=K'_0=K_0$ and view $(A,\lambda)$ as in the $(\Ort)$ case. Via the diagonal embedding $k_0\hookrightarrow k$,  $\RH^1(\BA_F,A[\fp_0])$ is a $k_0$-space and acquires a metabolic orthogonal $k_0$-space structure which is trace compatible with the split unitary $k$-structure. Note we don't need the additional hypothesis in Theorem \ref{quadratic-structure} to ensure this orthogonal $k_0$-structure to exist. Indeed, this orthogonal structure and the one given in Theorem \ref{quadratic-structure} have the same hermitian form $\Tr_{k/k_0}h_v$ (in notation here).  But the quadratic forms may be different if $p=2$ (cf. Section \ref{split-unitary}). Since the Kummer subspace $\CL=\prod_{v\in \Omega}\CL_{A,\fp_0,v}$ and the image of the global chomology $\RH^1(F,A[\fp_0])$ are all $k$-subspaces and self-orthogonal under $\Tr_{k/k_0}h_v$, they are  also self-orthogonal under $h_v$, and therefore maximal isotropic for the split unitary structure (cf. Section \ref{split-unitary}). In particular, as in Theorem \ref{intersection}, the Selmer group $\Sel_{\fp_0}(A)/\Sha^1(F,A[\fp_0])$ is also an intersection of two maximal isotropic subspaces in the split unitary $k$-space $\RH^1(\BA_F,A[\fp_0])$.

\section{Shafarevich-Tate groups}
\subsection{The primary components of Shafarevich-Tate groups}
The Shafarevich-Tate group $\Sha(A)$ of  the abelian variety $A$ is defined  as 
\[\Sha(A)=\Ker\left(\RH^1(F,A)\lra \prod_{v\in \Omega} \RH^1(F_v,A)\right),\]
and denote $\Sha(A)_{/\div}$ to be the quotient group of $\Sha(A)$ modulo its divisible part. The Cassels-Tate pairing
\[\langle\ ,\ \rangle_{\mathrm{CT}}: \Sha(A)_{/\div}\times \Sha(A^\vee)_{/\div}\lra \BQ/\BZ\]
is a perfect pairing satisfying for any $\phi\in \End_F(A)$
\[\CTP{\phi(x)}{y}=\CTP{x}{\phi^\vee(y)}.\]
A detailed account for the Cassels-Tate pairings can be found in \cite{PS99} and \cite[I.6]{Milne-ADT06}. 

Let the notations be as in Section \ref{QC}. Composing the symmetric isogeny  $\lambda$, we obtain the pairing $\LP{x}{y}=\CTP{x}{\lambda(y)}$. Since $\lambda$ is $\dag$-sesquilinear under $\CO$ and symmetric, $\LP{\ }{\ }$ is $\dag$-adjoint. We restrict to the $\fp_0$-primary component: 
\[{\poLP{\ }{\ }}:\Sha(A)_{/\div}[\fp_0^\infty]\times \Sha(A)_{/\div}[\fp_0^\infty]\lra \BQ_p/\BZ_p.\]
By Corollary \ref{sesquilinear-adjoint}, there exists a unique $\dag$-sesquilinear pairing 
\[\Psi_{\lambda,\fp_0}:\Sha(A)_{/\div}[\fp^\infty_0]\times \Sha(A)_{/\div}[\fp_0^\infty]\lra K_{\fp_0}/{\CO_{\fp_0}}=\CO_{\fp_0}\otimes\BQ_p/\BZ_p\]
which is characterized by, for any $x,y\in \Sha(A)_{/\div}[\fp_0^\infty]$,
\[{\poLP{x}{y}}=\Tr_{K_{\fp_0}/\BQ_p}(d\Psi_{\lambda,\fp_0}(x,y)).\]

\begin{prop}\label{sha}
\begin{itemize}
\item[(1)] The pairing $\Psi_{\lambda,\fp_0}$ is a nondegenerate, $\dag$-sesquilinear and $(-d^{\dag-1})$-hermitian pairing of finite $\CO_{\fp_0}$-modules. 
\item[(2)] Suppose $\dag|_K=1$ and if $p=2$, we further assume that $2$ is unramified in $\CO$ and $\lambda$ is induced from some rational line bundle. Then $\Psi_{\lambda,\fp}$ is $\CO_{\fp}$-bilinear and alternating. In particular, $\Sha(A)_{/\div}[\fp^\infty]\cong M\bigoplus M$ for some finite $\CO_{\fp}$-module $M$. 
\end{itemize}
\end{prop}
\begin{proof}
(1) The finiteness of $\Sha(A)_{/\div}[\fp_0^\infty]$ can be found in  \cite[Proposition 5.14]{BK90}  or \cite[Theorem 1]{Flach90}. Since $p\nmid \deg(\lambda)$ and $\lambda$ is symmetric, the induced Weil pairing $E_{p}^\lambda$ 
is nondegenerate and alternating. By \cite[Theorem 1-2]{Flach90}, the pairing $\poLP{\ }{\ }$ is nondegenerate and anti-symmetric.  Then (1) follows from  Corollary \ref{sesquilinear-adjoint}. 

(2) If $\dag|_K=1$, then $\fp=\fp_0$ and by (1) $\Psi_{\lambda,\fp}$ is $\CO_\fp$-bilinear and anti-symmetric. If $p\neq 2$, then it is alternating. For $p=2$, by Proposition \ref{rational-line-bundle}, for any $a\in \CO$ the symmetric isogeny $\lambda a$  is induced from some rational line bundle. Then $\langle{\ },{\ }\rangle_{ \lambda a,\fp}$ is alternating (cf. \cite[Theorem 3.3]{Tateduality} and \cite[Corollary 7]{PS99}). In particular, for all $a\in \CO$ and $x\in \Sha(A)_{/\div}[\fp^\infty]$, 
\[\Tr_{K_\fp/\BQ_p}(da\Psi_{\lambda,\fp}(x,x))=\langle x,ax\rangle_{\lambda,\fp}=\langle x,x\rangle_{\lambda a,\fp}=0.\] 
By a similar argument of Proposition \ref{orthogonal}, we find $\Psi_{\lambda,\fp}$ is alternating. The last assertion follows by a well-known linear algebra argument. 
\end{proof}
\begin{remark}
When $p=2$, the assumption $2$ is unramified in $\CO$ seems sort of exact for $\Psi_{\lambda,\fp}$ to be alternating. One might hope to prove the alternatingness of  $\Psi_{\lambda,\fp}$ from that of $\langle\ ,\ \rangle_{\lambda,\fp}$, i.e. for all $a\in \CO_\fp$
\[\Tr_{K_\fp/\BQ_p}(da^2\Psi_{\lambda,\fp}(x,x))=\langle ax,ax\rangle_{\lambda,\fp}=0.\] 
This doesn't suffice. For example  if one takes $\CO=\BZ[\sqrt{2}]$ and $\fp=(\sqrt{2})$, then $\Tr_{K_\fp/\BQ_2}(d b \CO^2)\subset 2^m\BZ_2$ for $m\geq 0$ doesn't imply $b\in 2^m\CO_\fp$ (cf. the proof of Proposition \ref{orthogonal}).
\end{remark}

\begin{thm}\label{square-sha0}
\begin{itemize}
\item[(1)] If $p\neq 2$, then $\Sha(A)_{/\div}[p^\infty]\cong M\bigoplus M$ for some finite $\CO_{0}$-module $M$.
\item[(2)] Assume $p=2$, $\ell\neq 2$ and $2$ is prime to the discriminant of $\CO$. Suppose either $\dag|_K=1$ and $[K:\BQ]=\dim A$ or $\dag|_K\neq 1$. Then $\Sha(A)_{/\div}[2^\infty]\cong M\bigoplus M$ for some finite $\CO_{0}$-module $M$.
\end{itemize}
\end{thm}
\begin{proof}
If $p=2$, under the hypotheses, it follows from Theorem \ref{Rational-condition} that $\lambda$ is induced from a rational symmetric line bundle. Replacing $(K,\dag)$ by $(K_0,\dag|_{K_0})$, apply Proposition \ref{sha} to various primes $\fp_0$ of $\CO_0$  above $p$ and then the theorem follows from Proposition \ref{sha} and the decomposition
\[\Sha(A)_{/\div}[p^\infty]=\bigoplus_{\fp_0\mid p} \Sha(A)_{/\div}[\fp_0^\infty].\]
\end{proof}

\subsection{The $\fp$-parts of Shafarevich-Tate groups}
\begin{lem}\label{uniformizer}
Suppose $\dag|_K\neq 1$ and $\fp^\dag=\fp$. We can choose a uniformizer $\varpi$ of $\CO_\fp$ so that $\varpi^{\dag-1}=1$ resp. $\equiv -1\mod \fp$ if $\fp$ is inert resp. ramified in $K/K_0$.
\end{lem}

\begin{proof}
If $\fp$ is inert, just choose a uniformizer of $\CO_{0,\fp_0}$. Assume $K_\fp=K_{0,\fp_0}(\sqrt{A})$ is ramified where $A\in K_{0,\fp_0}$. If $A=\varpi_0$ a uniformizer of $K_{0,\fp_0}$, then take $\varpi=\sqrt{\varpi_0}$. If $p=2$ and the square-free normal form $\alpha$ of $A$ satisfies $\alpha-1=u \varpi_0^{2i+1}$  with $i<e$ where $u$ is a unit and $(2)=(\varpi_0)^e$, then take 
\[\varpi=\frac{1-\sqrt{\alpha}}{\varpi_0^i}\]
which is a root of the Eisenstein polynomial 
\[x^2+\left(\frac{2}{\varpi_0^i}\right)x+\left(\frac{1-\alpha}{\varpi_0^{2i}}\right)\]
and satisfies $\varpi^{\dag-1}\equiv -1\mod \fp$ (cf. \cite[\S2]{Casselman-quadratic}).
\end{proof}

The pairing $\Psi_{\lambda,\fp_0}$, via restricting,  induces a nondegenerate pairing 
\[\psi_{\lambda,\fp}:\Sha(A)_{/\div}[\fp]\times \Sha(A)_{/\div}[\fp^\dag]\xrightarrow{\Psi_{\lambda,\fp_0}} \fp^{-1}\CO_\fp/\CO_\fp\cong \CO_\fp/\fp\CO_\fp.\]
The last identification depends on a uniformizer $\varpi$ of $\CO_\fp$. If $\dag|_K\neq 1$ and $\fp^\dag=\fp$, we choose the uniformizer as in Lemma \ref{uniformizer}.

\begin{coro}\label{p-sha}
\begin{itemize}
\item[(1)]  In the $(\SU)$ case, $\left(\Sha(A)_{/\div}[\fp_0],\psi_{\lambda,\fp}\oplus \psi_{\lambda,\fp}^*\right)$ is a split unitary  $k$-space and 
\[\dim_{k_0} \Sha(A)_{/\div}[\fp]=\dim_{k_0} \Sha(A)_{/\div}(\fp^\dag).\]
\end{itemize}
The following statements $(2)$ and $(3)$ always hold for $p\neq 2$. 
\begin{itemize}
\item [(2)] In the $(\Ort)$ and $(\Sym)$ cases, $(\Sha(A)_{/\div}[\fp],\psi_{\lambda,\fp})$ is a symplectic $k$-space.
\item[(3)]  In the $(\Uni)$ case, $(\Sha(A)_{/\div}[\fp],\psi_{\lambda,\fp})$ is a unitary $k$-space.
\end{itemize}
Suppose $p=2$ and $\ell\neq 2$. If  
\begin{itemize}
\item in the $(\Ort)$ case, $[K:\BQ]=\dim A$ and the discriminant $d_\CO$ is odd,
\item in the $(\Sym)$ case, $[K:\BQ]=2\dim A$ and $d_{\CO_0}$ is odd,
\item in the $(\Uni)$ case, $d_\CO$ is odd, 
\end{itemize}
then the statements $(2)$ and $(3)$ still hold.
\end{coro}
\begin{proof}
(1) is clear and we suppose $\fp^\dag=\fp$. By Lemma \ref{uniformizer}, we always have $(d\varpi)^{\dag-1}\equiv 1\mod \fp$. Then it follows from Proposition \ref{sha} that $\psi_{\lambda,\fp}$ is nondegenerate, $\dag$-sesquilinear and $(-1)$-hermitian. For $x,y\in \Sha(A)_{/\div}[\fp]$, suppose $\psi_{\lambda,\fp}(x,y)=b\mod \fp.$ Through the identifications
\[\CO_\fp/\fp\CO_\fp\cong \fp^{-1}\CO_\fp/\CO_\fp\hookrightarrow \CO_\fp\otimes p^{-1}\BZ_p/\BZ_p,\]
one has $\Psi_{\lambda,\fp_0}(x,y)=\frac{pb}{\varpi}\otimes \frac{1}{p}$ and thus
\begin{eqnarray*}
\poLP{x}{y}&=&\Tr_{K_\fp/\BQ_p}(d\Psi_{\lambda,\fp_0}(x,y))=\Tr_{K_\fp/\BQ_p}\left(\frac{dp}{\varpi}b\right)\otimes \frac{1}{p}\\
&=&\Tr_{k/\BF_p}(\alpha \ov{b})=\Tr_{k/\BF_p}(\alpha \psi_{\lambda,\fp}(x,y))
\end{eqnarray*}
for some $\alpha \in k^\times$ and $\ov{b}$ denotes the reduction of $b\mod \fp$.  The first identity of the second row follows by Lemma \ref{trace-reduction} and the identification $p^{-1}\BZ_p/\BZ_p\cong \BF_p$. 

If $p\neq 2$, by \cite[Theorem 2]{Flach90}, $\poLP{\ }{\ }$ is alternating. If $p=2$, under the hypotheses, by Theorem \ref{Rational-condition}, $\lambda$ is induced from a rational line bundle. Here restricting the multiplication, we view the $(\Sym)$ case as the $(\Ort)$ case and apply Theorem \ref{Rational-condition}. Then by \cite[Theorem 3.3]{Tateduality} and \cite[Corollary 7]{PS99}, $\poLP{\ }{\ }$ is alternating. Invoking Theorem \ref{lifting-thm}, the $(-1)$-hermitian form $\psi_{\lambda,\fp}$ is even and determines the corresponding quadratic structure in each case. 
\end{proof}
\begin{coro}\label{p-sha1}
Suppose $\Sha(A)[\fp^\infty]$ is finite. In the $(\Ort)$ or $(\Sym)$ case of Corollary \ref{p-sha}, we have the parity
\[\dim_k A(F)\otimes_\CO k \equiv \dim_k \Sel_\fp(A)\mod 2.\]
\end{coro}
\begin{proof}
In the $(\Ort)$ or $(\Sym)$ case, $\Sha(A)_{/\div}[\fp]$ is a symplectic $k$-space and hence has even dimension. If $\Sha(A)[\fp^\infty]$ is finite, $\Sha(A)_{/\div}[\fp]=\Sha(A)[\fp]$. The congruence relation follows from the Kummer exact sequence
\[0\ra A(F)\otimes_\CO k\ra \Sel_\fp(A)\ra \Sha(A)[\fp]\ra0. \]
\end{proof}

\subsection{Examples}
\subsubsection{Abelian varieties with RM/CM}
Let $A$ be an abelian variety over a global field $F$ of characteristic $\ell\neq 2$. Suppose $A/F$ has RM by the ring of integers $\CO_K$ of a totally real number field $K$ and $2\nmid d_K$. There always exists an $\CO_K$-linear polarization $\lambda$ (cf. \cite[Proposition 1.10]{Rap78}). If  $\lambda$ is principal and $\Sha(A)$ is finite, then $\Sha(A)\cong M\bigoplus M$ for some finite $\CO_K$-module.  An abelian variety with CM by a CM field $K$ can be viewed as with RM by restricting the multiplication to $K_0$. 

Applying to elliptic curves $E/F$ together with their canonical theta polarizations, this recovers the well-known result of Cassels \cite{CA62} that if $\Sha(E)$ is finite, then $\Sha(E)\cong M\bigoplus M$ for some finite abelian group $M$. In particular, for any primes $p$, $\Sha(E)_{/\div}[p]$ aways have even $\BF_p$-dimensions. If $E$ has CM by the ring $\CO_K$ of integers of an imaginary quadratic field $K$, we can tell a little more that $\Sha(E)_{/\div}[\fp]$ also has even $\BF_p$-dimension for ramified primes $\fp$ of $K$ (cf. Corollary \ref{p-sha}-\ref{p-sha1}).

\subsubsection{Branched cyclic covering of $\BP^1$}
Let $q\geq 3$ be a prime. Let $F$ be a global field of characteristic $\ell\nmid 2q$ containing a primitive $q$-th root of unity $\omega$. Let $C\ra \BP^1$ be a branched cyclic covering of $\BP^1$ over $F$ of degree $q$ where $C$ is a smooth projective curve of genus $g\geq 1$.  By Kummer theory, $C$ admits an affine model $y^q=f(x)$ where $f(x)\in F[x]$ is a $q^{\mathrm{th}}$-power free polynomial.

Let $J=\Pic_{C/F}^0$ be the Jacobian variety. Let $K$ be the $q$-th cyclotomic field. In abuse of notation, let $\omega$ denote the automorphism $\omega (x,y)=(x,\omega y)$ on $C$. For any affine point $P\in C$, the divisor $\sum_{i=0}^{p-1} \omega^i P$ is trivial in $J$. Thus there is an embedding of the cyclotomic ring $\CO_K=\BZ[\omega]\hookrightarrow \End_F(J)$ via $\omega (x,y)=(x,\omega y)$.

\begin{prop}\label{theta}
Let $\lambda$ be the theta polarization of $J$ and  $\dag$ the associated Rosati involution. The Rosati involution $\dag$ stabilizes $\CO_K$ and restricts to the complex conjugation on $K$.
\end{prop}
\begin{proof}
The result doesn't depend on the base. So extending base, we may assume $C$ has an $F$-rational ramification point $Q$ and denote $D=[(g-1)Q]\in  \Pic^{g-1}_{C/F}(F)$. Recall the theta divisor $\Theta$ is the image of the natural map
\[C^{(g-1)}\lra \Pic^{g-1}_{C/F},\quad (P_1,\cdots,P_{g-1})\mapsto \left[\sum_{i=0}^{g-1} P_i\right].\]
We have the translation $t_D: \Pic_{C/F}^0\xrightarrow{\sim}\Pic^{g-1}_{C/F}$ and denote $\Theta_D=t^*_D \Theta$. Then the polarization $\lambda$ can be explicitly given as
\[\lambda(x)= (\Theta_D+x)-\Theta_D\in J^\vee.\]
This map is independent of the choice of $Q$ (cf. \cite[p261]{BLR90}). In order to check the stability of $\CO_K$, it suffices to show that $\lambda \omega^{-1}=\omega^\vee \lambda$. Since $\omega$ stabilizes $\Theta$ and $Q$, $\omega$ fixes the divisor $\Theta_D$. Then for any $x\in J$,
\[\omega^\vee\lambda(x)=\omega^*\lambda(x)=\omega^{-1}(\Theta_D+x)-\Theta_D)=(\Theta_D+\omega^{-1}x)-\Theta_D)=\lambda \omega^{-1}(x).\]
\end{proof}

A theta characteristic of the curve $C$ is a class of line bundles whose square is the canonical class $\omega_{C}\in \Pic_{C/F}^{2g-2}(F)$.
\begin{prop}
\begin{itemize}
\item[(1)] The curve $C$ has a rational theta characteristic. 
\item[(2)] If $\Sha(J)$ is finite, then $\Sha(J)\cong M\bigoplus M$ for some finite $\CO_{K_0}$-module. 
\item[(3)] If $\fq$ is a prime ideal of $\CO_K$ above $q$, then $\dim_{\BF_q} \Sha(J)_{/\div}[\fq]$ is even. 
\end{itemize}
\end{prop}
\begin{proof}
Since $q$ is odd, the discriminant $d_K$ is odd. By Theorem \ref{Rational-condition} and Proposition \ref{theta}, $\lambda=\varphi_L$ for some rational symmetric line bundle, or equivalently, $e_2^\lambda$ is alternating. Then (1) follows by \cite[Theorem 3.16]{PR11}. (2) follows by Theorem \ref{square-sha0} and Corollary \ref{square-sha} and (3) follows by Corollary \ref{p-sha} in the ($\Sym$) case.
\end{proof}
\begin{remark}
In \cite[Proposition 3.19]{PR11}, Poonen and Rains give some criterion for a hyper-ellipitc curve to have rational theta characteristics. 
\end{remark}

\appendix
\section{The $\underline{\Hom}$ sheaves and the Kummer descents}\label{Hom-construction}
Let $F$ be a field and let $A/F$ be an abelian variety with $A^\vee$ its dual abelian variety.  Let $\CO\subset \End_F(A)$ be a commutative integral domain containing the identity endowed with an involution $\dag$. Embed $\CO\hookrightarrow \End_F(A^\vee)$ by the assignment $a\mapsto a^\vee$. Let  $\lambda:A \ra A^\vee$ be a $\dag$-sesquilinear symmetric isogeny.  

For a finitely presented projective $\CO$-module $M$, the functor 
\[\shom{M}{A}(T):=\Hom_\CO(M,A(T))\]
on the fppf site over $F$ is representable by an abelian variety (cf. \cite[Proposition 7.33]{Milne-CMAV}). The action of $\CO$ on the abelian variety $\shom{M}{A}$ is given as $(af)(m)=f(a m)$. In particular, $\shom{\CO}{A}$ is identified with $A$ with compatible $\CO$-actions through the evaluation of a section at the identity. Endow $M^\vee=\Hom_\CO(M,\CO)$ with the $\CO$-action by $(as)(m)=s(am)$.  

Let $f:N\ra M$ be a homomorphism of finitely presented projective $\CO$-modules and  $\phi:A\ra B$ be an $\CO$-linear  homomorphism of abelian varieties. The induced homomorphism of abelian varieties 
\begin{equation}\label{morphism}
h(f,\phi): \shom{M}{A}\lra \shom{N}{B}, \quad s\mapsto \phi s f.
\end{equation}
is $\CO$-linear.
\begin{prop}\label{dual}
The dual abelian variety of $\shom{M}{A}$ is $\shom{M^\vee}{A^\vee}$ $($with compatible $\CO$-actions$)$ and the dual homomorphism $h(f,\phi)^\vee=h(f^\vee, \phi^\vee)$.
\end{prop}
\begin{proof}
Interpreting  in terms of Serre's tensor construction \cite{Serre94} or \cite[\S 7]{Conrad04}, one has 
\begin{equation}\label{serre-hom}
M^\vee\otimes_\CO A\cong \shom{M}{A},\quad s\otimes x\mapsto (m\mapsto  s(m)x)
\end{equation}
with compatible $\CO$-actions. Here the $\CO$-module structure of  the Serre tensor $M \otimes_\CO A$ is induced from that of $M$. Comparing with the normalization of $\CO$-module structures on $M^\vee$ and $A^\vee$ from \cite{Amir18}, the proposition follows from \cite[Proposition 5]{Amir18}.
\end{proof}
\begin{remark}
For any $\CO$-module $M$, denote ${^\dag M}$ the module $M$ with $\CO$-action twisted by the involution $\dag$. If both $f$ and $\phi$ are $\dag$-sesquilinear, the induced homomorphism $h(f,\phi)$ of (\ref{morphism}) is $\dag$-sesquilinear, which can also be viewed as an $\CO$-linear homomorphism
\[h(f,\phi):{^\dag \shom{M}{A}}=\shom{{^\dag M}}{{^\dag A}}\ra \shom{N}{B}.\]
\end{remark}

Let $K$ be the field of fractions of $\CO$. For $t\in K$, define an $\CO$-linear map $m_t$ and a $\dag$-sesquilinear map $m_t^\dag$ over $K$ by $m_t(x)=tx$ and $m_t^\dag(x)=tx^\dag$. Let $\fa\subset K$ be an invertible fractional ideal of $\CO$. Then $\fa$ is a projective finitely generated $\CO$-module. For $t\in \fa^{-1}$, $m_t|_{\fa}\in \fa^\vee$ and the assignment $t\mapsto m_t$ identifies $\fa^{-1}$ with $\fa^\vee$ as $\CO$-modules. In particular, $\CO$ is identified with $\CO^\vee$ in such a way.

\begin{coro}\label{sym-cri}
Let  $t\in K^\times$ be a nonzero element and let $\fa,\fb$ be invertible fractional ideals of $\CO$ satisfying $m^\dag_t(\fb)\subset \fa$.  The dual map of 
\[h(m^\dag_t, \lambda): \shom{\fa}{A}\lra \shom{\fb}{A^\vee}\]
is 
\[h(m^\dag_{t^\dag}, \lambda): \shom{\fb^{-1}}{A}\lra \shom{\fa^{-1}}{A^\vee}.\]
\end{coro}

In the following we assume $\CO$ is a one-dimensional noetherian intergral domain. Let $\fp$ be an invertible prime ideal of $\CO$.  Fix $\varpi\in \fp$ such that $\varpi$ generates $\fp/\fp^2$. Fix $t_0\in \fp^{-1}$ such that $t_0\varpi\equiv 1\mod \fp$ and necessarily $t_0$ generates $\fp^{-1}/\CO$. Let $i_\fp: \fp\ra \CO$ and $j_\fp:\CO\ra \fp^{-1}$ be the inclusion maps. Then $j_{\fp}$  is the dual map of $i_\fp$. 

The inclusion $i_\fp$ and the identity map $1_A$ induce an isogeny
\[h(i_\fp,1_A): A=\shom{\CO}{A}\lra \shom{\fp}{A}.\]
The kernel of this map is $\shom{\CO/\fp}{A}$. Associating to a section $\CO/\fp\ra A$ the image of the identity identifies $\shom{\CO/\fp}{A}$ with $A[\fp]$. We thus have an exact sequence of sheaves
\[0\lra A[\fp]\lra A\lra \shom{\fp}{A}\lra 0,\]
which induces the Kummer $\fp$-descent exact sequence
\[0\ra \Hom(\fp,A(F))\ra \RH^1(F,A[\fp])\ra \RH^1(F,A)[\fp]\ra 0.\]
The dual map of $h(i_\fp,1_A)$ is
\[h(j_{\fp}, 1_{A^\vee}): \shom{\fp^{-1}}{A^\vee}\lra \shom{\CO}{A^\vee}=A^\vee.\]
The kernel of $h(j_{\fp},1_{A^\vee})$ is $\shom{\fp^{-1}/\CO}{A^\vee}$ which is identified with $A^\vee[\fp]$ by  associating a section $\fp^{-1}/\CO\ra A^\vee$ the image of  $t_0\in \fp^{-1}/\CO$. To ease the notation, we denote 
\[g=h(i_\fp,1_A): A=\shom{\CO}{A}\ra B=\shom{\fp}{A}\]
and
\[g'=h(i_{\fp^\dag},1_A): A=\shom{\CO}{A}\ra B'=\shom{\fp^\dag}{A}.\]
We have the following commutative diagram
\begin{equation}\label{one-side}
\xymatrix{
0\ar[r]&A[\fp^\dag]\ar[r]\ar[d]^{\lambda}&\shom{\CO}{A}\ar[r]^{g'}\ar[d]^{h(m^\dag_{\varpi^\dag},\lambda)}&\shom{\fp^\dag}{A}\ar[r]\ar[d]^{h(m^\dag_{\varpi^\dag}, \lambda)}&0\\
0\ar[r]&A^\vee[\fp]\ar[r]&\shom{\fp^{-1}}{A^\vee}\ar[r]^{g^\vee}&\shom{\CO}{A^\vee}\ar[r]&0.
}
\end{equation}
Taking cohomologies, we have a commutative diagram of Kummer descent exact sequences
\begin{equation}\label{com-kummer}
\xymatrix{
0\ar[r]&B'(F)/g'(A(F))\ar[d]^{h(m^\dag_{\varpi^\dag},\lambda)}\ar[r]&\RH^1(F,A[\fp^\dag])\ar[d]^{\lambda}\ar[r]&\RH^1(F,A)[\fp^\dag]\ar[r]\ar[d]^{h(m^\dag_{\varpi^\dag},\lambda)}&0\\
0\ar[r]&A^\vee(F)/g^\vee(B^\vee(F))\ar[r]&\RH^1(F,A^\vee[\fp])\ar[r]&\RH^1(F,B^\vee)[\fp]\ar[r]&0.
}
\end{equation}

\begin{prop}\label{dual-kummer}
Suppose $\deg(\lambda)$ is invertible in  $\CO/\fp$. The vertical morphisms in the diagram $($\ref{com-kummer}$)$ are all isomorphisms.
\end{prop}
\begin{proof}
Since $\deg(\lambda)$ is invertible in $\CO/\fp$, $\lambda$ induces an isomorphism $A[\fp^\dag]\cong A^\vee[\fp]$. In particular, the middle vertical map is an isomorphism. It remains to prove the right vertical map is an isomorphism by snake lemma.  

There is an isogeny $\mu:A^\vee\ra A$ such that $\mu\lambda=\deg(\lambda)_{A}$ and $\lambda \mu=\deg(\lambda)_{A^\vee}$. Since $\lambda$ is $\dag$-sesquilinear, $\mu$ is also $\dag$-sesquilinear ( see \cite[Lemma 5.4, Proposition 5.12]{Edixhoven-AV}).  The diagram (\ref{one-side}) and the following diagram
\[\xymatrix{
\shom{\fp^{-1}}{A^\vee}\ar[r]^{g^\vee}\ar[d]^{h(m^\dag_{t_0},\mu)}&\shom{\CO}{A^\vee}\ar[d]^{h(m^\dag_{t_0}, \mu)}\\
\shom{\CO}{A}\ar[r]^{g'}&\shom{\fp^\dag}{A},
}\]
together with the canonical isomorphisms of (\ref{serre-hom})
\[{\fp^{\dag}}^{-1}\otimes A(F)\cong B'(F) \text{ and } \fp\otimes A^\vee(F)\cong B^\vee(F)\]
induce the following commutative diagrams
\begin{equation}\label{compose}
\xymatrix@C=6pt{
{\fp^{\dag}}^{-1}/\CO\otimes A (F)\ar@{}[r]|*=0[@]{\cong{\ }}\ar[d]^{m^\dag_{\varpi}\otimes \lambda}&B'(F)/g'(A(F))\ar[d]^{h(m^\dag_{\varpi^\dag},\lambda)}&
\CO/\fp\otimes A^\vee(F)\ar@{}[r]|*=0[@]{\cong{\ }}\ar[d]^{m^\dag_{t_0^\dag}\otimes \mu}&A^\vee(F)/g^\vee(B^\vee(F))\ar[d]^{h(m^\dag_{t_0},\mu)}\\
\CO/\fp\otimes A^\vee(F)\ar@{}[r]|*=0[@]{\cong{\ }} &A^\vee(F)/g^\vee(B^\vee(F)),&
{\fp^{\dag}}^{-1}/\CO\otimes A (F)\ar@{}[r]|*=0[@]{\cong{\ }} &B'(F)/g'(A(F)).
}
\end{equation}
By the choice of $t_0$ such that $\varpi t_0\equiv 1\mod \fp$, the compositions
\[m_{(\varpi t_0)^\dag}\otimes \deg(\lambda)_A=(m^\dag_{t_0^\dag}\otimes \mu)\circ(m^\dag_{\varpi}\otimes \lambda)\]
and 
\[m_{\varpi t_0}\otimes \deg(\lambda)_{A^\vee}=(m^\dag_{\varpi}\otimes \lambda)\circ (m^\dag_{t_0^\dag}\otimes \mu)\]
give the multiplication by $\deg(\lambda)$ on ${\fp^{\dag}}^{-1}/\CO\otimes A (F)$ and $\CO/\fp\otimes A^\vee(F)$ respectively. Let $p$ be the rational prime lying under $\fp$. Since both the groups are torsion of expoent $p$ and $p\nmid \deg(\lambda)$, both $m_{(\varpi t_0)^\dag}\otimes \deg(\lambda)_A$ and $m_{\varpi t_0}\otimes \deg(\lambda)_{A^\vee}$ are  isomorphisms, and thus all vertical maps in the (\ref{compose}) are isomorphisms. The proposition follows as desired.
\end{proof}

\section{A sheaf analogue of lifting adjoint pairings}
Let $L$ be a local field or a split quadratic extension of a local field over $\BQ_p$. Let $\CO$ be the ring of integers of $L$ endowed with an involution $\sigma$. If $L$ is a split quadratic extension, take $\sigma$ to be the automorphism switching the split components $\sigma(x,y)=(y,x)$. Let $\CC$ be a site and let $\mathfrak{Ab}_{\BZ_p}$ and $\Ab_\CO$ be the category of sheaves of $\BZ_p$-modules and $\CO$-modules over $\CC$ respectively. 
For $*=\BZ_p$ or $\CO$ and sheaves $M_1,M_2\in \Ab_*$, denote by $\bhom_*(M_1,M_2)$ the sheaf of homomorphisms $U\mapsto \Hom_*({M_1}|_{U},{M_2}|_{U})$ for any object $U\in \CC$. For $M_1,M_2\in \Ab_\CO$ and $N\in \Ab_{\BZ_p}$, the sheaf $\bhom_\CO(M_1,M_2)$ resp. $\bhom_{\BZ_p}(M_1, N)$ acquires an $\CO$-module structure: for any section $f\in \bhom_\CO(M_1,M_2)$ resp. $\bhom_{\BZ_p}(M_1,N)$ and $a\in \CO$, $(af)(m)=f(\sigma(a)m)$.

The ring $\CO$ is a free $\BZ_p$-module of finite rank. Let $d$ be a generator of the inverse different $\partial_{\CO/\BZ_p}^{-1}$. Then trace pairing 
\[\CO\times \CO\ra \BZ_p,\quad (x,y)\mapsto \Tr_{L/\BQ_p}(dxy)\]
is perfect. Denote $t_d:\CO\ra \BZ_p$ the trace map  $t_d(x)=\Tr_{L/\BQ_p}(dx)$. 

\begin{prop}\label{trace-isom}
Suppose $M_1\in \Ab_\CO$ and  $M_2\in \Ab_{\BZ_p}$. Composing $t_d\otimes_{\BZ_p} 1_{M_2}:\CO\otimes_{\BZ_p} M_2\ra M_2$ induces an isomorphism of sheaves of $\CO$-modules
\[\bhom_\CO(M_1,\CO\otimes_{\BZ_p}M_2)\xrightarrow{\cong} \bhom_{\BZ_p}(M_1,M_2).\]
\end{prop}
\begin{proof}
Taking $\CO$ in place of $\CR$ in \cite[Lemma A.3]{MR07}, the lemma follows verbatim. Then for any $U\in \CC$, 
\[\Hom_\CO(M_1(U),\CO\otimes_{\BZ_p}M_2(U))\xrightarrow{\cong} \Hom_{\BZ_p}(M_1(U),M_2(U))\]
as $\CO$-modules. Since the functoriality in $U$ is clear, the proposition follows.
\end{proof}

\begin{defn}
Let $M,M'\in \Ab_\CO$ and $N\in \Ab_{\BZ_p}$. 
\begin{itemize}
\item[(1)] A pairing $s:M\times M\ra M'$ is called $\sigma$-sesquilinear if, for any object $U\in \CC$, the pairing $s(U):M(U)\times M(U)\ra M'(U)$ satisfies, for any $a\in \CO$ and $x,y\in M(U)$,
\[s(U)(ax,y)=s(U)(x,\sigma(a)y)=as(U)(x,y).\] 
A $\sigma$-sesquilinear pairing $s$ is nondegenerate or perfect, if the induced map 
\[M\ra\Hom_\CO(M,M')\]
is an isomorphism of sheaves. 
\item[(2)] A pairing $s:M\times M\ra N$ is called $\sigma$-adjoint if, for any object $U\in \CC$, $s(U)$ is $\BZ_p$-bilinear and for any $a\in \CO$ and $x,y\in M(U)$,
\[s(U)(ax,y)=s(U)(x,\sigma(a)y).\] 
\item[(3)] The transpose $s^*$ of a pairing $s:M\times M\ra \CO\otimes_{\BZ_p}N$ is defined by
\[s^*(U)(x,y)=(\sigma\otimes 1_N)(s(U)(y,x)).\] 
For $\lambda=\pm 1$, a $\sigma$-sesquilinear pairing $s$ is called $\lambda$-hermitian if $s^*=\lambda s$.
\end{itemize}
 \end{defn}
The above definition applies to the case $\CO=\BZ_p$ with $\sigma$ the trivial involution. If $\sigma=1$, then $d^{\sigma-1}:=\sigma(d)/d=1$; Otherwise, as in (\ref{gen-dif}), we fix $d\in \partial_{\CO/\BZ_p}^{-1}$ so that $d^{\sigma-1}=+1$ resp. $-1$ if $L/L_0$ is ramified resp. unramified where $L_0$ is the fixed field of $\sigma$.

 \begin{coro}\label{sesquilinear-adjoint}
Suppose $M\in \Ab_{\CO}$ and $N\in \Ab_{\BZ_p}$. Composing $t_d\otimes_{\BZ_p} 1_N$ induces a bijection of the set of $\sigma$-sequilinear pairings $s:M\times M\ra \CO\otimes_{\BZ_p}N$ and the set of $\sigma$-adjoint pairings $s':M\times M\ra N$. This bijection preserves nondegeneracy and takes exactly $\lambda$-hermitian pairings to $(d^{\sigma-1}\lambda)$-hermitian parings. 
\end{coro}
\begin{proof}
It follows from Proposition \ref{trace-isom} that
\[\Hom_\CO(M,\bhom_\CO(M,\CO\otimes_{\BZ_p}N))\xrightarrow{\cong}\Hom_{\CO}(M,\bhom_{\BZ_p}(M,N)).\]
The LHS is the set of $\sigma$-sesquilinear pairings $s:M\times M\ra \CO\otimes_{\BZ_p}N$ while the RHS is the set of $\sigma$-adjoint pairings $s':M\times M\ra N$. By Proposition \ref{trace-isom}, this bijection preserves nondegeneracy. If this bijection sends $s$ to $s'$, it also sends ${s}^*$ to $d^{\sigma-1}{s'}^*$. Thus $s+\lambda s^*=0$ if and only if $s'+\lambda d^{\sigma-1}{s'}^*=0$.
\end{proof}

\begin{prop}\label{orthogonal}
Suppose $M\in \Ab_{\CO}$ and $N=\BZ_p$ or $\BQ_p/\BZ_p$.  Let $s:M\times M \ra \CO_{\BZ_p}\otimes N$ be a nondegenerate $\dag$-sesquilinear pairing and let $s'$ be the associated $\sigma$-adjoint pairing. Let $W$ be an $\CO$-submodule of $M$. The the orthogonal components of $W$ under $s$ and $s'$ coincide. 
\end{prop}
\begin{proof}
Denote the orthogonal complements by $W^{\perp,s}$ and $W^{\perp,s'}$ respectively.  Since $s'(x,y)=(t_d\otimes 1_N)(s(x,y))$, $W^{\perp,s}\subset W^{\perp,s'}$. For the converse, suppose $x\in W^{\perp, s'}$ and $y\in W$. Since $N=\BZ_p$ or $\BQ_p/\BZ_p$, we can write $s(x,y)=b\otimes n$ for some $b\in \CO$ and $0\neq n\in N$. Let $I$ be the annihilator of $n$ in $\BZ_p$. If $N=\BZ_p$, then $I=0$; Otherwise $I=p^m\BZ_p$ for some $m>0$. Since $W$ is an $\CO$-module, $ay\in W$ for all $a\in \CO$. Then for all $a\in \CO$,
\[s'(x, ay)=(t_d\otimes 1_N)(da^\dag s(x,y))=\Tr_{L/\BQ_p}(da^\dag b) n=0.\]
 Thus $\Tr_{L/\BQ_p}(\partial^{-1}_{\CO/\BZ_p}b)\subset I$. If $I=0$, then $b=0$. If $I=p^m\BZ_p$, then $\Tr_{L/\BQ_p}(\partial^{-1}_{\CO/\BZ_p}p^{-m}b)\subset \BZ_p$. Since $\partial^{-1}_{\CO/\BZ_p}$ is the maximal fractional ideal with traces in $\BZ_p$, $\partial^{-1}_{\CO/\BZ_p}p^{-m}b\subset \partial^{-1}_{\CO/\BZ_p}$ and then $b\in p^m\CO$. Thus in any case $s(x,y)=b\otimes n=0$ as deisred. 
 \end{proof}

\bibliographystyle{plain}
\bibliography{reference}
\end{document}